\documentstyle[epsf,stmaryrd,amscd,amstex,amssymb,
12pt
]{amsart}
\input xy

\newcommand{\nc}{\newcommand}
\newcommand{\rc}{\renewcommand}

\nc{\Aut}{{	\operatorname{Aut}	}}
\nc{\codim}{{	\operatorname{codim}	}}
\nc{\Ob}{{	\operatorname{Ob}	}}
\nc{\PGL}{{	\operatorname{PGL}	}}
\nc{\supp}{{	\operatorname{supp}	}}
\nc{\tr}{{	\operatorname{tr}	}}

\nc{\Rep}{{	{\cal{R}}ep		}}
\nc{\one}{{	\mbox{\bf{1}}		}}
\nc{\iso}{	\overset{\sim}{\lra}	}

\nc{\nen}{\newenvironment}

\nc{\pr}{\protect}

\nc{\nn}{{\newline}}
\nc{\np}{{\newpage}}	

\nc{\lab}{	\label}
\nc{\npp}{{	\newpage\setcounter{page}{0}	}}
\nc{\setpa}{		\setcounter{part}		}
\nc{\setse}{		\setcounter{section}	}
\nc{\setsus}{		\setcounter{subsection}		}
\nc{\setsss}{		\setcounter{subsubsection}	}
\nc{\setpage}{		\setcounter{page}	}

\nc{\nfd}{ $$\text{ This version is preliminary and approximate,
		             it is not for distribution. }$$	}

\nc{\noi}{{\noindent}}
\rc{\pf}{{	\noindent {\em Proof.}		}}
					
\nc{\epf}{\square}

\nc{\s}{{\spadesuit}}
\nc{\heart}{{\tiny \cen{\tiny $\heartsuit $ }	}} 
\nc{\cont}{\tableofcontents}

\nc{\sbr}{{	\smallpagebreak	}}
\nc{\mbr}{{	\medpagebreak	}}
\nc{\bbr}{{	\bigpagebreak	}}
\nc{\bbb}{ 	\boldsymbol 	}

\nc{\bib}{		}
\nc{\bit}[1]{	\bibitem[#1]{#1} 		}

\rc{\b}{ 	\big         			}
\nc{\lam}[1]{{ 	\text{\large $#1$	}	}}
\nc{\smm}[1]{{ 	\text{\small $#1$	}	}}
\nc{\fom}[1]{{ 	\text{\footnotesize $#1$	}	}}
\nc{\tinm}[1]{{ \text{\tiny $#1$	}	}}

\nc{\bu}{ \bullet         }  			
\nc{\bbu}{ \aa{\bbb \bullet}         }  	
\nc{\bus}{{	^\bullet	}}	 	
\nc{\bui}{{	_\bullet	}}	 	

\nc{\bem}{{	\begin{em}	}}
\nc{\eem}{{	\end{em} 	}}

\nc{\bbox}{{	\blackbox	}}	
\nc{\bx}{	\boxed	}		

\nc{\mmbox}[1]{{	\mbox{$#1$}	}}	
\nc{\tbox}[1]{{		\mbox{\tx{#1}}	}}
 	
\nc{\ot}{		\leftarrow			}
\nc{\tto}{		\longrightarrow			}
\nc{\ott}{		\longleftarrow			}
\nc{\too}[1]{{		\aa{#1}\rightarrow			}}
\nc{\oot}[1]{{		\aa{#1}\leftarrow			}}
\nc{\ttoo}[1]{{		\aa{#1}\longrightarrow			}}
\nc{\oott}[1]{{		\aa{#1}\longleftarrow			}}

\nc{\Too}[2]{{		\aa{#1}{\bb{#2}\rightarrow}		}}
\nc{\ooT}[2]{{		\aa{#1}{\bbb{#2}\leftarrow}		}}
\nc{\TToo}[2]{{		\aa{#1}{\bb{#2}\longrightarrow}		}}
\nc{\ooTT}[2]{{		\aa{#1}{\bbb{#2}\longleftarrow}		}}
\nc{\toot}[2]{{		\aa{#1}{\bb{#2}\rightleftarrows}	}}
\nc{\ttoot}[2]{{	\aa{#1}{\bb{#2}\rightleftrightarrows}	}}

\nc{\ra}{{	\rightarrow		}}
\nc{\laa}{{	\leftarrow	}}	
\nc{\lra}{{\longrightarrow}}

\nc{\lr}{{\leftrightarrow}}     	
\nc{\lrs}{{\rightleftarrows}}     	

\nc{\imp}{{\Rightarrow}}        	
\nc{\impp}{{\Leftarrow}}        	
\nc{\eq}{{\Leftrightarrow}}        	
\nc{\impl}{{\Longrightarrow}}        	
\nc{\imppl}{{\Longleftarrow}}        	
\nc{\eql}{{\Longleftrightarrow}}        	
	\nc{\Ra}{{\Rightarrow}}         	
	\nc{\LRa}{{\Leftrightarrow}}        	

\nc{\inj}{{\pr	\hookrightarrow	}}    		
\nc{\injj}{{\pr	\hookleftarrow	}}    		

\nc{\sur}{{	\twoheadrightarrow	}}	
\nc{\surr}{{	\twoheadleftarrow	}}	
			
\nc{\mm}{{	\mapsto		}}     		
\nc{\mmm}{{	\leftarrow\shortmid }}		
\nc{\ainj}[1]{{\aa{#1}{\pr\hookrightarrow}	}}    	
\nc{\ainjj}[1]{{\aa{#1}{\pr\hookleftarrow}	}}    	

\nc{\asur}[1]{{	\aa{#1}\twoheadrightarrow	}}	
\nc{\asurr}[1]{{\aa{#1}\twoheadleftarrow	}}	
			
\nc{\amm}[1]{{	\aa{#1}\mapsto		}}     	
\nc{\ammm}[1]{{	\aa{#1}\leftarrow\shortmid }}	

\nc{\va}{{\uparrow}}              		

\nc{\syp}[1]{	^{ (#1) }		} 	
\nc{\up}[1]{	^{ (#1) }		} 	
\nc{\lp}[1]{	_{ (#1) }		}	
\nc{\hp}[1]{	^{ [#1] }		}	
\nc{\cle}{\preceq}		
\nc{\cl}{\prec}			
\nc{\cge}{\succeq}		
\nc{\cg}{\succ}			

\nc{\bb}{	\pr\underset 	}           
\rc{\aa}{ 	\pr\overset 	}            

\nc{\indd}{{ ${} \ \ \ \ \  \ \        {} $	}}	
\nc{\inddd}{{ 	\indd\indd			}}	
\nc{\nnd}{{ 	\nn  \indd 			}}	
\nc{\nndb}{{ 	\nn  \indd $\bullet$		}}	
		
\nc{\bce}{	\begin{center}	}
\nc{\ece}{	\end  {center}	}
\nc{\cen}[1]{	\begin{center}	{\em  #1}	\end  {center}	}

\nc{\bss}{{\backslash}}           		
\nc{\barr}{ 	\overline 	}      		
\nc{\ud}{	\underline	}		

\nc{\ti}{\tilde}              
\nc{\tii}{\widetilde}         

\nc{\hatt}{\widehat}				
\nc{\hata}{{	\bbb{ \hat{} }		}}	

\nc{\ch}{\check}              			
\nc{\cha}{{ 	\bbb{ \check{} }	}}      

\nc{\sub}{{	\subseteq	}}         
\nc{\subb}{{	\supseteq	}}         
\nc{\nsub}{{	\nsubseteq	}}         
\nc{\nsubb}{{	\nsupseteq	}}         %

\nc{\nin}{{	\notin	}}

\nc{\lb}{\langle}             				
\nc{\rb}{\rangle}

\nc{\lB}{	\left(	}             			
\nc{\rB}{	\right)	}
\nc{\BBl}{{	\bbb{ \left( \right.}	}}             	
\nc{\BBr}{{	\bbb{ \left. \right)}	}}

\nc{\Pa}[2]{ {\lb} #1 {,} #2 {\rb} }				

\nc{\cD}[1]{ \tx{ $$\CD {#1} \endCD $$ }  }		

\nc{\mat} {		\left(		\matrix	}	
\nc{\emat}{		\endmatrix	\right)	}
\nc{\sm} {		\left(		\smallmatrix	}	
\nc{\esm}{		\endsmallmatrix	\right)	}
\nc{\smat} {		\left(		\smallmatrix	}	
\nc{\esmat}{		\endsmallmatrix	\right)	}

\nc{\matr} {		\left[		\matrix	}	
\nc{\ematr}{		\endmatrix	\right]	}
\nc{\smr} {		\left[		\smallmatrix	}	
\nc{\esmr}{		\endsmallmatrix	\right]	}
\nc{\smatr} {		\left[		\smallmatrix	}	
\nc{\esmatr}{		\endsmallmatrix	\right]	}

\nc{\imat} {		\left.		\matrix	}	
\nc{\eimat}{		\endmatrix	\right.	}
\nc{\ism} {		\left.		\smallmatrix	}	
\nc{\eism}{		\endsmallmatrix	\right.	}

\nc{\ca}{		\left\{		\smallmatrix	}	
\nc{\eca}{		\endsmallmatrix	\right\}	}
\nc{\Ca}{		\left\{		\matrix		}	
\nc{\Eca}{		\endmatrix	\right.		}	
\nc{\eCa}{		\endmatrix	\right\}	}	

\nc{\com}{	\begin{diagram}	}
\nc{\ecom}{	  \end{diagram}	}

\nc{\tab}{	\begin{tabular}		}
\nc{\etab}{	\end{tabular}		}	
\nc{\hl}{{	\hline			}}

\nc{\Eq}{	\begin{equation}	}
\nc{\Eeq}{	\end{equation}	}
\nc{\aln}{	\begin{align}	}
\nc{\ealn}{	\end{align}	}

\nc{\Rpart}{	\rc{\thepart}{\Roman{part}}	}
\nc{\Apart}{	\rc{\thepart}{\arabic{part}}	}
		
\nc{\rref}[2]{\ref{#1}.\ref{#2}}

\nc{\pa}[1]{ 	\part{#1}		}

\nc{\se}[1]{ 	\section{\bf#1}		}
\nc{\ses}[1]{ 	\section*{\bf#1}		}
\nc{\sus}{ 	\subsection		}
\nc{\sss}{ 	\subsubsection		}

\nc{\Lem}{ 	\subsection{Lemma}		}
\nc{\lem}{ 	\subsubsection{Lemma}		}
\nc{\slem}{ 	\subsubsection*{Lemma}		}
\nc{\sublem}{ 	\subsubsection{ Sublemma}	}
\nc{\ssublem}{ \subsubsection*{ Sublemma}	}

\nc{\Lemm}{ 	\subsection{Lemma}		}
\nc{\lemm}{ 	\subsubsection{Lemma}		}
\nc{\slemm}{ 	\subsubsection*{Lemma}		}
\nc{\sublemm}{ 	\subsubsection{ Sublemma}	}
\nc{\ssublemm}{ \subsubsection*{ Sublemma}	}

\nc{\Pro}{ 	\subsection{Proposition}	}
\nc{\pro}{ 	\subsubsection{Proposition}	}
\nc{\spro}{ 	\subsubsection*{Proposition}	}

\nc{\Cor}{ 	\subsection{Corollary}		}
\nc{\cor}{ 	\subsubsection{Corollary}	}
\nc{\scor}{ 	\subsubsection*{Corollary}	}

\nc{\Corr}{ 	\subsection{Corollary}		}
\nc{\corr}{ 	\subsubsection{Corollary}	}
\nc{\scorr}{ 	\subsubsection*{Corollary}	}

\nc{\Theo}{ 	\subsection{Theorem}		}		
\nc{\theo}{ 	\subsubsection{Theorem}		}
\nc{\stheo}{ 	\subsubsection*{Theorem}	}

\nc{\pretheo}{ 	\subsubsection{Pretheorem}	}

\nc{\rem}{ 	\subsubsection{Remark}		}
\nc{\srem}{ 	\subsubsection*{Remark}	}

\nc{\rems}{ 	\subsubsection{Remarks}		}

\nc{\srems}{ 	\subsubsection*{Remarks}	}

\nc{\Def}{ 	\subsection{Definition}		}
\nc{\ddef}{ 	\subsubsection{Definition}	}
\nc{\sdef}{ 	\subsubsection*{Definition}	}

\nc{\comm}{ 	\subsubsection{Comment}		}
\nc{\scomm}{ 	\subsubsection*{Comment}	}
\nc{\comms}{ 	\subsubsection{Comments}		}
\nc{\scomms}{ 	\subsubsection*{Comments}	}

\nc{\claim}{ 	\subsubsection{Claim}		}
\nc{\sclaim}{ 	\subsubsection*{Claim}	}

\nc{\nota}{ 	\subsubsection{Notation}	}

\nc{\conj}{ 	\subsubsection{Conjecture}	}
\nc{\sconj}{ 	\subsubsection*{Conjecture}	}

\nc{\ex}{ 	\subsubsection{Example}		}
\nc{\sex}{ 	\subsubsection*{Example}	}
\nc{\exs}{ 	\subsubsection{Examples}	}
\nc{\sexs}{ 	\subsubsection*{Examples}	}
\nc{\Ex}{ 	\subsection{Example}		}
\nc{\sEx}{ 	\subsection*{Example}	}
\nc{\Exs}{ 	\subsection{Examples}	}
\nc{\sExs}{ 	\subsection*{Examples}	}

\nc{\que}{ 	\subsubsection{Question}	}
\nc{\ques}{ 	\subsubsection{Questions}	}
\nc{\sque}{ 	\subsubsection*{Question}	}
\nc{\sques}{ 	\subsubsection*{Questions}	}

\nc{\bi}{	\begin{itemize}\item		}
\rc{\i}{	\item			}
\nc{\ei}{ \end{itemize}	}
\nc{\ben}{	\begin{enumerate}\item		}
\nc{\een}{	\end{enumerate}			}
\nc{\ftt}[1]{{\footnote{#1}}}
\nc{\fttt}[1]{{$^($\footnote{#1}$^)$}}
\nc{\bftt}[1]{\footnote{#1}}
\nc{\bububu}{{\bb{\bu} { \aa{\bu}\bu }}}
\nc{\ftx}[1]{{${\bububu}{}^($\footnote{#1}$^)\bububu $}}

\nc{\adc}[1]{{\addtocounter{#1}{1}}}

\newcounter{roma}
\nc{\r}{ \fbox { \fbox { R$\theroma$.\adc{roma} }}}
\nc{\rrr}[1]{\fbox{\fbox{READ:\ #1}R$\theroma$.\adc{roma}}}

\newcounter{suggestions}

\nc{\fff}[1]{ \fbox{ \fbox{#1}S$\thesuggestions$.\adc{suggestions}} }
\nc{\f}[1]{\ftx{ \fbox{S$\thesuggestions.$}\adc{suggestions} \newline  }}
\nc{\sugg}{\tx{\bf SUGGESTION}}
\nc{\SU}[1]{\fff{\sugg$@>>>$}}
	\nc{\dv}{\tx{\bf A dying version}}
\nc{\DV}[1]{\fff{{\dv$@>>>$}}}
	\nc{\cracra}{\tx{\bf{END}}}
\nc{\ENDD}{\begin{flushright}{\fff{$@<<<$\cracra}}\end{flushright}}
\nc{\END}{{\fff{$@<<<$\cracra}}}
\nc{\ENDDV}{\begin{flushright}{\fff{$@<<<$\tx{\bf END of \dv}}}\end{flushright}}

\nc{\IN}[1]{\fbox{IN:\ $@>>>$\ \fbox{#1}}}
\nc{\OUT}[1]{\fbox{\fbox{#1}\ {$@>>>$\ OUT}}}

\nc{\Ao}{{	\A^1	}}
\nc{\Po}{{	\pP^1	}}

\nc{\h}{{	\hslash	}}	

\nc{\All}{{	\forall		}}
\nc{\Exx}{{	\exists 	}}
\nc{\yy}{\infty}
\nc{\ys}{{  \frac{\infty}{2}  }}

\nc{\ii}{{i\in I}}
\nc{\ww}{{w\in W}}
\nc{\SES}[5]{{	0 @>>> {#1} @>{#2}>> {#3} @>{#4}>> {#5} @>>> 0	}}
\nc{\Ses}[3] {{	0 @>>> {#1} @>>>     {#2} @>>>     {#3} @>>> 0	}}
\nc{\pl}{{\oplus}}              		
\nc{\tim}{{\times}}
\nc{\timL}{{\aa{L}\times}}
\nc{\timR}{{\aa{R}\times}}
\nc{\Rtim}{{\aa{R}\times}}
\nc{\btim}{{\boxtimes}}
\nc{\ltim}{\ltimes}                  	%
\nc{\rtim}{\rtimes}			%
\nc{\ltr}{\triangleleft}        %
\nc{\rtr}{\triangleright}       %

\nc{\ten}{{	\otimes		}}
\nc{\Lten}{{	\aa{L}\otimes	}}            
\nc{\Ltim}{{	\aa{L}\times	}}            
\nc{\Lcap}{{	\aa{L}\cap	}}            
\nc{\tenA}{	\bb{A}\ten	}
\nc{\tenB}{	\bb{B}\ten	}
\nc{\tenZ}{	\bb{\Z}\ten	}
\nc{\tenR}{	\bb{\R}\ten	}
\nc{\tenC}{	\bb{\C}\ten	}
\nc{\tenk}{	\bb{\k}\ten	}
\nc{\bten}{{\boxtimes}}         		

\nc{\con}{{ @>>{\protect\cong}> }}  	
\nc{\conl}{{ 	@>>{\cong}>	}}  	
\nc{\conn}{{    @<{\cong}<<  	}}  	
\nc{\Con}{{	\equiv		}}	
\nc{\appr}{{	\sim		}}	
\nc{\eqr}{{	\sim		}}	
\nc{\equi}{{	\sim		}}	

\nc{\fra}{ 	\frac	}     	
\nc{\ffr}[2]{{ 	\text{\footnotesize $\frac{#1}{#2}$	}	}}

\nc{\ha}{{ \frac{1}{2} }}     		
	\nc{\half}{{ \frac{1}{2} }}    	

\nc{\ci}{{\circ}}               
\nc{\cd }{{\cdot}}            	
\nc{\cddd}{{\cdot\cdot\cdot}}	

\nc{\ox}{{	\OO_X		}}               
\nc{\omx}{{	\om_X		}}               
\nc{\Omx}{{	\Om_X^1		}}               
\nc{\Coh}{{	\CC oh		}}               %
\nc{\qcoh}{{	q\Coh		}}               %

\nc{\xt}{{	X_*(T)		}}
\nc{\Xt}{{	X^*(T)		}}

\nc{\cfm}{{	co\fm		}}	

\nc{\cupp}{\bigcup}             
\nc{\capp}{\bigcap}
\nc{\pll}{\bigoplus}

\nc{\pii}{\prod}                
\nc{\ppii}{\bigprod}            
\nc{\cci}{\sqcup}              
\nc{\ccii}{\bigsqcup}

\nc{\wwe}{\bigwedge}            
\nc{\cce}{\bigcoprod}           

\nc{\aaa}{	\stackerel	}	

\nc{\edd}{{ \end{document}	}}

\nc{\tx}{	\text		}		
\nc{\df}{{ \protect\overset{ \text{def}}= 	}}		
\nc{\dff}{{ \ \df\				}}		
\nc{\inv}{{ {}^{-1}      }}			

\nc{\thh}{	^{\text{th}}	}                     	
\nc{\st}{	^{\text{st}}	}                     	
\nc{\nd}{	^{\text{nd}}	}                     	
\nc{\rd}{	^{\text{rd}}	}                     	

\nc{\pmo}{{ 	\pm 1		}}
\nc{\mpo}{{ 	\mp 1		}}

\nc{\htt}{  \text{ht}}				

\nc{\emp}{{   \emptyset}}      			

\nc{\cowe}{{	\vee	}}			
\nc{\we}{{\wedge}}				
\nc{\wee}{{	\aa{\bullet}\wedge	}}		
\nc{\wetwo}{{     \pr\overset{2}\wedge       }}	

\nc{\limp}{{	\pr\underset {\leftarrow} \lim		}}	
\nc{\Limp}{{	\pr\underset {\leftarrow} {\bbb\lim}	}}	
\nc{\limi}{{	\pr\underset {\rightarrow}\lim		}}      
\nc{\Limi}{{	\pr\underset {\rightarrow}{\bbb\lim}	}}	
\nc{\llim}[1]{	 \bb{#1}\lim        	}   
\nc{\llimp}[1]{ \bb{#1}{ \pr\underset {\leftarrow} \lim       } }
\nc{\LLimp}[1]{ \bb{#1}{ \pr\underset {\leftarrow} {\bbb\lim} } } 	
\nc{\llimi}[1]{ \bb{#1}{ \pr\underset {\rightarrow}\lim       } }
\nc{\LLimi}[1]{ \bb{#1}{ \pr\underset {\rightarrow}{\bbb\lim} } }	
						
\nc{\ppp}{{ {\Bbb P}^1 }}            		
\nc{\ppn}{{ {\Bbb P}^n }}            		
\nc{\pt}{	{ \text{pt} }	}		
				
\nc{\qlb}{{ \barr{{\Bbb Q}_l} }}      		
\nc{\ffq}{{  {\Bbb F}_q  }}           		
\nc{\ffp}{{  {\Bbb F}_p  }}           		
\nc{\tw}{   {}^{(1)}	}		

\nc{\Ab}{{ 	\AA b 		}}      		%
\nc{\Set}{{ 	\SS et 		}}      		%
\nc{\Top}{{ 	\TT op 		}}      		%
\nc{\Pic}{{ 	\tx{Pic}	}}      		%

\nc{\del}{{\partial }}
\nc{\delb}{{\partial }}
\nc{\dd}[2]{	\fra{d{#1}}{d{#2}}		}
\nc{\ddel}[2]{	\fra{\del{#1}}{\del {#2}}	}

\nc{\Spec}{{ 	\text{Spec}      		}}
\nc{\Specf}{{ 	\text{Specf}      		}}
\nc{\Spf}{{ 	\text{Spf}      		}}
\nc{\hk}{{     \text{hyperk$\ddot{a}$hler} 	}}
\nc{\susy}{{\text{supersymmetry}}}

\nc{\ie}{{,\ \     \text{i.e.,}\ \ 	}}
\nc{\iif}{{\ \     \text{if}\ \ 	}}
\nc{\aand}{{\ \ \  \text{and}\ \ \ 	}}
\nc{\hence}{{\ \ \ \text{hence}\ \ \ 	}}
\nc{\while}{{\ \ \ \text{while}\ \ \ 	}}
\nc{\with}{{\ \ \  \text{with}\ \ \ 	}}
\nc{\oor}{{\ \     \text{or}\ \ 	}}
\nc{\foor}{{\ \     \text{for}\ \ 	}}
\nc{\suchthat}{{\ \     \text{such that}\ \ 	}}

\nc{\rk}{{\operatorname{rk}}}
\nc{\Ker}{{\operatorname{Ker}}}
\nc{\Coker}{{\operatorname{Coker}}}
\rc{\Im}{{ 	\text{Im} 	}}
\nc{\rank}{{	\ \text{rank} 	}}
\nc{\Res}{{	\  \text{Res}   }}
\nc{\Hom}{{\operatorname{Hom}}}
\nc{\End}{{	\text{End}	}}
\nc{\RHom}{{	\text{RHom}	}}
\nc{\HHom}{{	\text{$\HH$om}	}}
\nc{\RHHom}{{	\text{R$\HH$om} }}
\nc{\RGa}{{	\text{R$\Ga$}	}}
\nc{\EEnd}{{	\text{$\EE nd$}	}}
\nc{\AAut}{{	\text{$\AA ut$}	}}
\nc{\Ext}{{\operatorname{Ext}}}
\nc{\Tor}{{\operatorname{Tor}}}
\nc{\Der}{{	\text{Der}	}}

\nc{\ord	}{{ \text{ord} }}			
\nc{\divv	}{{ \text{div} }}			
\nc{\Lie	}{{ \text{Lie} }}

\nc{\timA} {{   \pr\underset{A}\tim             }}
\nc{\timB} {{   \pr\underset{B}\tim             }}
\nc{\timC} {{   \pr\underset{C}\tim             }}
\nc{\timG} {{   \pr\underset{G}\tim             }}
\nc{\timH} {{   \pr\underset{H}\tim             }}
\nc{\timN} {{   \pr\underset{N}\tim             }}
\nc{\timP}{{    \pr\underset{P}\tim             }}
\nc{\timQ}{{    \pr\underset{Q}\tim             }}
\nc{\timS} {{   \pr\underset{S}\tim             }}
\nc{\timT} {{   \pr\underset{T}\tim             }}
\nc{\timU} {{   \pr\underset{U}\tim             }}
\nc{\timV} {{   \pr\underset{V}\tim             }}
\nc{\timX} {{   \pr\underset{X}\tim             }}
\nc{\timY} {{   \pr\underset{Y}\tim             }}
\nc{\timZ} {{   \pr\underset{Z}\tim             }}

\nc{\ab}{{       ^{\text{ab}}   		}}
\nc{\af}{{       ^{\text{aff}}  		}}

\nc{\cod}{\text{codim}}	

\rc{\AA}{{\cal A}}
\nc{\BB}{{\cal B}}
\nc{\CC}{{\cal C}}
\nc{\DD}{{\cal D}}
\nc{\EE}{{\cal E}}
\nc{\FF}{{\cal F}}
\nc{\GG}{{\cal G}}
\nc{\HH}{{\cal H}}
\nc{\II}{{\cal I}}
\nc{\JJ}{{\cal J}}
\nc{\KK}{{\cal K}}
\nc{\LL}{{\cal L}}
\nc{\MM}{{\cal M}}
\nc{\NN}{{\cal N}}
\nc{\OO}{{\cal O}}
\nc{\PP}{{\cal P}}
\nc{\QQ}{{\cal Q}}
\nc{\RR}{{\cal R}}
\rc{\SS}{{\cal S}}
\nc{\TT}{{\cal T}}
\nc{\UU}{{\cal U}}
\nc{\VV}{{\cal V}}
\nc{\WW}{{\cal W}}
\nc{\ZZ}{{\cal Z}}
\nc{\XX}{{\cal X}}
\nc{\YY}{{\cal Y}}

\nc{\A}{{\Bbb A }}
\nc{\B}{{\Bbb B}}
\nc{\C}{{\Bbb C}}
		\nc{\cc}{{\Bbb C}}
\nc{\Cs}{{\Bbb C^*}}
		\nc{\cs}{{\Bbb C^*}}
		\nc{\ccs}{{\Bbb C^*}}
\nc{\D}{{\Bbb D}}
\nc{\E}{{\Bbb E}}
\nc{\F}{{\Bbb F}}
\nc{\G}{{\Bbb G}}
	\nc{\hH}{{\Bbb H}}
\nc{\I}{{\Bbb I}}
\nc{\J}{{\Bbb J}}
\nc{\K}{{\Bbb K}}
	\nc{\lL}{{\Bbb L}}
\nc{\M}{{\Bbb M}}
\nc{\N}{{\Bbb N}}
	\nc{\oO}{{\Bbb O}}
	\nc{\pP}{{\Bbb P}}
\nc{\Q}{{\Bbb Q}}
\nc{\R}{{\Bbb R}}
	\nc{\sS}{{\Bbb S}}
\nc{\T}{{\Bbb T}}
\nc{\U}{{\Bbb U}}
\nc{\V}{{\Bbb V}}
\nc{\W}{{\Bbb W}}
\nc{\Z}{{\Bbb Z}}
\nc{\X}{{\Bbb X}}
\nc{\Y}{{\Bbb Y}}

\nc{\k}{\Bbbk}

\let\L\lL
\let\O\oO
\let\P\pP
\let\S\sS

\nc{\fA}{{\frak A}}
\nc{\fB}{{\frak B}}
\nc{\fC}{{\frak C}}
\nc{\fD}{{\frak D}}
\nc{\fE}{{\frak E}}
\nc{\fF}{{\frak F}}
\nc{\fG}{{\frak G}}
\nc{\fH}{{\frak H}}
\nc{\fI}{{\frak I}}
\nc{\fJ}{{\frak J}}
\nc{\fK}{{\frak K}}
\nc{\fL}{{\frak L}}
\nc{\fM}{{\frak M}}
\nc{\fN}{{\frak N}}
\nc{\fO}{{\frak O}}
\nc{\fP}{{\frak P}}
\nc{\fQ}{{\frak Q}}
\nc{\fR}{{\frak R}}
\nc{\fS}{{\frak S}}
\nc{\fT}{{\frak T}}
\nc{\fU}{{\frak U}}
\nc{\fV}{{\frak V}}
\nc{\fW}{{\frak W}}
\nc{\fZ}{{\frak Z}}
\nc{\fX}{{\frak X}}
\nc{\fY}{{\frak Y}}
\nc{\fa}{{\frak a}}
\nc{\fb}{{\frak b}}
\nc{\fc}{{\frak c}}
\nc{\fd}{{\frak d}}
\nc{\fe}{{\frak e}}
\nc{\fg}{{\frak g}}
\nc{\fh}{{\frak h}}
\nc{\fiI}{{\frak i}}  
	\nc{\ffi}{{\frak i}}  
\nc{\fj}{{\frak j}}
\nc{\fk}{{\frak k}}
\nc{\fl}{{\frak{l}}}
\nc{\fm}{{\frak m}}
\nc{\fn}{{\frak n}}
\nc{\fo}{{\frak o}}
\nc{\fp}{{\frak p}}
\nc{\fq}{{\frak q}}
\nc{\fr}{{\frak r}}
\nc{\fs}{{\frak s}}
\nc{\ft}{{\frak t}}
\nc{\fu}{{\frak u}}
\nc{\fv}{{\frak v}}
\nc{\fw}{{\frak w}}
\nc{\fz}{{\frak z}}
\nc{\fx}{{\frak x}}
\nc{\fy}{{\frak y}}

\nc{\al}{{\alpha }}
\nc{\be}{{\beta }}
\nc{\ga}{{\gamma }}
\nc{\de}{{\delta }}
\nc{\ep}{{\varepsilon }}
\nc{\vap}{{\epsilon }}

\nc{\ze}{{\zeta }}
\nc{\et}{{\eta }}
\rc{\th}{{\theta }}
\nc{\vth}{{\vartheta }}

\nc{\io}{{\iota }}
\nc{\ka}{{\kappa }}
\nc{\la}{{\lambda }}
\nc{\vpi}{{	\varpi		}}
\nc{\vrho}{{	\varrho		}}
\nc{\si}{{	\sigma 		}}
\nc{\ups}{{	\upsilon 	}}
\nc{\vphi}{{	\varphi 	}}
\nc{\om}{{	\omega 		}}

\nc{\Ga}{{\Gamma }}
\nc{\De}{{\Delta }}
\nc{\nab}{{\nabla}}
\nc{\Th}{{\Theta }}
\nc{\La}{{\Lambda }}
\nc{\Si}{{\Sigma }}
\nc{\Ups}{{\Upsilon }}
\nc{\Om}{{\Omega }}

\nc{\Aa}{{	\text{A}	}}
\nc{\Bb}{{	\text{B}	}}
\nc{\Cc}{{	\text{C}	}}
\nc{\Dd}{{	\text{D}	}}
\nc{\Ee}{{	\text{E}	}}
\nc{\Ff}{{	\text{F}	}}
\nc{\Gg}{{	\text{G}	}}
\nc{\Hh}{{	\text{H}	}}
\nc{\Ii}{{	\text{I}	}}
\nc{\Jj}{{	\text{J}	}}
\nc{\Kk}{{	\text{K}	}}
\nc{\Ll}{{	\text{L}	}}
\nc{\Mm}{{	\text{M}	}}
\nc{\Nn}{{	\text{N}	}}
\nc{\Oo}{{	\text{O}	}}
\nc{\Pp}{{	\text{P}	}}
\nc{\Qq}{{	\text{Q}	}}
\nc{\Rr}{{	\text{R}	}}
\nc{\Ss}{{	\text{S}	}}
\nc{\Tt}{{	\text{T}	}}
\nc{\Uu}{{	\text{U}	}}
\nc{\Vv}{{	\text{V}	}}
\nc{\Ww}{{	\text{W}	}}
\nc{\Zz}{{	\text{Z}	}}
\nc{\Xx}{{	\text{X}	}}
\nc{\Yy}{{	\text{Y}	}}

\nc{\bGa}{{	\bbb{\Ga}	}}
\nc{\bA}{{	\bbb{A}		}}
\nc{\bB}{{	\bbb{B}		}}
\nc{\bC}{{	\bbb{C}		}}
\nc{\bD}{{	\bbb{D}		}}
\nc{\bE}{{	\bbb{E}	}}
\nc{\bF}{{	\bbb{F}	}}
\nc{\bG}{{	\bbb{G}	}}
\nc{\bH}{{	\bbb{H}	}}
\nc{\bI}{{	\bbb{I}	}}
\nc{\bJ}{{	\bbb{J}	}}
\nc{\bK}{{	\bbb{K}	}}
\nc{\bL}{{	\bbb{L}	}}
\nc{\bM}{{	\bbb{M}	}}
\nc{\bN}{{	\bbb{N}	}}
\nc{\bO}{{	\bbb{O}	}}
\nc{\bP}{{	\bbb{P}	}}
\nc{\bQ}{{	\bbb{Q}	}}
\nc{\bR}{{	\bbb{R}	}}
\nc{\bS}{{	\bbb{S}	}}
\nc{\bT}{{	\bbb{T}	}}
\nc{\bU}{{	\bbb{U}	}}
\nc{\bV}{{	\bbb{V}	}}
\nc{\bW}{{	\bbb{W}	}}
\nc{\bX}{{	\bbb{X}	}}
\nc{\bY}{{	\bbb{Y}	}}
\nc{\bZ}{{	\bbb{Z}	}}
\nc{\ba}{{	\bbb{a}	}}
			\nc{\bbbb}{{	\bbb{b}	}}
\nc{\bc}{{	\bbb{c}	}}
\nc{\bd}{{	\bbb{d}	}}
			\nc{\bbe}{{	\bbb{e}	}}
			\nc{\bbf}{{	\bbb{f}	}}
\nc{\bg}{{	\bbb{g}	}}
\nc{\bh}{{	\bbb{h}	}}
			\nc{\bbi}{{	\bbb{i}	}}
\nc{\bj}{{	\bbb{j}	}}
			
\nc{\bbk}{{	\k	}}

\nc{\bl}{{	\bbb{l}	}}
\nc{\bm}{{	\bbb{m}	}}
\nc{\bn}{{	\bbb{n}	}}
\nc{\bo}{{	\bbb{o}	}}
\nc{\bp}{{	\bbb{p}	}}
\nc{\bq}{{	\bbb{q}	}}
\nc{\br}{{	\bbb{r}	}}
\nc{\bs}{{	\bbb{s}	}}
\nc{\bt}{{	\bbb{t}	}}
			\nc{\bbbu}{{	\bbb{u}	}}
\nc{\bv}{{	\bbb{v}	}}
\nc{\bw}{{	\bbb{w}	}}
\nc{\bxx}{{	\bbb{x}	}}
\nc{\by}{{	\bbb{y}	}}
\nc{\bz}{{	\bbb{z}	}}

\nc{\sA}{{\mathsf A}}
\nc{\sB}{{\mathsf B}}
\nc{\sC}{{\mathsf C}}
\nc{\sD}{{\mathsf D}}
\nc{\sE}{{\mathsf E}}
\nc{\sF}{{\mathsf F}}
\nc{\sG}{{\mathsf G}}
\nc{\sH}{{\mathsf H}}
\nc{\sI}{{\mathsf I}}
\nc{\sJ}{{\mathsf J}}
\nc{\sK}{{\mathsf K}}
\nc{\sL}{{\mathsf L}}
\nc{\sM}{{\mathsf M}}
\nc{\sN}{{\mathsf N}}
\nc{\sO}{{\mathsf O}}
\nc{\sP}{{\mathsf P}}
\nc{\sQ}{{\mathsf Q}}
\nc{\sR}{{\mathsf R}}
\rc{\sS}{{\mathsf S}}
\nc{\sT}{{\mathsf T}}
\nc{\sU}{{\mathsf U}}
\nc{\sV}{{\mathsf V}}
\nc{\sW}{{\mathsf W}}
\nc{\sX}{{\mathsf X}}
\nc{\sY}{{\mathsf Y}}
\nc{\sZ}{{\mathsf R}}
\nc{\sa}{{\mathsf a}}
\rc{\sb}{{\mathsf b}}
\rc{\sc}{{\mathsf c}}
\nc{\sd}{{\mathsf d}}
\nc{\sg}{{\mathsf g}}
\nc{\sh}{{\mathsf h}}
\nc{\sj}{{\mathsf j}}
\nc{\sk}{{\mathsf k}}
\nc{\sn}{{\mathsf n}}
\nc{\so}{{\mathsf o}}
\nc{\sq}{{\mathsf q}}
\nc{\sr}{{\mathsf r}}
\nc{\su}{{\mathsf u}}
\nc{\sv}{{\mathsf v}}
\rc{\sv}{{\bbb{\mathsf v}}}
\nc{\sw}{{\mathsf w}}
\nc{\sx}{{\mathsf x}}
\nc{\sy}{{\mathsf y}}
\nc{\sz}{{\mathsf z}}

\nc{\toc}{{ 	\small{\tableofcontents} }}
\nc{\addl}{	\addcontentsline{toc}{subsection}	}

\nc{\all}{{ ^{(\alpha)} }}
\nc{\bee}{{ ^{(\beta)} }}
\nc{\gaa}{{ ^{(\gamma)} }}
\nc{\nnnn}{{ ^{ ( n ) } }}
\nc{\nnn}{{ ^{ [ n ] } }}
\nc{\GK}{{  	G(\KK)		}}
\nc{\GO}{{  	G(\OO)		}}

\begin{document}

\title[Modular representations 
 and noncommutative Springer resolution]{
Representations of semi-simple Lie algebras in prime characteristic
 and noncommutative Springer resolution
}

\author{
Roman Bezrukavnikov
}
\address{\small
Department of Mathematics, Massachusetts Institute of Technology, 77 Massachusetts ave.,
Cambridge, MA 02139, USA
}
\email{
bezrukav@@math.mit.edu
}
\author{    Ivan Mirkovi\'c     }
\address{\small
Department of Mathematics and Statistics,
University of Massachusetts, Amherst, MA 01003, USA
}
\email{                mirkovic@@math.umass.edu        }


\thanks{}
\begin{flushright}
\end{flushright}
\maketitle

\nc{\Iaff}{{	I_\aff		}}
\nc{\Waff}{	W_\aff		}
\nc{\Baff}{	\B_\aff		}
\nc{\Waffp}{	W_\aff'		}
\nc{\Baffp}{{	\B_\aff'	}}
\nc{\Waffsc}{{	W_\aff^{sc}	}}
\nc{\Baffsc}{{	\B_\aff^{sc}	}}
\nc{\Bf}{{	\B		}}	
\nc{\Fl}{{	\FF l		}}

\nc{\bio}{{    \bbb{\io}   }}
\nc{\unr}{  _{\mathrm{unr}}  }
\nc{\PPunr}{    _{\PP-{\mathrm{unr}}}  }

\nc{\AZ}{{AZ???}}
\nc{\prr}{{\mathrm{pr}}}
\newcommand{\bbA}{{\mathbb A}}
\nc{\Fr}{{\operatorname{Fr} }}
\nc{\FN}{{   \mathrm{FN}  }}
\newcommand{\Pn}{{{\mathbb P}^n}}
\nc{\WT}{{    \operatorname{wt} }}
\newcommand{\LW}{{^\Theta W}}
\newcommand{\LB}{{^\Theta B}}
\newcommand{\LI}{{^\Theta I}}
\newcommand{\LC}{{^\Theta C}}
\newcommand{\bbI}{{\mathbb I}}
\newcommand{\bbD}{{\mathbb D}}
\nc{\Cone}{{\operatorname{Cone}}}
\nc{\dom}{{dominant???}}
\rc{\mod}{{	\text{mod{\ }}	}}
\nc{\isla}{{????}}

\setpa{-1}

\newcommand{\rnc}{\renewcommand}

\nc{\ad}{{\mbox{\bf{ad}}}}
\nc{\AJ}{{\operatorname{aj}}}
\rc{\Aut}{{\operatorname{Aut}}}
\nc{\Bls}{{{\cal B}ls}}
\nc{\Boxtimes}{{\fbox{$\times$}}}
\nc{\blt}{{\bullet}}
\nc{\bSt}{{\mbox{\bf{St}}}}
\nc{\card}{{\operatorname{card}}}
\nc{\Cch}{{\check{C}}}
\nc{\chara}{{\operatorname{char}}}
\nc{\CHom}{{\cal{H}om}}
\rc{\Coker}{{\operatorname{Coker}}}
\rc{\codim}{{\operatorname{codim}}}

\nc{\cSgn}{{\cal{S}gn}}
\nc{\depth}{{\operatorname{depth}}}
\nc{\dirlim}{{\underset{\rightarrow}{\operatorname{lim}}}}
\nc{\dotbox}{{\overset{\bullet}{\boxtimes}}}
\nc{\dotimes}{{\overset{\bullet}{\otimes}}}
\nc{\Ed}{{\operatorname{Edge}}}
\rc{\Ext}{{\operatorname{Ext}}}
\nc{\Fac}{{\cal{F}ac}}
\nc{\Fun}{{\operatorname{F}}}
\nc{\FS}{{\cal{FS}}}
\rc{\Hom}{{\operatorname{Hom}}}
\nc{\had}{{{\hat{\mbox{\bf{ad}}}}}}
\nc{\hgt}{{\operatorname{ht}}}
\nc{\Id}{{\operatorname{Id}}}
\nc{\id}{{\operatorname{id}}}
\nc{\Ima}{{\operatorname{Im}}}
\nc{\ind}{{\operatorname{ind}}}
\nc{\Ind}{{\operatorname{Ind}}}
\nc{\infi}{{\operatorname{inf}}}
\nc{\infh}{{\frac{\infty}{2}}}
\nc{\invlim}{{\underset{\leftarrow}{\operatorname{lim}}}}
\nc{\Jac}{{{\cal J}ac}}
\rc{\Ker}{{\operatorname{Ker}}}
\nc{\lcm}{{\operatorname{lcm}}}
\nc{\Locsys}{{{\cal L}ocsys}}
\nc{\Map}{{{\cal M}ap}}
\nc{\Mor}{{\operatorname{Mor}}}
\nc{\MS}{{\cal{MS}}}
\rc{\Ob}{{\operatorname{Ob}}}
\nc{\opp}{{\operatorname{opp}}}
\nc{\Or}{{{\cal O}r}}
\nc{\Ord}{{{\cal O}rd}}
\nc{\Part}{{{\cal P}art}}
\rc{\PGL}{{\operatorname{PGL}}}
\rc{\Pic}{{\operatorname{Pic}}}
\rc{\Rep}{{{\cal{R}}ep}}
\rc{\rk}{{\operatorname{rk}}}
\nc{\Sets}{{{\cal{S}}ets}}
\nc{\Sew}{{{\cal{S}}ew}}
\nc{\sgn}{{\operatorname{sgn}}}
\nc{\Sh}{{{\cal S}h}}
\nc{\Sign}{{{\cal S}ign}}
\nc{\Spe}{{\mbox{\bf{Sp}}}}
\nc{\supr}{{\operatorname{sup}}}
\nc{\Supp}{{\operatorname{Supp}}}
\rc{\supp}{{\operatorname{supp}}}
\nc{\Teich}{{{\cal{T}}eich}}
\nc{\tFS}{{\widetilde{\cal{FS}}}}
\rc{\Tor}{{\operatorname{Tor}}}
\nc{\totimes}{{\tilde{\otimes}}}
\rc{\tr}{{\operatorname{tr}}}
\nc{\tRep}{{\widetilde{{\cal R}ep}}}
\nc{\tTeich}{{\widetilde{{\cal T}eich}}}
\nc{\Vect}{{{\cal V}ect}}
\nc{\Ve}{{\operatorname{Vert}}}
\nc{\wt}{{\widetilde}}


\rc{\bo}{{\mbox{\bf{0}}}}
\nc{\One}{{\mbox{\bf{1}}}}
\rc{\one}{{\mbox{\bf{1}}}}

\nc{\BA}{{\Bbb A}}
\rc{\bA}{{\overline{A}}}
\rc{\ba}{{\mbox{\bf{a}}}}

\nc{\baB}{{\overline{B}}}
\nc{\baeta}{{\bar{\eta}}}
\nc{\baJ}{{\bar{J}}}
\rc{\bB}{{\mbox{\bf{B}}}}
\rc{\bc}{{\mbox{\bf{c}}}}
\rc{\bC}{{\overline{C}}}
\nc{\BC}{{\Bbb{C}}}
\nc{\bCC}{{\overline{\cal{C}}}}
\nc{\bCM}{{\overline{\cal{M}}}}
\nc{\BD}{{\overline{D}}}
\rc{\bd}{{\mbox{\bf{d}}}}
\nc{\BE}{{\overline{E}}}
\nc{\BF}{{\overline{F}}}
\rc{\bF}{{\mbox{\bf{F}}}}
\rc{\bg}{{\mbox{\bf{g}}}}
\rc{\bG}{{\mbox{\bf{G}}}}
\rc{\bG}{{\bf{G}}}
\nc{\BG}{{\Bbb G}}
\nc{\bGamma}{{\overline{\Gamma}}}
\nc{\bbH}{{\bar{\mbox{\bf{H}}}}}
\rc{\bH}{{\mbox{\bf{H}}}}
\rc{\bbi}{{\mbox{\bf{i}}}}
\rc{\bI}{{\mbox{\bf{I}}}}
\rc{\bL}{{\mbox{\bf{L}}}}
\nc{\BL}{{\Bbb{L}}}
\nc{\blambda}{{\bar{\lambda}}}
\rc{\bM}{{\mbox{\bf{M}}}}
\nc{\bmu}{{\vec{\mu}}}
\rc{\bN}{{\mbox{\bf{N}}}}
\nc{\BN}{{\Bbb{N}}}
\nc{\bnu}{{\mbox{\boldmath{${\nu}$}}}}
\nc{\bof}{{\mbox{\bf{f}}}}
\nc{\BP}{{\Bbb P}}
\rc{\bP}{{\mbox{\bf{P}}}}
\nc{\BPO}{{\overset{\circ}{\BP}}}
\nc{\BQ}{{\Bbb Q}}
\rc{\bq}{{\mbox{\bf{q}}}}
\nc{\BR}{{\Bbb{R}}}
\rc{\bR}{{\mbox{\bf{R}}}}
\rc{\br}{{\mbox{\bf{r}}}}
\nc{\breta}{{\bar{\eta}}}
\rc{\bs}{{\mbox{\bf{s}}}}
\rc{\bS}{{\mbox{\bf{S}}}}
\rc{\bt}{{\mbox{\bf{t}}}}
\rc{\bU}{{\mbox{\bf{U}}}}
\rc{\bV}{{\mbox{\bf{V}}}}
\rc{\bu}{{\mbox{\bf{u}}}}
\nc{\BUpsilon}{{\bar{\Upsilon}}}
\rc{\bw}{{\mbox{\bf{w}}}}
\rc{\bX}{{\mbox{\bf{X}}}}
\nc{\BZ}{{\Bbb{Z}}}
\rc{\bz}{{\mbox{\bf{z}}}}
\rc{\bZ}{{\mbox{\bf{Z}}}}
\nc{\bzero}{\mbox{\boldmath{$0$}}}

\nc{\CA}{{\cal A}}
\nc{\CAD}{{\overset{\bullet}{\cal{A}}}}
\nc{\CAO}{{\overset{\circ}{\cal{A}}}}
\nc{\CB}{{\cal B}}
\nc{\CalD}{{\cal D}}
\nc{\CE}{{\cal E}}
\nc{\CF}{{\cal F}}
\nc{\CG}{{\cal G}}
\nc{\CH}{{\cal H}}
\nc{\CI}{{\cal I}}
\nc{\CID}{{\overset{\bullet}{\cal{I}}}}
\nc{\CJ}{{\cal J}}
\nc{\CK}{{\cal K}}
\nc{\CL}{{\cal L}}
\nc{\CM}{{\cal M}}
\nc{\CN}{{\cal N}}
\nc{\CO}{{\cal O}}
\nc{\CP}{{\cal P}}
\nc{\CPO}{{\overset{\circ}{\cal{P}}}}
\nc{\CQ}{{\cal Q}}
\nc{\CR}{{\cal R}}
\nc{\CS}{{\cal S}}
\nc{\CT}{{\cal T}}
\nc{\CTD}{{\overset{\bullet}{\cal{T}}}}
\nc{\CTPO}{{\overset{\circ}{\cal{T}\cal{P}}}}
\nc{\CU}{{\cal{U}}}
\nc{\CV}{{\cal V}}
\nc{\CW}{{\cal W}}
\nc{\CX}{{\cal X}}
\nc{\CY}{{\cal Y}}
\nc{\CZ}{{\cal Z}}

\nc{\dCL}{{\overset{\bullet}{\cal{L}}}}
\nc{\ddelta}{{\overset{\bullet}{\delta}}}
\nc{\dfu}{{\overset{\bullet}{\frak{u}}}}
\nc{\dlambda}{{\overset{\bullet}{\lambda}}}
\nc{\DO}{{\overset{\circ}{D}}}
\nc{\dpar}{{\partial}}
\nc{\dS}{{\overset{\bullet}{S}}}
\nc{\dT}{{\overset{\bullet}{T}}}


\nc{\hCH}{{\hat{\cal{H}}}}
\nc{\hCI}{{\hat{\cal{I}}}}
\nc{\hfC}{{\hat{\frak{C}}}}
\nc{\hfg}{{\hat{\frak{g}}}}
\nc{\hL}{{\hat{L}}}
\nc{\OH}{{\overset{\circ}{H}}}
\nc{\hpsi}{{\hat{\psi}}}
\nc{\hx}{{\hat{x}}}

\nc{\jo}{{\overset{\circ}{j}}}

\nc{\phid}{{\overset{\bullet}{\phi}}}

\nc{\tA}{{\tilde{A}}}
\nc{\ta}{{\tilde{a}}}
\nc{\tB}{{\tilde{B}}}
\nc{\tb}{{\tilde{b}}}
\nc{\tBP}{{\tilde{\BP}}}
\nc{\tC}{{\tilde{C}}}
\nc{\tc}{{\tilde{c}}}
\nc{\tCA}{{\tilde{\cal{A}}}}
\nc{\tCC}{{\tilde{\cal{C}}}}
\nc{\tCH}{{\tilde{\cal{H}}}}
\nc{\tCI}{{\tilde{\cal{I}}}}
\nc{\tCO}{{\tilde{\cal{O}}}}
\nc{\tCP}{{\tilde{\cal{P}}}}
\nc{\tCT}{{\tilde{\cal{T}}}}
\nc{\tD}{{\tilde{D}}}
\nc{\tDelta}{{\tilde{\Delta}}}
\nc{\tE}{{\tilde E}}
\nc{\tF}{{\tilde F}}
\nc{\tfD}{{\tilde{\frak{D}}}}
\nc{\tfF}{{\tilde{\frak{F}}}}
\nc{\tff}{{\tilde{\frak{f}}}}
\nc{\tfu}{{\tilde{\frak{u}}}}
\nc{\tJ}{{\tilde{J}}}
\nc{\tj}{{\tilde{j}}}
\nc{\tK}{{\tilde K}}
\nc{\tL}{{\tilde{L}}}
\nc{\tM}{{\tilde{M}}}
\nc{\tP}{{\tilde{P}}}
\nc{\tPhi}{{\tilde{\Phi}}}
\nc{\tpi}{\tilde{\pi}}
\nc{\TPO}{{\overset{\circ}{T\BP}}}
\nc{\tR}{{\tilde{R}}}
\nc{\tS}{{\tilde S}}
\nc{\tT}{{\tilde{T}}}
\nc{\ttau}{{\tilde{\tau}}}
\nc{\ttheta}{{\tilde{\theta}}}
\nc{\tU}{{\tilde{U}}}
\nc{\tUpsilon}{{\tilde{\Upsilon}}}
\nc{\tW}{{\tilde W}}
\nc{\ty}{{\tilde y}}
\nc{\tY}{{\tilde Y}}
\nc{\txi}{{\tilde{\xi}}}

\nc{\UD}{{\overset{\bullet}{U}}}
\nc{\UO}{{\overset{\circ}{U}}}

\nc{\vA}{{\vec{A}}}
\nc{\valpha}{{\vec{\alpha}}}
\nc{\vbeta}{{\vec{\beta}}}
\nc{\vc}{{\vec{c}}}
\nc{\vD}{{\vec{D}}}
\nc{\vd}{{\vec{d}}}
\nc{\vgamma}{{\vec{\gamma}}}
\nc{\vK}{{\vec{K}}}
\nc{\vlambda}{{\vec{\lambda}}}
\nc{\vmu}{{\vec{\mu}}}
\nc{\vnu}{{\vec{\nu}}}
\nc{\vo}{{\vec{0}}}
\nc{\vu}{{\vec{u}}}
\nc{\vx}{{\vec{x}}}
\nc{\vy}{\vec{y}}
\nc{\vzero}{\vec{0}}

\nc{\XO}{{\overset{\circ}{X}}}

\nc{\ya}{{\operatorname{aj}}}


\rc{\nen}{\newenvironment}
\nc{\ol}{\overline}
\nc{\ul}{\underline}

\nc{\Lra}{\Longrightarrow}

\nc{\Llra}{\Longleftrightarrow}
\nc{\hra}{\hookrightarrow}
\rc{\iso}{\overset{\sim}{\lra}}
\nc{\rlh}{\rightleftharpoons}


\nc{\PS}{{\cal{PS}}}
\nc{\oCG}{{\overline{\cal G}}}
\nc{\oCQ}{{\overline{\cal Q}}}
\nc{\oCZ}{{\overline{\cal Z}}}
\nc{\dZ}{{\overset{\bullet}{\cal Z}}{}}
\nc{\ddZ}{{\ddot{\cal Z}}{}}
\nc{\oZ}{{\overset{\circ}{\cal Z}}{}}

\nc{\dP}{{\overset{\bullet}{\cal P}}{}}
\nc{\oP}{{\overset{\circ}{\cal P}}{}}
\nc{\oQ}{{\overset{\circ}{\cal Q}}{}}
\nc{\obp}{{\overset{\circ}{{\bf p}}}{}}
\nc{\tbj}{{\tilde{\bf j}}{}}
\nc{\tbp}{{\tilde{\bf p}}{}}
\nc{\tfC}{{\widetilde{\frak C}}{}}
\nc{\tfE}{{\widetilde{\frak E}}{}}
\nc{\tfj}{{\widetilde{\frak j}}{}}
\nc{\tfQ}{{\widetilde{\frak Q}}{}}
\nc{\tfp}{{\widetilde{\frak p}}{}}
\nc{\ofQ}{{\overset{\circ}{{\frak Q}}}{}}
\nc{\tGQ}{{\widetilde{\cal{GQ}}}{}}
\nc{\oGQ}{{\overset{\circ}{\cal{GQ}}}{}}
\nc{\ooGQ}{{\overset{\circ\circ}{\cal{GQ}}}{}}
\nc{\oGZ}{{\overset{\circ}{\cal{GZ}}}{}}
\nc{\tGZ}{{\widetilde{\cal{GZ}}}{}}

\nc{\Ue}{{U_\varepsilon}}
\nc{\Upe}{{\Upsilon_\varepsilon}}
\nc{\crho}{{\check{\rho}}}
\nc{\ctheta}{{\check{\theta}}}

\unitlength=1pt
\setlength{\baselineskip}{18pt}
\setlength{\parindent}{0cm}
\setlength{\parskip}{6pt}
\setlength{\textwidth}{16cm}
\setlength{\textheight}{21cm}
\setlength{\oddsidemargin}{0.1in}
\setlength{\evensidemargin}{0.1in}
\setlength{\headheight}{30pt}
\setlength{\headsep}{40pt}
\setlength{\topmargin}{10pt}
\setlength{\marginparwidth}{0pt}


%



\nc{\tww}{{ {}^{*,1}    }}
\nc{\zhc}{{ \fZ_{HC} }}
\nc{\zfr}{{ \fZ_{Fr} }}
\nc{\fzp}{{ \fz_{p} }}
\nc{\pp}{ ^{[p]} }

\nc{\utx}{{     \DD_X       }}
\nc{\utb}{{     \DD_\BB     }}

\nc{\dx}{   \DD_X       }
\nc{\db}{   \DD_\BB     }

\nc{\z}{{   ^{\bbb 0}   }}



\nc{\PPhi}{{	\bbb\Phi	}}

\nc{\hzc}{{ \hatt{(\chi,0)} }}
\nc{\hrc}{{ \hatt{\chi,-\rho}   }}
\nc{\mr}{{  {-\rho} }}

\nc{\ob}{{  \OO_\BB     }}               
\nc{\uc}{{  u_\chi      }}
\nc{\Fp}{{  {\Bbb F}_p  }}
\nc{\frx}{{ Fr_X        }}

\nc{\Fnew}{ \mbox{Fr}_X }
\nc{\Frx}{  \Fnew       }

\newtheorem{Thm}{Theorem}
\newtheorem{Prop}{Proposition}
\newtheorem{Sublem}{Sublemma}
\newtheorem{Fact}{Fact}
\newtheorem{Claim}{Claim}
\theoremstyle{remark}
\newtheorem{Rem}{Remark}

\nc{\proof}{{\it Proof. }}

\newcommand{\imbed}{\hookrightarrow}
\renewcommand{\iso}{{\tii \longrightarrow}}
\newcommand{\isol}{{\tii \longleftarrow}}
\newcommand{\To}{\longrightarrow}

\newcommand{\Lotimes}{\overset{\rm L}{\otimes}}
\newcommand{\Ltimes}{\overset{\rm L}{\times}}
\newcommand{\Rtimes}{\overset{\rm R}{\times}}

\newcommand{\oplusl}{\bigoplus\limits}
\newcommand{\cupl}{\bigcup\limits}

\def\square{\hbox{\vrule\vbox{\hrule\phantom{o}\hrule}\vrule}}

\renewcommand{\E}{{\mathcal E}}
\renewcommand{\V}{{\mathcal V}}
\renewcommand{\N}{{\mathcal N}}
\renewcommand{\O}{{\mathcal O}}

\newcommand{\Etil}{\widetilde{\mathcal E}}

\renewcommand{\G}{{\mathcal G}}
\rc{\L}{{\tx{L}}}
\renewcommand{\bO}{{\bf O}}

\newcommand{\Ntil}{{    \tii\NN         }}
\newcommand{\Ptil}{ \tilde P        }
	\newcommand{\Phitil}{{  \tilde\Phi      }}
\newcommand{\Ftil}{{    \tilde F        }}
\newcommand{\Dtil}{{    \tii \D         }}

\newcommand{\Nhat}{{        \hatt\NN    }}
\newcommand{\Lambdahat}{{   \hatt \La   }}
\newcommand{\Chat}{{\hat C}}
\newcommand{\Uhat}{{\hat U}}
\newcommand{\Dhat}{{\hat D}}
\newcommand{\Bhat}{{\hat \B}}
\newcommand{\Shat}{{\hat S}}

\newcommand{\Ubar}{{\overline U}}

\newcommand{\Db}{   \text{D}^{b}    }
\newcommand{\Dmin}{ \text{D}^-  }
\let\Dm\Dmin
\newcommand{\Dmt}{{ D^b_\D(\Dtil)       }}
\newcommand{\Um}{{  D^b(U_0)        }}
\newcommand{\Umt}{{ D^b_0(U)        }}
\newcommand{\Loc}{{ \mathcal L      }}

\nc{\til}{  \tilde  }

\nc{\Mod}{{ mod }}
\newcommand{\Dmod}{{        \Mod(\DD)       }}
\newcommand{\Dtmodall}{{    \Mod(\Dtil) }}      
\newcommand{\Dmodall}{{     \Mod(\Dtil) }}
\newcommand{\Dtmod}{{       \Mod(\D_{\hat 0})   }}

\newcommand{\Umodall}{{     \Mod(U)             }}
\newcommand{\Umodo}{{       \Mod(U_0)       }}
\newcommand{\Umod}{{        \Mod(U_{\hat 0})    }}
\newcommand{\Umodallfg}{{   \Mod^{fg}(U)                }}
\newcommand{\Umodofg}{{     \Mod^{fg}(U_0)      }}
\newcommand{\Umodfg}{{      \Mod^{fg}(U_{\hat 0})   }}

\newcommand{\Wa}{{      W_{\aff}     }}
\newcommand{\Wap}{{     W_{\aff}'    }}

\renewcommand{\t}{{\mathfrak t}}
\newcommand{\n}{{\mathfrak n}}
\newcommand{\g}{{\mathfrak g}}
\renewcommand{\h}{{\mathfrak h}}
\renewcommand{\b}{{\mathfrak b}}
\newcommand{\m}{{\mathfrak m}}

\renewcommand{\bu}{{        \bullet     }}

\newcommand{\gal}{\check{\ }}

\newcommand{\Zet}{{\mathbb Z}}

\renewcommand{\sur}{\twoheadrightarrow}

\newcommand{\eps}{\epsilon}

\renewcommand{\M}{{\mathcal M}}
\renewcommand{\N}{{\mathcal N}}
\renewcommand{\O}{{\mathcal O}}

\newcommand{\gtil}{{\tii{\mathfrak g}}}

\nc{\diff}{{    \bbb{\bss}  }}

\nc{\hla}{{\hatt\la}}

\renewcommand{\S}{{\mathcal S}}

\renewcommand{\Fp}{{{\mathbb F}_p}}

\newcommand{\Gm}{{{\mathbb G}_m}}

\rc{\U}{{   \tx{U}  }}
\nc{\uU}{{  \U      }}

\rc{\i}{\item}

\newcommand{\gt}{{\tilde \g }}
\renewcommand{\bY}{{\bf Y}}
\newcommand{\bpi}{{\bf \pi}}

\renewcommand{\P}{{\mathcal P}}

\nc{\Kbar}{{\bar{\KK}}}

\rc{\aa}{\pr\overset}            




\renewcommand{\iso}{{\tii \longrightarrow}}

\def\square{\hbox{\vrule\vbox{\hrule\phantom{o}\hrule}\vrule}}

\renewcommand{\N}{{\mathcal N}}

\renewcommand{\O}{{\mathcal O}}

\renewcommand{\G}{{\mathcal G}}

\nc{\hattt}{{               }}

\renewcommand{\h}{{\mathfrak h}}
\renewcommand{\t}{{\mathfrak t}}
\renewcommand{\bu}{{            \bullet         }}

\newcommand{\Aone}{{\mathbb A}^1}
\newcommand{\Pone}{{\mathbb P}^1}
\newcommand{\Ql}{{  {\mathbb Q}_l       }}
\newcommand{\Qlb}{{ \bar{\mathbb Q}_l   }}
\newcommand{\Qu}{{\mathbb Q}}
\newcommand{\Qubar}{{\overline{\mathbb Q}}}
\newcommand{\Ce}{{\mathbb C}}
\newcommand{\bbS}{{\mathbb S}}
\newcommand{\RE}{{\mathbb R}}

\renewcommand{\M}{{\mathcal M}}
\renewcommand{\N}{{\mathcal N}}
\renewcommand{\O}{{\mathcal O}}
\renewcommand{\F}{{\mathcal F}}

\renewcommand{\b}{{\mathfrak b}}
\renewcommand{\t}{{\mathfrak t}}


\nc{\aff}{{ \tx{aff}    }}
\rc{\RGa}{{ \Rr\Gamma   }}

\nc{\uHom}{{ \underline{\rm Hom}   }}
\nc{\uExt}{{ \underline{\rm Ext}   }}
\nc{\uuHom}{{ \bf Hom  }}

\nc{\PR}{{\bf pr}}
\nc{\AS}{AS}




\nc{\bDD}{{ \bbb{\DD}   }}




\theoremstyle{remark}









\renewcommand{\A}{{\mathcal A}}
\renewcommand{\cD}{{\mathcal D}}
\renewcommand{\O}{{\mathcal O}}
\newcommand{\capl}{\bigcap\limits}

\newcommand{\cH}{{\mathcal H}}
\newcommand{\Nt}{{\widetilde{\N}}}
\renewcommand{\t}{{\mathfrak t}}
\renewcommand{\b}{{\mathfrak b}}
 \newcommand{\p}{{\mathfrak p}}

\newcommand{\Rhat}{{\widehat R}}

\renewcommand{\P}{{\mathcal P}}
\newcommand{\cS}{{\mathcal S}}


\rc{\bl}{{ \fbox{$\blacklozenge\to$} }}
\nc{\wl}{{ \fbox{$@<<<\lozenge$}     }}

\let\iso\con
\nc{\IIII}{I}   \let\Sigma\IIII

\nc{\mud}{{ \aa{\bu}\mu     }}
\nc{\tI}{{  \text{\bf I}        }}
\nc{\mn}{   _{\mu|\nu}      }

\nc{\wti}{\widetilde w}

\nc{\Sti}{\widetilde S}

\nc{\Pb}{{	\barr P	}}

\nc{\kk}{\k}
\nc{\Ch}{\bf Check!}

\nc{\exoticsheaves}{{exotic\ sheaves}}

\newcommand{\Gbu}{{	\aa{\bu}G		}}
\newcommand{\buu}[1]{{	\aa{\bu}{#1}		}}

\newcommand{\Ztil}{{	\aa{\bu}G_e		}}
\newcommand{\Gatil}{{	\aa{\bu} \Gamma		}}

\newcommand{\bfB}{{\bf B}}
\newcommand{\bfBe}{\bfB_e}
	\newcommand{\bfBp}{\bfB'}
\newcommand{\bfBep}{\bfB_e'}

\newcommand{\bfBSe}{\bfB_{S}}
\newcommand{\bfBS}{\bfB_{S}}

\newcommand{\LG}{{^LG}}

\nc{\bfT}{{\bf T}}

\rc{\BC}{{	\BB\CC		}}
\nc{\FGP}{{	\tx{FGP}	}}

\nc{\hz}{{ 	\hatt{0} 	}}
\nc{\hzk}{{ 	^\hz_k 		}}
\nc{\Ehzk}{{ 	\E^\hz_k 		}}

\nc{\NO}{{	\sN\sO	}}

\newcommand{\Alc}{\operatorname{Alc}}

\rc{\bmu}{ \begin{multline*} }
\newcommand{\fmu}{ \end{multline*}   }

\nc{\WaffC}{{ 	W_\aff^{Cox}	}}	
\nc{\BaffC}{{ 	{\mathbb B}_\aff^{Cox}	}}

\nc{\beS}{{	\be_S	}}
\rc{\bee}{{	\be_e	}}

\nc{\cG}{{	\ch G	}}

\nc{\GR}{{	G_R	}}
\nc{\gR}{{	_R	}}

\nc{\Tau}{{	\TT	}}

\rc{\st}{{ 	\bbb{\star}			}}
\nc{\stt}{{ 	\st\st				}}
\nc{\sttt}{{ 	\st\st\st			}}
\nc{\zsttt}{{ 	(0_\sttt)			}}
\nc{\osttt}{{ 	(1_\sttt)			}}
\nc{\tsttt}{{ 	(2_\sttt)			}}

\nc{\bst}{{ 	\bbb{(\st)}			}}
\nc{\bstt}{{ 	\bbb{(\stt)}			}}
\nc{\bsttt}{{ 	\bbb{(\sttt)}			}}

\nc{\Fq}{{	\F_q	}}

\nc{\fGm}{{	\fG_\fm	}}
\rc{\bz}{{	\bbb{z}		}}

\nc{\Ske}{{	S_{\k,e}		}}
\nc{\Stilke}{{	\tii{S_{\k,e}}		}}
\nc{\Stilkep}{{	\tii{S_{\k,e}}'		}}
\nc{\BBke}{{	\BB_{\k,e}		}}
\nc{\BBkep}{{	\BB_{\k,e}'		}}
\nc{\hBBke}{{	\hatt{\BB_{\k,e}}	}}
\nc{\hBBkep}{{	\hatt{\BB_{\k,e}}'	}}

\nc{\SCe}{{	S_{\C,e}		}}
\nc{\StilCe}{{	\tii{S_{\C,e}}		}}
\nc{\StilCep}{{	\tii{S_{\C,e}}'		}}
\nc{\BBCe}{{	\BB_{\C,e}		}}
\nc{\BBCep}{{	\BB_{\C,e}'		}}
\nc{\hBBCe}{{	\hatt{\BB_{\C,e}}	}}
\nc{\hBBCep}{{	\hatt{\BB_{\C,e}}'	}}

\newcommand{\Ctil}{{\widetilde C}}
\newcommand{\Ltil}{\widetilde {\LL}}

\newcommand{\Stil}{\tii S}
\newcommand{\Stilp}{\tii S'}
\nc{\Se}{{	S_{e}		}}
\nc{\Stile}{{	\tii{S_{e}}		}}
\nc{\Stilep}{{	\tii{S_{e}}'		}}
\nc{\BBe}{{	\BB_{e}		}}
\nc{\BBep}{{	\BB_{e}'		}}
\nc{\hBBe}{{	\hatt{\BB_{e}}	}}
\nc{\hBBep}{{	\hatt{\BB_{e}}'	}}

\nc{\RCp}{{	R_C^+	}}
\nc{\RC}{{	R_C	}}

\nc{\Cohe}{{	\EE{}coh			}}

\nc{\sff}[1]{{\sf (#1)}}

\nc{\Abig}{	A		}	
\nc{\Az}{	A^0		}	

\rc{\bD}{{	\bbb{\sD}		}}

\nc{\chro}{{	\ch\rho	}}

\nc{\Fc}{{\aa{o}\sF}}
\nc{\Gc}{{\aa{o}\sG}}
\nc{\bAA}{{	\bbb\AA	}}
\nc{\IC}{{  \tx{IC}   }}

\begin{flushright}
{\em
To Joseph Bernstein with
admiration
and
gratitude
}
\end{flushright}

\toc



\setse{-1}	
\se{	
Introduction	
}

\begin{abstract}

We prove most of
Lusztig's conjectures   from \cite{Kth2},
including the
existence of a canonical basis in the Grothendieck group of a Springer
fiber. The conjectures also predict that this basis controls numerics of
 representations  of the Lie algebra of a
 semi-simple algebraic group over an algebraically closed field of positive
characteristic. We check this for almost all
characteristics.

To this end we construct a non-commutative resolution of the nilpotent cone
which is derived equivalent to the Springer resolution.
On the one hand, this noncommutative resolution is shown
to be compatible with the positive characteristic
localization equivalences of \cite{BMR1}. On the other hand, it is compatible
with the t-structure arising from an equivalence with the derived
category of perverse sheaves on the affine flag variety of the
Langlands dual group
\cite{AB} inspired by local geometric Langlands duality. This allows
one to apply Frobenius purity theorem of \cite{BBD}
to deduce
the desired properties of the basis.

We expect the noncommutative counterpart of the Springer resolution
to be of independent interest from the perspectives
of  algebraic geometry and geometric Langlands duality.

\end{abstract}


Let $G$ be a reductive group over an algebraically closed field
 $\k$ of characteristic $p>h$ and $\g$ be its Lie algebra (here $h$ denotes
the Coxeter number of $G$).
Let $\PP$ be a partial flag variety and consider the space $\gt_\PP$
of pairs of a parabolic subalgebra $\fp\in\PP$ and an element in it.
For the full flag variety $\BB$ we usually denote $\tii\fg_\BB$
simply by
$\gt$.
We have a map
$\mu_\PP:\gt_\PP\to \g$.

In \cite{BMR1}, \cite{BMR2} we have shown that
the derived category of $\g$-modules
with a fixed generalized central character is equivalent to
the derived category of coherent sheaves on $\gt  
_\PP$ set-theoretically
supported on
$\PP_e=\mu_\PP^{-1}(e)$;
 here the partial flag variety $\PP$ and  $e\in \g   
$ depend on the central
character.
A numerical consequence of this equivalence is an isomorphism
between the Grothendieck groups of the
abelian categories $mod^{fg}(U(\g)_\sigma)$ and $Coh(\PP_e)$,
where $U(\g)_\sigma$ is the
quotient of the enveloping algebra by a central character $\sigma$
and $mod^{fg}$ denotes
the category of finite dimensional (equivalently, finitely generated) modules.
 This
implies, in particular, that the number of irreducible representations
with a fixed central character $\si$
equals the rank of the Grothendieck group
of $Coh(\PP_e)$, which is known to coincide with the sum of Betti numbers
of $\PP_e$.

To derive  more precise information
about numerical invariants of $\g$-modules
one needs a characterization of the elements in  $K^0(Coh(\PP_e))$
which correspond to irreducible $\g$-modules
and their  projective covers.
Such a characterization is suggested
by the work of Lusztig
\cite{Kth2}. In
{\em loc. cit.} he describes certain properties of a
basis in the Grothendieck group of a Springer fiber and conjectures
that a basis with such properties exists and controls (in a certain precise
sense) numerical invariants of irreducible $U(\g)_\sigma$ modules.
(He also shows that a basis with such properties is essentially unique).
The properties of a basis are similar to those enjoyed by
Kazhdan-Lusztig bases of a  Hecke algebra and canonical  bases
in modules over a quantum group, for this reason we will refer to a basis
satisfying Lusztig's axioms as a {\em canonical basis}.

In the present paper we prove most of the
conjectures from \cite{Kth2}.
The first step is the construction of a non-commutative counterpart of the
Springer resolution as a lift of modular representation categories
to characteristic zero. By this we mean a certain noncommutative algebra
$\Az$ defined canonically up to a Morita equivalence. The 
 center of the algebra is identified
with the ring $\O(\N)$ of regular functions on the nilpotent cone
$\N\subset \g^*$, where $\g$ is taken over $R=\Zet[\frac{1}{h!}]$.
This noncommutative resolution is canonically derived equivalent
to the ordinary Springer resolution, i.e. it comes with a canonical
equivalence of triangulated categories $D^b(mod^{fg}(\Az))
\cong D^b(Coh(\Nt))$
where $\Nt$ is the cotangent bundle to the flag variety. Furthermore,
for $\k$ as above and
any $e\in \N(\k)$, the base change $\Az\otimes _{\O(\N)}\k_e$
is canonically Morita equivalent to a central reduction of
$U(\g_\k)$.

The above properties of  $\Az$ imply that
the numerics of non-restricted modular
representation categories is independent of
(sufficiently large) characteristic and show
that $\Az$ provides
a lifting of such representation categories to characteristic zero.
To put things into perspective,
 recall that a similar construction for representations
of the algebraic group $G_\k$ (this setting is very close
 to restricted representations
of the Lie algebra $\g_\k$) was obtained in \cite{AJS}. In that case
the resulting category in characteristic zero turns out to be equivalent
to representations of a quantum group at a root of unity.
We expect that a similar statement holds for non-restricted
Lie algebra modules considered in the present work, see
 Conjecture \ref{conj_Uq} below. Apart from that Conjecture, we avoid
 quantum groups in this paper.

Our method of construction of the noncommutative resolution $\Az$ is
based on an action of the affine braid group $\Baff$ on the derived
categories $D^b(Coh(\Nt))$, $D^b(Coh(\gt_\BB))$.
Here the action
of the generators of $\Baff$ is described by certain simple
correspondences.
The fact that the corresponding functors
obey the relations of $\Baff$ is proven in \cite{BR}.
The algebra $\Az$ is determined (uniquely up to a Morita equivalence)
by
the t-structure on $D^b(Coh(\Nt))$
corresponding to the tautological one under the equivalence
with $D^b(mod^{fg}(\Az))$. This t-structure
 is characterized in terms of the action of $\Baff$.
The comparison with modular localization and the proof of  existence
of a t-structure with required properties is based on compatibility
of the $\Baff$ action with {\em intertwining} (or {\em shuffling}) functors
on the derived categories of modular representations.
Notice that the latter are closely connected with ``translation
through the wall'' functors, thus translation functors play a prominent
role in our argument.
The use of translation functors to establish independence of
the category of modular representations of characteristic goes
back (at least) to \cite{AJS}.

From the arguments alluded to above one can derive that the basis
in the Grothendieck group of a Springer fiber corresponding to
irreducible $\g_\k$ modules satisfies all the axioms of a canonical
basis except for one, the so-called {\em asymptotic orthogonality}
property. The latter is reduced to
certain compatibility between the above t-structures and  the
 multiplicative group action on  Slodowy slices.
It says that the grading on the slice algebras, i.e. the algebras
``controlling''
 the derived category of coherent sheaves on the resolution of a
Slodowy slice,
can be arranged to be positive. By this we mean that
  components of negative degrees in the algebra vanish, while
 the degree zero component is semi-simple.
 An analogous reformulation of Kazhdan-Lusztig conjectures
 is due to Soergel. Another feature parallel to Kazhdan-Lusztig theory
is Koszul property of the slice algebras, see 
 \cite{BGS} for the corresponding facts about category $O$.

Properties of this type are usually deduced from a Theorem of \cite{BBD} about
weights of Frobenius acting on the stalks of $l$-adic intersection cohomology
sheaves. Our proof also follows this strategy. The $l$-adic sheaves
are brought into the picture by the result of \cite{AB} which
provides an equivalence between the derived category of $G$-equivariant
sheaves on $\Nt$ (over a field of characteristic zero) and a
certain subcategory of the derived category
of  constructible sheaves on
 the  affine flag variety $\Fl$ of the Langlands dual group.
This result is a categorical counterpart of one of the key
ingredients in the proof
of the tamely ramified local Langlands conjecture.

We show that the t-structure of perverse sheaves on $\Fl$
is compatible with the t-structure coming from the equivalence
with $D^b(mod^{fg}(\Az))$.
This is achieved by interpreting the $\Baff$
action on the perverse sheaves side as the geometric counterpart
of the action of elements in the standard basis of the affine Hecke
algebra on the anti-spherical module.

\medskip
Thus the key step in our argument is compatibility between the two
t-structures on $D^b(Coh(\Nt))$, one coming from modular representations
via the equivalence of \cite{BMR1} and another from
 perverse sheaves on $\Fl$ via the equivalence of \cite{AB}.
An indication of such a compatibility can be (and has been)
found by unraveling
 logical connections between the works of G.~Lusztig.
 However, we do not claim to have arrived at a conceptual
explanation of this coincidence.

A possible conceptual approach to the material presented in this paper is via
the local geometric Langlands duality formalism.
Recall \cite{LocGeomLang,Fr}
 that the latter theory seeks to attach to
a (geometric) local Langlands parameter a certain triangulated
category, a categorification of a
representation of a $p$-adic group attached to the Langlands
parameter by the classical local Langlands conjectures.
 According to \cite{LocGeomLang}
this triangulated category should arise as the derived category of an
abelian
category. That abelian category can conjecturally
be identified with the category of modules over an affine Lie algebra
at the critical level  with a fixed central character.
We propose the category of modules over the above algebra $\Az$ with a fixed
central character as another construction for the so-called category of
{\em Iwahori equivariant objects in a local Langlands category},
see Conjecture \ref{conj_crit} (proven in \cite{Blin})
 for a concrete statement arising from
comparing  our results with that of \cite{LocGeomLang,LocGeomLang1}.

We also hope that the t-structures on the derived categories
of coherent sheaves (in particular, those on derived
categories of coherent sheaves on varieties over $\Ce$) constructed below
 are of interest from the algebro-geometric
point of view. We expect that the construction generalizes to other
symplectic resolutions of singularities (cf.
\cite{BeKa},
\cite{Ka}) and is related to Bridgeland stability conditions, see
e.g. \cite{AnnoB}.


 \medskip

 The paper is organized as follows.  In section \ref{sect1}
 we describe the affine braid group
 action on the derived categories of coherent sheaves and state existence
 and uniqueness of a t-structure characterized in terms of this action.
We refer to \cite{BR} for construction of the $\Baff$ actions and a proof of its properties.

Section \ref{sect2} presents a proof of the facts about the
t-structures.
Uniqueness is deduced directly from a categorical counterpart of the
 quadratic relations
satisfied by the action of a simple reflection $\tii{s_\al}\in \Baff$
(the action of $\tii{s_\al}$
on the
corresponding Grothendieck groups
satisfies quadratic relations
because this action of $\Z[\Baff]$ factors
through the affine Hecke algebra).
Existence is shown by reduction
 to positive characteristic, where the statement is deduced from localization in positive characteristic
 \cite{BMR1}, \cite{BMR2}.

 Sections \ref{sect3} and \ref{sect4} present parabolic versions
 of the construction of t-structures. (They are not
 needed for the proof of Lusztig's conjectures:\
sections \ref{sect5} and \ref{sect6} are logically independent of sections
 \ref{sect3} and \ref{sect4}).

 Section \ref{sect5} recalls Lusztig's conjectures  \cite{Kth2} and reduces
 them to a positivity property of a grading on the slice algebras,
stated in detail in \ref{Property bst}.
 We finish the section by showing that positivity of the grading
implies Koszul
 property of the graded algebras.

 Section \ref{sect6} proves this compatibility by relating the t-structure 
to perverse sheaves
 on affine flag variety of the Langlands dual group.
The relation between complexes of constructible sheaves on the affine flag
 varieties  and derived categories of coherent sheaves 
comes from the result of \cite{AB}.
 Once the relationship between our abelian categories and perverse 
sheaves is established,
 the desired property of the grading follows from
the purity theorem, similarly to the proof
  of Kazhdan-Lusztig conjecture.

Appendix
\ref{Involutions on homology of Springer fibers}
contains a proof of a technical statement about compatibility of the Springer
 representation of the Weyl group on cohomology of a Springer fiber with a certain involution
 on the cohomology space (the compatibility also follows from a recent
preprint \cite{Kato}). 
This is needed in analysis of the involution of the (equivariant)
 cohomology space appearing in Lusztig's formulation of his conjectures.
Appendix B by Eric Sommers establishes a property
of the central element in an $SL(2)$ subgroup of $G$, which also enters
comparison of our categorical picture with the formulas from
\cite{Kth2}.

 {\bf Acknowledgments.}
The paper is an outgrowth of the ideas conceived in 1999 during the
Special Year on Representation Theory at Princeton Institute for Advanced
Study led by George Lusztig, we are very grateful to IAS and to
Lusztig for
the inspiring atmosphere.
Also, this work is a development of ideas found in Lusztig's papers,
we are happy to use another opportunity
 to acknowledge our intellectual debt to him.

During the decade of the project's hibernation it benefited a lot from
our communication with various people, an incomplete list includes
Michel van den Bergh, Jim Humphreys, Jens Jantzen, Simon Riche,
 Eric Sommers, David Vogan.
We are
very grateful to all of them.


We thank Dmitry Kaledin who explained to one of us the proof of Proposition
\ref{Kos_conj}, and  Valery Lunts who single-handedly organized a
mathematical gathering at his Moscow country house where
that communication has taken place.

 R.B. was supported
by DARPA via AFOSR grant FA9550-08-1-0315 and NSF grants
DMS-1102434, DMS-0854764.
I.M. was supported by NSF grant DMS 0901768.

\sus{Notations and conventions}
Let $G_\Z$ be a split reductive group over $\Z$.
We work over the base ring $R=\Z[\frac{1}{h!}]$
where $h$ is the maximum of  Coxeter numbers of simple factors.
So we denote
by $G=G_R$ the base change  of $G_\Z$ to $R$
and its Lie algebra by $\fg=\fg_R$.

We will use the notation $\k$ for {\em geometric points} of $R$\ie
maps $R\to\k$ where $\k$ is an algebraically closed field.
We will use an abbreviation
FGP for the set of geometric points
of
$R$ that have finite
characteristic.
Let
$\NN\sub\fg$
be the nilpotent cone and
$\BB$ the flag variety. We denote by $\Nt=T^*\BB@>>>\NN$
the Springer resolution and by  $\gt@>>>\fg$ the
Grothendieck map.
For convenience,
we fix a nondegenerate  invariant quadratic form on $\fg$
and use it to identify $\fg$ and $\fg^*$, hence
also $\gt^*$ and $\gt$.

Let $H$ be the abstract Cartan group of $G$ with Lie algebra $\fh$.
Let $\La=X^*(H)$
be the weight lattice of $G$, $Q\subset \La$ be the
root lattice and $W$ the Weyl group.
Our choice of positive roots
is such that for a
Borel subalgebra $\fb$ with a Cartan subalgebra $\ft$, the isomorphism
$\ft\cong \fh$ determined by $\fb$ carries
roots in $\fb$ into
negative roots.
Let
$I\sub\Iaff$ be the vertices
of the Dynkin diagram
for the Langlands dual group
$\cG$
and
of the affine Dynkin diagram
for $\cG$, we consider them  as affine-linear functionals on $\fh^*$.

Set  $W_\aff=W\ltimes \La$, $\WaffC = W\ltimes Q$.
Then $\WaffC$ is a Coxeter group corresponding to the affine Dynkin
graph of the Langlands dual group $\cG$, also $\WaffC\subset \Waff$
is a normal subgroup with an abelian quotient $\Waff/\WaffC\cong
\Lambda/Q \cong \pi_1(\cG)$. Thus $\Waff$ is the extended affine
Weyl group for $\cG$.
Let $\B\subset \BaffC\sub \Baff$ denote  the
braid groups attached to $W$, $\WaffC$, $\Waff$ respectively.
Let $\Waffsc\subb\WaffC$
and
$\Baffsc\subb\BaffC$
correspond to the simply connected cover of the derived subgroup of
$G$.

Thus $\Baff$ contains  reduced expressions
 $\wti$ for $w\in W_\aff$, and also a subgroup isomorphic to $\La$
consisting of the elements $\theta_\la$,
$\la\in\La$,
such that
 $\theta_\la =\widetilde \la$ when  $\la$
 is a dominant weight.
Denote by
$\Baff^+\sub\BaffC$ the semigroup generated by lifts
$\tii{s_\al}$ of all simple reflections $s_\al$ in $\BaffC$.

We consider the categories $\Coh(X)\sub
\qcoh(X)$ of coherent and
quasicoherent sheaves on $X$.
For a noetherian scheme $Y$ we sometimes denote
$\RHom_Y\dff \RHom_{D^b(\Coh(Y))}$.
The fiber products in this paper
are taken in the category of schemes
(as opposed to fiber product of varieties with the reduced scheme structure),
unless stated otherwise; more general derived fiber product 
is discussed in section \ref{DGconv}.

For a closed subscheme $Y\sub X$ we denote by
$\Coh_Y( X)$ the category of coherent sheaves on $X$ supported
set theoretically on $Y$.
In this paper
we will consider formal neighborhood of $Y$ in $X$
only in the case when
$Y\sub X$ is a base change of an affine closed embedding $S'\to S$
for a projective morphism $X\to S$.
In this situation, by the formal neighborhood
$\hatt Y$ we will understand  the {\em scheme }
 $\hatt Y = X\times _{S} \hatt {S'}$ where $\hatt{S'}$ is the
spectrum of the $I_{S'}$-adic completion of the ring $\O(S)$;
here $I_{S'}$ denotes the ideal of $S'$.
Notice that by \cite[Th\'  eoreme 5.4.1]{EGA} $\hatt Y$ is the inductive limit of nilpotent
thickenings of $Y$ in $X$  in the category
of $S$-schemes. Also by \cite[Th\' eoreme 5.1.4]{EGA}\fttt{We
thank Michael Temkin for providing this reference.} the category 
 of coherent sheaves on $\hatt Y$
is equivalent to the category of coherent sheaves on the formal
scheme completion of $Y$ in $X$. 


For any abelian category $\CC$ we denote its Grothendieck group
by $K^0(\CC)$ and in a particular case of coherent sheaves on
a scheme $X$ or finitely generated modules over
an algebra $A$ we denote
$K(X)=K^0[\Coh(X)]$
and
$K(A)=K^0[mod^{fg}(A)]$.

The pull-back or push-forward functors on sheaves
are
understood to be
the derived
functors, and
$Hom^i(x,y)$ means $Hom(x,y[i])$.

The base changes of $\gt$ and $\Nt$ with respect to
a $\fg$-scheme $S@>>>\g$
will be denoted by
$\tii S=S\times _{\g} \gt$ and $\tii S'=S\times _{\g} \Nt$.
For a complex of
coherent sheaves $\E$ on $\gt$ (respectively, $\Nt$) we let $\E_S$
(respectively, $\E_S'$)
denote its pull-back to $\tii S$ (respectively, $\tii S'$).

\se{
t-structures
 on  cotangent bundles of flag varieties:\
statements and preliminaries
}
\label{sect1}

As stated above, our basic object is the
base change
$G=G_R$
of a split reductive group $G_\Z$ over $\Z$
to the base ring $R=\Z[\frac{1}{h!}]$,
where $h$ is the maximum of  Coxeter numbers of simple factors.

Our main goal in the first two sections is to construct a
certain t-structure
$\Tau^{ex}$
on
$D^b(\Coh(\tii\fg))$, called the {\em exotic}
t-structure. 										
The induced t-structure on $D^b(Coh(\gt_\k))$
for a field $\k$ of positive characteristic
is related to representations of the Lie algebra $\g_\k$.
In this section we
state the results on $\Tau^{ex}$
after recalling the key ingredients:\
the action of the affine braid group
on $D^b(\Coh(\tii\fg))$,
tilting
generators
in
$D^b(\Coh(\tii\fg))$
and representation theoretic t-structures.
Some proofs are postponed to later sections.

The next three
subsections are devoted to a certain action of $\Baff$
on the derived categories of (equivariant) coherent sheaves.
In \ref{act_R} we explain a basic formalism of convolutions on derived
categories of coherent sheaves available under certain flatness
assumptions, use it to define geometric action of a group on
a derived category of coherent sheaves and state
the existence of a certain geometric action of $\Baff$ on
the derived category of coherent sheaves on $\gt$, $\Nt$.

A strengthened version of this result is presented in \ref{BCact}
where existence
of a compatible collection
of geometric actions of $\Baff$ on the fiber product spaces $\gt_S$, $\Nt_S$
(under some conditions on $S$) is stated.
In fact, such compatible collections of actions arise naturally from
a more transparent structure which is a direct generalization of
the notion of a geometric
action to the case when the space
 is not necessarily flat over the base.
This generalization involves basics of DG-schemes theory. In an attempt
to make the statements more transparent we present an informal discussion
of this more general construction in \ref{DGconv}.

Thus from the formal
 point of view subsection
\ref{DGconv} and theorem
\ref{Braid group action with base R}
are not needed. We have included them in an attempt to make
the exposition more transparent.

\sus{
Geometric action of the affine braid  group
} \label{act_R}

{\bf Definition.}
By a {\em weak homomorphism} from a group to a monoidal category
we will mean a homomorphism from the group to the group of isomorphism
classes of invertible objects.
A {\em weak action} of a group on a category $\CC$ is a weak homomorphism
from the group to the monoidal category of endo-functors of $\CC$.

\medskip

Let $X$ be a finite type flat scheme  over
a Noetherian base $S$. Then the category $D^-(\qcoh(X\times _S X))$ is a
monoidal category where the monoidal structure comes from convolution:
$\F_1*\F_2
=\
pr_{13*}(pr_{12}^*(\F_1)\Lotimes pr_{23}^*(\F_2))$ where
$pr_{12}$,
$pr_{23}$,
$pr_{13}$ are the three
projections $X\times_S X\times _SX\to X\times_SX$.
This monoidal category acts on $D^-(\qcoh(X))$
by
$\F:\GG\mapsto
pr_{1*}(\F\Lotimes pr_2^*(\G))$.

By a weak {\em geometric action} of a group on $X$
over $S$,
we will understand
a weak homomorphism from the group to  $D^-(\qcoh(X\times _S X))$.

We will say that the action is {\em finite} if its image is contained in the
full subcategory $D^b(\Coh(X\times_S X))$ and the corresponding action
on the derived category of sheaves on $X$
preserves $D^b(\Coh(X))\subset D^-(\qcoh(X))$.

For a map $S'\to S$ we can base change the above structures in a straightforward
way. Namely, the pull-back functor $D^-(\qcoh(X\times_S X))\to D^-
(\qcoh(X_{S'}\times _{S'} X_{S'})$ is monoidal and the pull-back
functor $D^-(\qcoh(X))\to D^-(\qcoh(X_{S'}))$ is compatible with the
action of the monoidal categories. Thus a weak geometric action
of a group on $X$ over $S$ induces a weak geometric action of the same
group on the fiber product space $X_{S'}$ over $S'$.

\sss{
Action of $\Baff$ on $\gt$, $\Nt$ over $R$
}\label{111}
 The Weyl group  $W$ acts on
$\gt^{reg}
\dff\gt\tim_\fg\fg^{reg}$.
Let $\Gamma_w\subset \gt\times _\g \gt$ be the closure
 of the graph of the action of $w\in W$ and set
$\Gamma_w'=\Gamma_w\cap \Nt^2$.

\stheo {\em
\lab{Braid group action with base R}
 There exists  a unique finite weak
geometric
action
 of $\Baff$ on $\gt$
 (respectively, on $\Nt$)
over $R$, such that:

i) for $\la\in \La$,  $\theta_\la$
corresponds to the direct image of the line bundle
$\O_{\gt}(\la)$ (respectively, $\O_{\Nt}(\la)$) under the diagonal embedding.

 ii) for a finite simple reflection $s_\al\in W$, $\widetilde{s_\al}\in \B$
corresponds to
the structure sheaf
$\O_{\Gamma_{s_\al}}$ (respectively, $\widetilde{s_\al} \mapsto \O_{\Gamma_{s_\al}'}$).
}

The proof  appears in \cite{BR}.
 We denote the  weak
geometric action on $\gt$ by $\Baff \ni b \mm\ K_b\in
D^b[\Coh(\gt\tim_R\gt)]$.

\srem By the discussion preceding the Theorem, we also get
geometric actions of $\Baff$ on, say, $\gt_\k$, $\Nt_\k$ where $k$
is a field mapping to $R$. For applications below we need to consider
more general base changes, these are dealt with in \ref{BCact} below.

\rem
It is possible to deduce the Theorem from the results
of \cite{BMR1} which provide an action
of $\Baff$ on the derived
 category of modular representations using the reduction to
 prime characteristic techniques
of section \ref{sect2} below. This would make the
 series of \cite{BMR1}, \cite{BMR2} and the present paper
self-contained. However, this would further increase the amount of
technical details without
adding new conceptual features to the picture.
For this reason we opted for a reference to a
more satisfactory proof in \cite{BR}, see also
\cite{Ribraid} for a partial result
in this direction.


\sus{Digression: convolution operation via DG-schemes
}
\lab{DGconv}
This subsection serves the purpose of motivating the formulation
in the next Theorem \ref{Action of braid group on  base changes}
 which is a strengthening of Theorem
\ref{Braid group action with base R}.
It relies on some basic elements of the formalism of DG-schemes.
Neither that formalism nor the statements of the present subsection
will be used in the rest of the paper (except for Remark \ref{remDGf}).
See \cite{BR} for details.

Let $X\to S$ be again a morphism of finite type with $S$ Noetherian,
but let us no longer assume that it is flat.
Then one can consider the
{\em derived fiber product}\fttt{It may be more logical
to denote the fiber product by $\Rtimes$ as it can be thought of
as a right derived functor in the category of schemes, corresponding
to the left derived functor $\Lotimes$ in the category of rings
which is opposite to the category of affine schemes.}
$X^2_S=X\Ltimes _S X$,
this is
a differential graded scheme whose structure sheaf is the derived
tensor product  $\O(X)\Lotimes _{\O(S)} \O(X)$.

Definitions similar to the ones presented in \ref{act_R}
work also in this context providing
the triangulated category $DG\Coh(X^2_S) $
of coherent $\O_{X^2_S}$-modules with the convolution
monoidal structure. This monoidal category acts on
the category $D^b(Coh(X))$.

The example relevant for us is when $S=\g_R$ and
$X=\gt_R$ or $X=\Nt_R$.
Notice that in the first case one can show that
$Tor_{>0}^{\O_\g}(\O_\gt,\O_\gt)=0$ which implies that
the derived fiber product reduces to the ordinary fiber product and
$DG\Coh(\gt\Ltimes _\g \gt)\cong D^b(\Coh (\gt\times _\g \gt))$.
However, even in this case the definition of monoidal structure
can not (to our knowledge) be given without using derived schemes,
as it involves the triple fiber product
$\gt\Ltimes _\g \gt \Ltimes_\g \gt$
 where higher Tor vanishing does not hold.

Given a pair of morphisms $X\to S \to U$
we get a natural morphism $i_{U}^S:X^2_S\to X^2_U$.
It turns out that the functor of direct image $(i_{U}^S)_*$ can be equipped
with a natural monoidal structure and the action of $DG\Coh(X^2_S)$
on $D^b(\Coh(X))$ factors through $DG\Coh(X^2_U)$.

For example, we can take $X=\gt$ or $\Nt$, $S=\g$, $U=R$.
The composed map $X\to U$ is flat, so the construction of the monoidal
structure and the action in this case reduces to the more elementary case
described in Theorem
\ref{Braid group action with base R}.

The advantage of considering
the finer structure of a geometric action on
$X$
over $S$
rather than the weaker
structure of a geometric action on
$X$
over $U$
(which
in our example
happens
 to be more elementary) is the possibility to perform the base change
construction
for the base $S$.

Namely,
given a morphism $S'\to S$ consider
$X_{S'}=S'\Ltimes _S X$ and $X^2_{S'}:=
X_{S'}\Ltimes _{S'} X_{S'}\cong (X\Ltimes_S X)\Ltimes
_S S'$. Then
$D^b[\Coh(X^2_{S'})]$
is a monoidal category acting on
$D^b[\Coh(X_{S'})]$.

The functor of pull-back under the morphism $X^2_{S'}\to X^2_S$  turns out
to be monoidal while the pull-back and push-forward functors for the morphism
$X_{S'}\to X_S$ are compatible with the module category structure.
In particular a (weak) geometric action of a group $\Ga$ on $X$ over $S$
yields weak actions of $\Ga$ on $DG\Coh(X_{S'})$ for any $S'\to S$.
This way our action of $\Baff$ yields some actions considered by other
authors, see Remark \ref{other_people}.

The
geometric actions of $\Baff$ on $\gt$, $\Nt$
from Theorem
\ref{Braid group action with base R},
actually lift to
geometric actions over
$\fg$
and this provides a rich supply of interesting
base changes of the action.
Rather than spelling out the details on geometric
actions over base $\fg$
we
will here
record a
collection of
actions of
$\Baff$ on the derived categories of a  class
of {\em exact} base change varieties
and the compatibilities they enjoy.
The exactness condition on the base change $S\to \fg$ guarantees
that $\gt_{S}$ (or $\Nt_S$)
is an ordinary scheme rather than a DG-scheme.
It
excludes some examples natural from representation-theoretic
perspective, see Remark \ref{remDGf}, but
is still sufficient for our present purposes.

\sus{$\Baff$ actions on exact base changes}
\label{BCact}
\label{BCactus}

We say that a fiber product $X_1\times_Y X_2$ is {\em exact}
 if
\begin{equation}\label{Torvan}
Tor_{>0}^{\O_Y}(\O_{X_1},\O_{X_2})=0.
\end{equation}
(We also say that
the base change $X_1@>>>X$ of $X_2$ is exact.)
We let $\BC$ (respectively, $\BC'$) denote the category of affine Noetherian
$\g$ schemes
$S\to \g$ such that the base change of $\gt$ (respectively, $\Nt$)
to $S$ is exact.
We set $\Stil=S\times _{\g} \gt$ and
$\Stilp=S\times _{\g}\Nt$.
\fttt{It may be possible to treat the two cases uniformly by considering also base changes
with respect to the morphism $\gt\to \h$. We do not develop this approach here.}




\lem {\em
\label{Sl_van}
\label{Exact base change examples}
Base changes
$\tii S$, $\tii S'$,  $\Gamma_{s_\al}\times_\g S$ and
   $\Gamma_{s_\al}'\times_\g S$
are exact
for the following maps $S@>>>\fg$:\

i) $\g_\RR$ for any Noetherian $R$-scheme $\RR$;

ii)
the spectrum $\hatt \fX$ of a completion of $\O_{\g_\RR}$ at any closed
$\fX\sub \fg_\RR$,

iii) any normal  slice
$S\subset \g_\RR$ to a nilpotent orbit in $\fg_\RR$.
}

\pf Parts (i) and (ii) are clear. Part
(iii) follows from smoothness of
the conjugation map $G\tim S@>>>\fg$, which is
clear since the differential in the direction
of $G$ produces orbital directions and
the one in direction of $S$
produces the normal directions. \qed

\sss{
Action of braid group on  base changes
}\label{132}
We will use the action of  $G\tim \Gm$ on $\fg$ (and
all related objects), where $\Gm$ acts on $\fg$ by dilations.

\stheo {\em
\lab{Action of braid group on  base changes}
Let $\GG$ be a group with a fixed homomorphism to $G\tim\Gm$.

a) Let $S$ be a scheme with a $\GG$ action and
 $S@>>>\fg$ be a  $\GG$-equivariant affine map.
If it is  in $\BC$ (respectively, in $\BC '$),
then the category $D^b(\Coh^\GG(\Stil))$
(respectively, $D^b(\Coh^\GG(\Stilp))$),
carries a canonical weak
action
 of $\Baff$
such that


\bi
(i)
For a finite simple reflection $s_\al\in W$
the generator $\tii{s_\al}$ acts by convolution with
$\O_{\Gamma_{s_\al}\times_\g S}$,
respectively  $\O_{\Gamma_{s_\al}'\times_\g S}$,
{\em provided} that
the fiber product $\Gamma_{s_\al}\times_\g S$,
 respectively  $\Gamma_{s_\al}'\times_\g S$,\
is exact.
\i
(ii)
The generators
$\theta_\la$, $\la \in \La$ act by tensoring with the line bundle
$\O(\la)$. 
\ei

For a $\GG$-morphism $S_1\to S_2$ in $\BC$
(respectively, $\BC'$), the pull-back and
 push-forward functors are compatible with the $\Baff$ action.
The
change of equivariance functors for $\GG'\to\GG$
commute with the $\Baff$
action.

b)
Let $\k$ be an algebraically closed field
of characteristic zero or $p>h$ and $e\in \g^*_\k$ be a
nilpotent element.
If the group $\GG$ fixes $e$,
then the  induced action
of $\Baff$ on
$K^0(Coh_{\BBke}^\GG(\gt_\k)) = K^0(Coh_{\BBke}^\GG(\Nt_\k))=K^\GG(\BBke)$
factors through the standard action of the affine Hecke algebra \cite{Kth2}
in the following way.
\bi
(i)
For a finite simple reflection $s_\al$, the action of $\tii{s_\al}$
on the K-group of $\BBke$ is by
 \begin{equation}\label{K_of_s_al}
\tii {s_\al}
\ =\
(-v)^{-1}
T_{s_\al},
 \end{equation}
where $ T_{s_\al}$
is the action
(from \cite{Kth2}),
of the Hecke algebra on the
K-group.

\i

(ii)
For
$\la\in\La$  the action of
 $\theta_\la \in \Baff$ is compatible
 with the action of $\theta_\la$ in the affine Hecke algebra
 defined in \cite{Kth2}.
\ei

In particular, under the Chern character map\fttt{Here
by homology we mean $l$-adic homology ($l\ne \chara(\k)$), or the
classical homology with rational coefficients
if $\k=\Ce$.}
 $K^0(Coh_{\BBke}(\gt_\k))\to H_*(\BBke)$,
the action of $\B\sub\Baff$ on
the source factors through the
Springer representation of $W$ on the target.}

\srem
Notice that the statement involving $\Nt$ is {\em not}
a particular case of the statement
about $\gt$, because (in particular), the fiber product
  $\NN\tim_\fg\gt$ is not reduced, so is not isomorphic
 to $\Nt=(\NN\tim_\fg\gt)^{red}$
  as a scheme.

\sss{Examples}
\label{other_people}
(1)
When $S\subset  \g_\Ce$ is the slice
to the subregular orbit, then $\Stilp$ is the minimal
resolution of a Kleinian singularity.
The $\Baff$ action in this case is generated
by reflections at spherical objects,
(see \cite{Bridgeland} or references therein).

(2)
\lab{Extension of the Cautis-Kamnitzer action to the affine braid group}
Let us notice a relation to an action
on coherent sheaves on  affine Grassmannians.

Let $S$ be a normal slice to a nilpotent $e_n$ in $sl(2n)$, with
two equal Jordan blocks.
Then
by the result of \cite{Anno}
the restriction to $\B\sub \Baff$ of the above action on  $\tii S'$,
coincides
with the action constructed by Cautis and Kamnitzer \cite{CK}
(up to a possible
 change of normalization).


\sss{
Some properties of the action
}
Let
$b\to \Pi (b)$
denote the composed map $\Baff \to \Waff \to \Waff/\La=W$.
Let $i^\Delta:\gt\to \gt\times_\g\gt$ and $pr:\gt \to \h$ be the diagonal
embedding and the projection.

\slem {\em
\label{conjugation of h}
a)
For $\F\in D^b(Coh(\h))$ and $b\in \Baff$ we have
$$
K_b*i^\Delta_*pr^*(\F)*K_{b^{-1}}
\cong i^\Delta_*pr^*(\Pi(b)_*\F)
.$$


b) For $\al\in I$, $K_{\tii{s_\al}^{-1}}\cong \Omega^{top}_{\Gamma_{s_\al}}
\cong \O_{\Gamma_{s_\al}}\lb-\rho, - \alpha+\rho\rb$.

c)
$\tii w(\O)\cong \O$ for $w\in W$.}

\pf a) It suffices to construct the isomorphism for the generators
of $\Baff$.
These isomorphisms come
from the fact that
$K_{\theta_\lambda}$ and
$K_{\tii {s_\al}}$ are supported on the preimage under the map
$\gt\times _\g \gt\to \h\times _{\h/W} \h$ of, respectively,
the diagonal and the graph of $s_\al$.

b) is proved in \cite{Ribraid}.




c)
This reduces immediately to the case of $SL_2$ where it follows from
the description of
$\Ga_{s_\al}$ as the blow up of $\gt$ along the zero section $\BB$.
 \qed

\sus{
Certain classes of t-structures on coherent sheaves
}
\sss{
Braid positive and exotic t-structures on $T^*(G/B)$
}
\lab{Braid positive and exotic t-structures on T*(G/B)}
A t-structure
 on $D^b(\Coh(\Sti))$,
 is called {\em braid positive} if
for any
vertex  $\alpha$ of the affine Dynkin graph of the dual group
the action of $\tii{s_\alpha}\in\Baff$ is right exact.
It is called {\em exotic} if
it is braid positive and also
the functor
 of direct image to $S$ is
exact with respect to this t-structure on  $D^b(\Coh(\Sti))$ and ordinary
 t-structure on $D^b(\Coh(S))$.

\sss{
Locally free t-structures
}
\lab{Locally free t-structures}
Here we isolate a class of t-structures
which admit certain simple construction. One
advantage is that such t-structures can be pulled-back under
reasonable base changes (say, base changes that are
affine and exact, see lemma
\ref{Locally free}.a).

For a map of
Noetherian
schemes
$f:X\to S$ we will say that  a coherent sheaf $\E$ on
$X$ is  a (relative)
{\em tilting generator}
if
the functor from $D^b(\Coh(X))$ to
$D^b(\Coh( f_*\EEnd(\EE)^{op}))$
given by $\F\mapsto\ \tx{R}f_* \RHHom (\EE,\F)$ is an equivalence.
This in particular implies
that
$f_*\EEnd(\E)$ is a coherent sheaf of rings, that the
functor lands in the bounded derived category of
coherent modules,
and that
$\tx{R}f_*\HHom(\E,\E)= f_*\HHom(\E,\E)$.

If $\E\in D^b(\Coh(X))$ is a relative tilting generator,
then the tautological t-structure on the derived category
$D^b[\Coh( f_*(\EEnd(\EE)^{op}))]$
induces a t-structure
$\Tau_\E$
on $D^b(\Coh(X))$. We call it the $\E$ t-structure.
We say that a t-structure is
{\em locally free over $S$} if it is
of the  form
$\Tau_\E$
where the relative tilting generator
$\EE$ is a vector bundle. Then
$\Tau_\E$
is given by
$$
\F\in D^{\geq 0} \iff (\tx{R} f_* \RHHom)^{<0}(\E,\F)=0
\aand
\F\in D^{\leq 0} \iff  (\tx{R} f_* \RHHom)^{>0}(\E,\F)=0
.$$

If
$S$ is affine we omit ``relative''
and say that $\EE$ is a tilting generator of $D^b[\Coh(X)]$.
Then $\EE$ is a projective generator for the heart of
$\Tau_\EE$.
In particular,
two tilting generators $\E$, $\E'$ define the same t-structure iff
they are {\em equiconstituted},
where two objects $M_1, M_2$ of an additive
category are called {\em equiconstituted} if for
 $k=1,2$ we have
$M_k\cong \oplusl N_i^{\oplus d^i_k}$
for some $N_i$ and $d^i_k> 0$.

\sss{Weak generators and tilting generators}\label{vdb}
We say that an object $X$ of a triangulated category $\cD$
is a {\em weak generator} if $X^\perp=0$, i.e. if\ \
$\Hom^\bu(X,S)=0\Rightarrow S=0$.


For future reference we recall the following

\stheo {\em
\label{gen crit}
\cite[Thm 7.6]{HvdB}\fttt{
In {\em loc. cit.}
this statement is stated under the running assumption that the scheme is
of finite type over $\Ce$.
 However, the same proof works in the present generality.
}
Assume that the scheme $X$ is projective over an affine Noetherian scheme.
Then a coherent sheaf $\EE$ is a tilting generator if and only if
$\EE$ is a weak generator for $D(\qcoh(X))$ and it is
a quasi-exceptional object\ie
$\Ext^i(\EE,\EE)=0$ for $i\ne 0$.}

\sus{
Exotic t-structures and noncommutative Springer resolution
}
\label{NC_exist}

\theo {\em
\label{tstr}
\lab{Vector bundle E}
Let $S@>>>\fg$ be an exact base change
of $\tii\fg$ (resp. of $\Nt$) with
affine Noetherian $S$.

a) There exists a unique exotic t-structure
$\Tau^{ex}_S$
on
the derived category of coherent sheaves on
$\Sti$ (resp $\Sti'$).
It is
given
by:
$$
D^{\ge 0, ex}
 =\
\{\F\ ;\ \tx{R}\,pr_* (b^{-1}\F)\in D^{\ge 0} (\Coh(S))\
\forall b\in \Baff^+
\} ;
$$
$$
D^{\le 0, ex}
\ =\
\{\F\ ;\ \tx{R}\,
pr_* (b\F)\in D^{\le 0}(\Coh(S)) \ \forall b\in \Baff^+
\}
.$$

b) This t-structure is locally free over
$S$.
In fact, there exists a
$G\times \Gm$-equivariant
vector bundle $\EE$ on $\gt$ such that for any
$S$ as above, its pull-back
$\EE_S$  to $\Sti$ (resp. $\EE_S'$ to $\Sti'$),
is
a tilting generator over $S$, and the  corresponding t-structure
is the exotic structure $\Tau^{ex}$.
In particular, the  pull-back  $\EE_S$
is a  projective generator of the heart of $\Tau^{ex}_S$.

}

\pf
In proposition \ref{uni} we
check that, for $S$ as above, any exotic t-structure
satisfies
the description
from (a). This proves
uniqueness.
The existence of a
vector bundle  $\EE$
on
$\gt$
whose pull backs
produce exotic t-structures
for any $S$ as above is proved in
\ref{Existence of an exotic tilting generator}.
 \qed


We denote the heart of the
exotic t-structure
$\Tau^{ex}_S$
by
$\Cohe(\Sti)$ (resp.
$\Cohe(\Sti')$).

\rem
\lab{stabilizer Omega of the fundamental alcove}
While the definition of an exotic t-structure
involves only the {\em non-extended} affine Weyl group $\WaffC$,
the Theorem shows that the same property
--
the right exactness
of
the canonical lifts
$\wti\in\Baff$ -- also holds for
all $w$ in the extended
affine Weyl group $\Waff$.
In particular, the stabilizer $\Om$ of the fundamental alcove
in $\Waff$ acts by  $t$-exact automorphisms of $D^b(\Coh(\Sti))$
(this is an  abelian group
$\Om\cong \La/Q \cong \Waff/\WaffC\cong \pi_1(\cG)$).

\sss{
Algebras
$\Abig$
and
$\Az$
}
\label{NCSpr}
\lab{Algebras Abig and Az}
 The
exotic t-structures
described in Theorem
  \ref{tstr}
can also be recorded as follows.
For any $S$ as in the Theorem we get an associative algebra $\Abig_S\dff
End(\EE_S)^{op}$ (respectively, $\Az_S\dff End(\EE_S')^{op}$),\
together
with an equivalence $D^b[Coh(\gt_S)]\cong D^b[mod^{fg}(\Abig_S)]$
(respectively,
 $D^b[Coh(\Nt_S)]\cong D^b[mod^{fg}(\Az_S)]$),
sending $\EE_S$ to the free rank one module.
It is clear that the algebra together with the equivalence of derived
categories determines the t-structure, while the t-structure determines
the algebra uniquely up to a Morita equivalence.

According to the terminology of, say,
\cite{BoOr}, the noncommutative
$\O(\N)$-algebra $\Az=\End(\E|_\Nt)$
is a {\em noncommutative resolution}
of singularities of the singular affine algebraic variety $\N$,
while $\Abig=\End(\E)$ is a noncommutative resolution of the affinization
$\fg\tim_{\fh/W}\fh$ of $\gt$.
In view of its close relation to Springer resolution, we
call $\Az$
a
{\em noncommutative Springer resolution},
while $\Abig$ will be called a
{\em noncommutative Grothendieck resolution}
of $\fg\tim_{\fh/W}\fh$,
cf. \cite{ICM}.

For future reference we record some properties of the algebra $A$ that directly follow from the Theorem.

\slem {\em
(a)
$\Abig$ is a vector bundle and a Frobenius algebra over $\fg$.

(b)
$\Az\dff\End(\E|_\Nt)$ is
a base change $\Az\cong\Abig\Lten_{\OO(\fh)}\OO_0$
of $\Abig=\End(\E)$ in the direction of $\fh$.

(c)
Algebras associated to base changes
$S\to \fg$
are themselves
base changes in the direction of $\fg$:\ \
$\Abig_S\cong \Abig\otimes_{\O(\g)}\O(S) $\
(resp.\
$\Az_S\cong \Az\otimes_{\O(\g)} \O(S)$).
}

\pf
(a)
The sheaf of algebras
$\AA=\EEnd(\EE)$
is Frobenius
for the trace functional $tr$.
Since
$\gt$ and $\g$  are Calabi-Yau and have the same dimension,
Grothendieck duality implies that the sheaf
$A=(\gt@>\pi>>\fg)_*\AA$ is self-dual.
In particular, it is a Cohen-Macaulay sheaf, so since $\fg$ is smooth
this implies that it is a vector bundle.
Moreover, since
$\gt$ is finite and flat over the regular locus $\g_r\sub\fg$,
Frobenius structure $tr$ on $\AA$
induces a Frobenius structure on
$\Abig|_{\fg_r}$.
Now, since the complement is of codimension three,
this
extends to a Frobenius structure
on the algebra bundle $\Abig$.

For
(b) we have
$$
\Az
\ \aa{(1)}= \
\RGa(\AA|_\Nt)
\ \aa{(2)}= \
\RGa(\AA\Lten_{\OO(\fh)}\ \OO_0)
=\ \RGa(\AA)\Lten_{\OO(\fh)}\ \OO_0
\ \aa{(3)}=\
\Abig\Lten_{\OO(\fh)}\ \OO_0
,$$
here vanishing statements
(1)
and (3)  come from $\EE$ and $\EE|_\Nt$ being tilting generators,
while
(2) follows from $\Nt=\gt\tim_\fh 0$
and flatness of $\gt@>>>\fh$.

(c)
follows from base change isomorphisms
$ \RGa (\EEnd(\EE_S))\cong \O(S)\Lotimes_{\O(\g)} \Abig
$\
(and similarly for $S'$), which follow from the exactness
assumption on base change to $S$.
Since the space in the right hand side belongs to $D^{\leq 0}$
while the space in the left hand side
lies in $D^{\geq 0}$, both in fact lie in homological degree zero
and the above isomorphisms hold.

\sss{
Remark on DG-version of the theorem
}
\label{remDGf}
We have required exactness of base change to avoid dealing with DG-schemes.
Using some basic elements of that formalism one can derive the following
generalization of
Theorems
\ref{Action of braid group on  base changes} and  \ref{tstr}.
Let $S$ be an arbitrary affine Noetherian
 scheme equipped with a morphism
$S\to \fg$.
First,
$\Baff$ acts naturally on the triangulated category
 $DG\Coh(\gt\Ltimes   _{\g}S)$
of differential graded coherent sheaves
on the derived fiber product  $DG\Coh(\gt\Ltimes   _{\g}S)$.
Then this allows us
to extend the definition of an exotic t-structure to that context.

Finally,
let $\E$ be as in the theorem and
$\Abig=End(\EE)$ as above, then
we have
an equivalence of triangulated
categories
(the tensor product
$\Abig\otimes _{\O_{\g}}\O(S)$
does not have to be derived since  $\Abig$ is flat over
$\O(\g)$ by Lemma
\ref{Algebras Abig and Az}) :
$$
DG\Coh(\gt\Ltimes   _{\g}S)
\ \cong \
D^b(mod^{fg}[\Abig\otimes _{\O_{\g}}\O(S)])
.$$
The  t-structure on
$DG\Coh(\gt\Ltimes   _{\g}S)$ corresponding to the tautological
t-structure on  $D^b(mod^{fg}[\Abig\otimes _{\O_{\g}}\O(S)])$ is the unique
exotic t-structure.

In particular, when $S=\{e\}$ is a $\k$ point of $\g$ where $\k$ is an
 algebraically closed field of characteristic $p>h$, then
the category
$mod^{fg}[\Abig\otimes _{\O_{\g}}\O_e]$
is identified with a regular
block in the category of $U_{\k,e}$-modules, where the subscript
denotes reduction of $U_\k$ by the corresponding
maximal ideal in the Frobenius center.
The category
$
DG\Coh(\gt\Ltimes   _{\g}e)
$ of coherent sheaves on the
DG Springer fiber
is studied in \cite{Simonchik}.

\sus{
Representation theoretic
t-structures on derived categories of coherent sheaves
}
\label{RTtstr}
We now record a particular case of Theorem \ref{tstr}
that follows from the results of \cite{BMR1}, \cite{BMR2}.
In the next section
we will deduce the general case
from this particular case.

Fix
$\k\in \FGP$ and a  nilpotent $e\in \N(\k)$.

\sss{
The center of $
U\fg_\k
$
}
\lab{The center of Ug}
The description of the center of enveloping algebra
in characteristic $p>h$ is
$
Z(U\fg_\k)
\cong
\OO_{ \fg_\k^*\tw\tim_{\fh_\k^*/W\tw}\fh_\k^*/W }
$,
where
$X\tw$ denotes the Frobenius twist of a $\k$-scheme $X$
and the
map
$\fh_\k^*\to
\fh_\k^*\tw$ is the Artin-Schreier map
\cite{BMR1}.
When  $X$ is one of $\fg_\k,\fh_\k$
we can use the canonical
${\Bbb F}_p$-rational structure to identify $\k$-scheme
$X\tw$ with $X$.
This gives
isomorphism
$
\fg_\k^*\tw\tim_{\fh_\k^*/W\tw}\fh_\k^*/W
\cong\
\fg_\k^*\tim_{\fh_\k^*/W}\fh_\k^*/W
\cong\
\fg_\k\tim_{\fh_\k/W}\fh_\k/W
$.


A compatible pair of $e\in\fg_\k$ and $\la\in \fh_\k$,
gives a central character
of $U_\k$ which we can then impose
on $U_\k$ or $mod(U_\k)$.
We denote
$U^\la_{\k,e}=U_\k\ten_{Z(U_\k)}\k_{\la,e}$,
while
$U_{\k,\hatt  e}^{\hatt \lambda}$
is the completion of
$\U_\k$
at $(\la,e)$
and
$mod^\la_e(U_\k)$ is the category of modules
with generalized character $(\la,e)$.
Similarly, we get $U_\k^\la$ or $U_{\k,e}$ by
imposing a central character condition for one 
of the two central subalgebras.
 This may be combined in other ways
to get objects like $mod^\la(U_{\k,e})$ etc.

\sss{
Representation theoretic t-structures
}
\lab{Representation theoretic t-structures}
By  the main result of
\cite{BMR1},  for integral regular
$\la$
we have canonical equivalence of categories of $\g$-modules
and coherent sheaves\fttt{
A priori such equivalences require
a choice
of a splitting bundle for certain Azumaya algebra,
by "canonical  equivalences''
we mean that
we use the {\em standard}  splitting bundle from
\cite[Remark 5.2.2.2]{BMR1}.
Also, we are suppressing Frobenius twist $X\tw$
from the notation using identifications
$X\tw\cong X$ that are available when $X$ is defined over the prime
subfield.
}
$$
\CD
D^b(\Coh(\hBBke))
@>>\cong>
D^b[mod^{fg}(\U^{\hatt\la}_{\k,\hatt e})]
\\
@A{\sub}AA
@A{\sub}AA
\\
D^b(\Coh_\BBke(\tii\fg_\k))
@>>\cong>
D^b[mod^{fg,\la}_{ e}(\U_\k)]
\endCD
\aand
\CD
D^b(\Coh(\hBBkep))
@>>\cong>
D^b[mod^{fg}(\U^\la_{\k,\hatt e})]
\\
@A{\sub}AA
@A{\sub}AA
\\
D^b( \Coh_\BBke(\Nt_\k))
@>>\cong>
D^b[mod^{fg}_{ e}(\U^\la_\k)]
\endCD
.$$
(Recall that the index
$\BBke$ refers to sheaves set-theoretically supported on
$\BBke$.)
Here, the second line  is Theorem 5.4.1 in \cite{BMR1}, the
first line is
stated in the footnote on the same page.

These equivalences provide each of the derived categories
of coherent sheaves with a t-structure -- the image of the tautological
t-structure on the derived category of modules.
According to
lemma  6.1.2.a in \cite{BMR1}, this t-structure depends only on the
alcove to which
$\frac{\la+\rho}{p}$ belongs, not  on
$\la $ itself.

We call the t-structure
obtained from $\la$ such that
$\frac{\la+\rho}{p}$ is
in the fundamental
alcove (e.g. $\la=0$),  the
{\em representation theoretic t-structure}
on the derived category of coherent sheaves
(RT t-structure for short).

Here by an {\em alcove} we mean a
 connected component of the complement to the affine coroot
hyperplanes
$H_{\ch\al,n}=\{\la | \langle \ch\alpha,\la\rangle
=n\},$
 in the dual space $\h^*_\RE$
to the real Cartan algebra
$\h_\RE$, where $\ch\al$ runs over the set of coroots and $n\in\Z $.
The {\em fundamental
alcove} $\fA_0$ is the locus of points where all positive coroots take
values between zero and one. Let $\Alc$ be the set of alcoves.

\theo
\label{RTplu}
\lab{RT t-structure is exotic}
{\em
For any $\k\in FGP$ and $ e\in\NN(\k)$
the RT t-structure on  $D^b[\Coh_\BBke(\gt_\k)]$,
 $D^b[\Coh_\BBke(\Nt_\k)]$ is exotic.
Therefore, for $\la\in\La$ such that $\fra{\la+\rho}{p}\in \fA_0$,
there are canonical
equivalences
of categories
$$
mod^{fg}(U^{\hatt\la}_{\hatt e})
\cong\Cohe(\hBBke)
\aand
mod^{fg}(U^\la_{\hatt e})
\cong\Cohe(\hBBkep)
.$$
}

The proof is based on

\pro
\lab{Baff equivariance of equivalence of reps and cohs}
{\em The equivalence of Theorem 5.4.1 of \cite{BMR1} is compatible with
the $\Baff$ action.}

{\em Proof\  }
follows from \cite[Section 5]{Ribraid}. \qed

\sss{
Proof
of Theorem \ref{RTplu}
}
In \cite[2.2.1]{BMR2}
the action of
Coxeter generators
$\tii {s_\al}$
of $\Baff$ (also denoted by $\tx{\bf I}^*_\al$
in {\em loc. cit.})
is defined through
a canonical distinguished triangle
$$
M\to R_\al(M)\to \tii {s_\al}(M),
\ \ \
M\in
D^b(mod^{0,fg}(\U_\k))
,$$
where $R_\al$ is the
so called {\em reflection functor}. Thus exactness of $R_\al$ implies
that $\tii {s_\al}$ acts by right exact functors.
Also, we have a commutative diagram
\cite[Lemma 2.2.5]{BMR2}:
$$
\begin{CD}
D^b(mod^{fg}_ e(U_\k^{0} )
@>>>
D^b[\Coh_\BBke(\Nt_\k)]
\\
@V{T_{0}^{-\rho}}VV
@VV{\RGa}V
\\
D^b(mod^{fg}_ e(U_\k^{-\rho}))
@>>>
D^b[\Coh_{e}(\N_\k)]
\end{CD}
$$
Here the horizontal arrows  are localization equivalences,
and $T_{0}^{-\rho}$ is the translation functor. Thus exactness of
 $T_{0}^{-\rho}$ implies that the RT t-structure satisfies the
normalization requirement in the definition of an exotic t-structure.
\qed

\sss{
Equivariant version of
representation theoretic t-structures
}
\label{163equivav}
\lab{Equivariant representation theoretic t-structures}

We will also need an equivariant version
of localization
Theorem of \cite{BMR1} and its relation to exotic t-structures.

Let
$\k,e,\la$ be as in  Theorem \ref{RTplu}
and let
$C$ be a torus with a fixed map to the centralizer of $e$ in $G$.

Recall a traditional  enhancement
of $mod(U^\la_e)$. Since $e$ vanishes on the image of $\fc=Lie(C)$
in $\g$, the action of $\fc$ on any object of
$mod(U^\la_e)$ has zero $p$-character.
The category of restricted $\fc$-modules is semi-simple
with simple objects indexed by $\fc^*(\Fp)=X^*(C)/p$; thus
every $M\in mod(U^\la_e)$ carries a canonical grading by $X^*(C)/p$.
One considers the category $mod_{gr}(U^\la_e)$ whose object is an
$U^\la_e$ module
together with a grading by $X^*(C)$. The grading should be compatible
with the natural $X^*(C)$ grading on $U^\la_e$ and the induced
$X^*(C)/p$ grading should coincide with the above canonical one.

The goal of this subsection is to describe a geometric realization
for  $mod_{gr}(U^\la_e)$. To simplify the statement of the derived
equivalences we need to enlarge the category (without changing
the set of irreducible objects nor the Grothendieck group).

Define the categories
$mod^{fg}(U^{\hatt{\la}}_{\hatt{e}},C)$,
 $mod^{fg}(U^{\la}_{\hatt{e}},C)$ as follows.
 An object $M$ of $mod^{fg}(U^{\hatt{\la}}_{\hatt{e}},C)$
(respectively, $mod^{fg}(U^{\la}_{\hatt{e}},C)$)
is an object of  $mod^{fg}(U^{\hatt{\la}}_{\hatt{e}})$
(respectively, $mod^{fg}(U^{\la}_{\hatt{e}})$) together
with a $C$ action (equivalently, an $X^*(C)$ grading) such that:

(i) The action map $U\to End(M)$ is $C$-equivariant.

 (ii)
Consider the two actions of $\fc$ on $M$: the
derivative $\alpha_C$
of the $C$ action
and the composition $\al_\g$ of the maps
$\fc\to \fg \to End(M)$. We require that the operator $\al_\g(x)-\al_C(x)$
is nilpotent for all $x\in \fc$.

\medskip

Notice that the actions $\al_\g$ and $\al_C$ commute; moreover,
condition (i) implies that
the difference  $\al_\g(x)-\al_C(x)$ commutes with the action of $\g$.

Also, if $M\in mod^{fg}(U^{\hatt{\la}}_{e})$ then the action $\al_\g$
is semi-simple, thus in this case conditions (i,ii) above imply
that $\al_\g=\al_C$ and $M\in  mod^{fg}_{gr}(U^{\hatt{\la}}_{\hatt{e}})$.
This applies in particular when $M$ is irreducible.

For future reference we mention also that one can consider the categories
$mod^{C,fg}(U^{\hatt{\la}}_{\hatt{e}})$, $mod^{C,fg}(U^{\la}_{\hatt{e}})$
of modules equipped with a $C$ action subject to the condition (i) above
only. A finite dimensional module $M$ in one of these categories splits as a direct sum
$M=\oplusl _{\eta\in \fc^*}M_\eta$
of generalized eigenspaces
of operators $\al_C(x)-\al_\g(x)$, $x\in \fc$;
moreover,  $M_\eta=0$ unless $\eta\in \fc^*(\Fp)$.

For a general $M\in mod^{C,fg}(U^{\hatt{\la}}_{\hatt{e}})$ the quotient  $M_n$ of $M$ by the
$n$-th power of the maximal ideal in $Z_{Fr}$ corresponding to $e$ is finite dimensional.
It is easy to see that the above decompositions for $M_n$ for different $n$ are compatible,
thus we get a decomposition of the category
\begin{equation}\label{modCspl}
mod^{C,fg}(U^{\hatt{\la}}_{\hatt{e}})=\oplusl_{\eta\in \fc^*(\Fp)}
mod^{C,fg}_\eta(U^{\hatt{\la}}_{\hatt{e}}),
\end{equation}
and similarly for $mod^{C,fg}(U^{\la}_{\hatt{e}})$.
Notice that
$mod^{C,fg}_0(U^{\hatt{\la}}_{\hatt{e}})=mod^{fg}(U^{\hatt{\la}}_{\hatt{e}},
C)$ and for $\tii \eta\in X^*(C)$ twisting the $C$-action by $\eta$
gives an equivalence
$mod^{C,fg}_{\theta}(U^{\hatt{\la}}_{\hatt{e}})\cong
mod^{C,fg}_{\theta+\eta}(U^{\hatt{\la}}_{\hatt{e}})$
where $\eta=\tii \eta\, \mod p\cdot X^*(C)$.

\theo \label{opportunity}{\em
a) There exist 
 compatible  equivalences of triangulated categories
 $$
D^b[\Coh^C(\hBBke)]
\cong
D^b[mod^{fg}(U^{\hatt{\la}}_{\hatt{e}},C)]
\aand
D^b[\Coh^C(\hBBkep)]
\cong
D^b[mod^{fg}(U^{\la}_{\hatt{e}},C)]
.$$

b) Under the functor
of forgetting the equivariant structure,
these equivariant equivalences  are compatible with the
equivalences of \cite{BMR1}
from \ref{Representation theoretic t-structures}.

c) The {\em representation theoretic t-structures}
that these equivalences define on categories of coherent sheaves
coincide with
the exotic t-structures, so we have induced equivalences
$$
mod^{fg}(U^{\hatt\la}_{\hatt e},C)
\cong\Cohe^C(\hBBke)
\aand
mod^{fg}(U^\la_{\hatt e},C)
\cong\Cohe^C(\hBBkep)
.$$}


The equivalences will be constructed in \ref{BMRequivav}. Compatibility with
forgetting the equivariance will be clear from the construction, while
compatibility with t-structures follows from compatibility with forgetting
 the equivariance.

\sss{Splitting vector bundles}
\lab{Splitting vector bundles}

Theorem
 \ref{RTplu}
has a geometric consequence.
Recall that for $\k\in \FGP$, $e\in \N_\k$ and $\la\in \La$
the Azumaya algebra coming from
$\la$-twisted differential operators
splits on the formal neighborhood
$\hBBkep$
 of the Springer fiber $\BBkep$ in $\Nt_\k$.
Let $\EE^{spl}_e(\la)$ be the splitting vector
bundle constructed in
\cite[Remark 5.2.2.2]{BMR1}
(the unramified shift of $\la$ that we use is $-\rho$).

\scor{\em
When $\frac{\la+\rho}{p}$ lies in the fundamental alcove,\fttt{
The result can be generalized
for an arbitrary alcove, see \ref{tstr_alc_sect} below.}
 the splitting vector bundle $\EE^{spl}_e(\la)$
does not depend on $p$
up to equiconstitutedness.
More precisely,
there exists a vector bundle
$\VV$ on $\gt$, defined over $R$, whose
base change to
$\hBBkep$
is equiconstituted with $\EE_e^{spl}(\la)$ for every $\k,\, e,\, \la$ as above.}

\pf
Take $\EE$ as in Theorem \ref{tstr} and $\VV=\EE^*$.
Then, in view of
Proposition \ref{EiLi},
both the base change
$\EE_{\hBBke}$
and the dual of the splitting bundle $\EE^{spl}_e(\la)^*$ are projective generators
for
$\Cohe(\hBBke)$,
thus they are
 equiconstituted.\fttt{
The reason that $\EE$ gets dualized is a difference of conventions.
For a tilting generator $\EE$
it is standard to use $\RHom(\EE,-)$ to get to modules over an algebra,
while for
a splitting bundle $V$ one uses $V\Lten -$.
The first functor produces a rank one free module over the algebra
when applied to $\E$ and the second when applied to $V^*$.
}
\qed

\sex
In particular,
 $Fr_*(\OO_{\BB_\k})$ does not depend on $p>h$ up to equiconstitutedness.
The reason is that for $e=0$ and $\la=0$
the splitting bundle for the restriction of
$\DD_{ \BB_\k}$ to $ \hBBkep$
can be chosen so that
its restriction to the zero section
$\BB_\k$ is $Fr_*(\OO_{\BB_\k})$
(cf. \cite{BMR1}).


\sus{Quantum groups,
affine Lie algebras
and
exotic sheaves}

We finish the subsection by stating two conjectures on other appearances
of the noncommutative Springer resolution
 (and therefore also of
exotic sheaves)
in representation theory.
Let $e\in\NN_\C$ and denote
$\Az_e\dff \Az_\C\ten_{\OO(\N_\C)}\OO_e$.

Let $U_\ze^{DK}$ be the {\em De Concini-Kac}
form of the quantum enveloping
algebra of $\g_\C$ at a root of unity $\ze$ of odd order $l>h$
\cite{DK}.
Recall that
 the center of
$U_\ze^{DK}$ contains a subalgebra $Z_l$, the so-called
{\em $l$-center}.
The spectrum of $Z_l$ contains the intersection of the variety of unipotent
elements in $G_\C$
with the big Bruhat cell  $B^+B^-$. We identify the varieties
of unipotent elements in $G_\C$ and nilpotent elements in $\g_\C$. Fix
$e\in \N(\C)$
such that the corresponding unipotent element lies in the big cell.

\conj  \label{conj_Uq} {\em
The category
$mod^{fg}(\Az_e)$
is equivalent to a regular
block in the category of $U_\ze^{DK}$-modules
with central character corresponding
to $e$.}

A possible way to approach this conjecture is by combining the
present techniques with localization for quantum groups
at roots of unity, \cite{BK}, \cite{Tanisa}.

\medskip

The next conjecture is motivated by the conjectures and results of
\cite{LocGeomLang}, see the Introduction.

\conj\fttt{This Conjecture has been proved in \cite{Blin}.}
  \label{conj_crit}
{\em Fix a nilpotent $G_\C$-oper
$O$ on the formal punctured disc with residue $e\in \N(\Ce)$.
The category $\AA_e$ of Iwahori-integrable modules
over the affine Lie algebra $\hatt{\g}_\Ce$ at the critical
level with central character corresponding to $O$, is
equivalent to the category of $\Az_e$-modules:
$$
\AA_e\ \cong\  mod(\Az_e).
$$}

\sus{
t-structures assigned to alcoves
}
\label{tstr_alc_sect}
In  \ref{Representation theoretic t-structures}
we used representation theory to  attach
 to each alcove a collection of
t-structures on formal neighborhoods of Springer fibers
in characteristic $p>h$.
The particular case of the fundamental alcove yields the exotic
t-structure (Theorem \ref{RTplu}).
The following generalization of Theorem \ref{tstr} shows that
all these t-structures can also be lifted to
$\gt_R$ and $\Nt_R$ and hence to zero characteristic.

 For $\fA_1, \, \fA_2\in \Alc$ we will say that $\fA_1$ {\em lies above}
 $\fA_2$ if for
any positive coroot $\check\alpha$ and $n\in \Zet$, such that the
affine hyperplane $H_{\check\al,n}=\{\la,\ |\, \langle \check
\alpha, \la\rangle =n \}$
 separates $\fA_1$ and $\fA_2$, alcove
$\fA_1$ lies above $H_{\al,n}$, while $\fA_2$ lies
below $H_{\ch\al,n}$, in the sense that  for
$\la_i\in \fA_i$ we have $\langle
\ch\al,\la_2\rangle < n < \langle \ch\al,\la_1\rangle$.

Recall the right action of $\Waff$ on the set of alcoves,
 it will be denoted by $w:\fA\mapsto \fA ^w$.

\lem
\lab{braid group elements giving transition between alcoves}
{\em a) There exists a unique map
$\Alc\times \Alc\to \BaffC\sub\Baff$,
 $(\fA_1,\fA_2)\mapsto b_{\fA_1,\fA_2}$, such that
\bi

i) $ b_{\fA_1,\fA_2}b_{\fA_2, \fA_3}=b_{\fA_1\fA_3}$
for any $\fA_1,\fA_2,\fA_3\in \Alc$.


\i
ii) $b_{\fA,\fA^w}=\tii w$\ for  $w\in \WaffC$ and
 $\fA\in\Alc$, provided that
$\fA^w$
lies above $\fA$.
\ei

b) This map satisfies:
$$
b_{\lambda+\fA_1,\lambda+\fA_2}=b_{\fA_1,\fA_2}\ \ \ {\text
{for}}\ \  \lambda\in \Lambda,$$
$$b_{\fA_0,\la+ \fA_0^{w\inv}}=\theta_\la \tii w\inv
\ \ \ {\text
{for}}\ \  w\in W,\, \la\in Q,$$
where $\lambda+\fA$ denotes the $\lambda$-shift of
 $\fA$.

c)
The element $b_{\fA_1,\fA_2}$ admits the following topological description.
Let
$\fh^*_{\C,reg}=\
\h^*_\C\setminus \cup_{\ch\al,n}\ (H_{\ch\al,n})_\C$.
Notice that a homotopy class of a path in
$\h^*_{\C,reg}$
connecting two alcoves
$\fA_1$ and $\fA_2$ in $\fh^*_\R$,
determines
an element in
$\BaffC= \pi_1(\h^*_{\C,reg}/\WaffC,\bu)$,
because
 each alcove is contractible, and alcoves are permuted transitively by
$\WaffC$, so they all give the same base point $\bu$ in
$\h^*_{\C,reg}/\WaffC$.

Let $\la\in \h^*$
be a regular dominant
weight. Then the subspace $-i\la
+\h^*_\RE$ does not intersect any of the affine coroot
hyperplanes in $\h^*_\Ce$. For two points $x\in \fA_1$, $y\in \fA_2$ consider
the path from $x$ to $y$ which is a composition of the following
three paths: $t\mapsto x- it\la$\ ($0\leq t \leq 1$),
any path from
$x-i\la$ to $y-i\la$
in $-i\la+\h_\RE^*$ and the path $t\mapsto y -i (1-t)\la$\
($0\leq t \leq 1$).
Then $b_{\fA_1,\fA_2}$ is represented by this path.}

\proof Uniqueness in (a) is clear since for any two alcoves there exists
an alcove which is above both of them. To check existence define
$b_{\fA_1,\fA_2}$ as in part c). Then property (i) is clear. To see
property (ii) it suffices to consider the case when $w=s_\al$ is a simple
reflection. Then $\tii s_\al$ is represented by the loop
which starts at the fundamental alcove $\fA_0$ and runs a half-circle
(in a complex line given by the direction $\al$)
around
the hyperplane of the affine coroot $\al$, in the
positive (counterclockwise) direction and
ending at  $s_\al(\fA_0)=\fA_0^{s_\al}$. The element $z$ of $\WaffC$
such that $z\fA_0=\fA$ sends $\fA_0^{s_\al}$ to $\fA^{s_\al}$ while
the simple affine coroot is sent to an affine linear functional taking
positive values on $\fA$ and negative values on $\fA^{s_\al}$.
Thus the two loops are manifestly homotopic.

b) The first property in (b) follows from uniqueness in (a), as translation
by $\lambda$ commutes with the right action of $\Waff$ and preserves the
partial order on alcoves. To check the second one, let us first consider
the case when either $\lambda=0$ or $w=1$. If $\lambda=0$ the statement
follows from (ii) as $\fA_0$ lies above $w(\fA_0)=\fA_0^w$ for $w\in W$.
When $w=1$ and $\lambda$ is dominant, then $\theta_\la=\tii \la$
and $\la+\fA_0$
lies above $\fA_0$, so the claim follows from (a,ii). Then the case
$w=1$ and arbitrary $\lambda$ follows from the first property in b). Finally,
the general case follows from (a,i) since
$b_{\la+\fA_0,\la+\fA_0^{w^{-1}}}=b_{\fA_0,\fA_0^{w^{-1}}}=\tii w^{-1}$ by the
first property in (b). $\square$

\sex
$
b_{\fA_0,-\fA_0}
=
b_{\fA_0,\fA_0^{w_0\inv} }
=\tii{w_0}\inv
$, hence
$
b_{-\fA_0,\fA_0}
=
(b_{\fA_0,\fA_0^{w_0\inv}})^{-1}
=\tii{w_0}
$.

\sss{}
\label{tstr_alcove}
\label{t-structures and alcoves}
The first part of the next Theorem is a reformulation
of Theorem \ref{tstr}.

\stheo (cf. \cite[2.1.5]{ICM}) {\em
a) Let $X=\gt$ and $S\in \BC$ or let $X=\Nt$ and $S\in \BC'$.

 There is a unique map from $\Alc$ to the set of
t-structures on $
D^b[\Coh(S\tim_\fg X)]
,$
$\fA\mapsto \Tau_\fA^{S,X}$
such that

\ben
 (Normalization)
 The derived global sections  functor $\RGa$ is $t$-exact with
respect to the t-structure  $\Tau_{\fA_0}^{S,X}$
corresponding to the fundamental alcove $\fA_0$.

\i (Compatibility with the braid action)
 The action of the element $b_{\fA_1,\fA_2}$ sends
$\Tau_{\fA_1}^{S,X}$ to $\Tau_{\fA_2}^{S,X}$.

\i (Monotonicity)
 If $\fA_1$ lies above $\fA_2$, then $D^{\ge 0}_{\fA_1}(X)
\ \subb\ D^{\ge 0}_{\fA_2}(X)$. 

\een

b) For a fixed alcove $\fA$, the t-structures $\Tau_{\fA}^{S,X}$
are compatible with base change $S_1\to S$ in
the sense that
the direct image functor is $t$-exact.

c) For each $S$ the t-structure $\Tau_\fA^{S,X}$ is locally free and
one can choose the corresponding
tilting generators
$\E_\fA^S,\ S\in\BB\CC_X$,  as pull-backs of a
$G\times \Gm$-equivariant locally free tilting generator $\E_\fA$ on $\gt$.

d) When $\fA=\fA_0$ is the fundamental alcove then $\Tau_{\fA_0}^{S,X}$ is the
exotic t-structure of Theorem \ref{tstr}.

e)
If vector bundle $\E_\fA^S$ is a tilting generator for
$\Tau^{S,X}_\fA$ then the
dual $(\E_\fA^S)^*$
is a tilting generator for $\Tau^{S,X}_{-\fA}$.}

The proof will be given in section \ref{proof_alcove}.

\se{
Construction of exotic t-structures
}
\label{sect2}

\sus{
Action of simple reflections $\tii{s_\al}^\pmo$ on coherent sheaves
}
\lab{Action of simple reflections on coherent sheaves}
 We will need an additional property
of the action, which can be viewed as a geometric version of the
quadratic relation in the affine Hecke algebra.

We start with an elementary preliminary Lemma.
\lem \label{aff_root_conj}
{\em For any $\al\in \Iaff$,\
$\tii {s_\al}$
is conjugate in the {\em extended}
affine braid group $\Baffsc$
to some $\tii {s_\be},\ \be\in I$.}

\pf
Consider first the case when $\al\in\Iaff$ is connected
in the affine Dynkin diagram to some root
$\be\in I$ of the same length. Then $u=s_\al s_\be$ has length two
and  $us_\al=s_\al s_\be s_\al
=s_\al s_\be s_\al= s_\be u$ has length three, so
$\tii u \tii {s_\al}=\tii {s_\be} \tii u$ in $\Baff$.
Therefore,
$\tii {s_\al}=\ ^{\tii u\inv}\tii {s_\be}$.
This observation suffices in all cases  but $C_n$.

For $C_n$ the affine Dynkin diagram is a line with
two roots of equal  length $\al\in \Iaff-I$ and $\be\in I$
at the ends. The  stabilizer $\Om$  of
the fundamental alcove in the
extended affine Weyl group $\Waffsc$ acts on the affine
Dynkin diagram
and an element $\om\in\Om$ realizes the symmetry that exchanges
$\al$ and $\be$. Since the length function
vanishes on $\Om$
we find that $s_\al=\ ^{\om} s_\be$ in $\Waffsc$
lifts to
$\tii {s_\al}=\ ^{\tii \om}\tii {s_\be}$ in $\Baffsc$. \qed

We are now ready to deduce the desired property of the action.
To state it we need some notations. As before, we denote the above weak
geometric action by $\Baff \ni b \mm\ K_b\in
D^b[\Coh(\gt\tim_{\fg}\gt)]$.
For a root $\alpha \in I_{\aff}$
let $H_\alpha \subset \h^*$ be the
 hyperplane passing through $0\in \h^*$ and parallel to the affine-linear
hyperplane of $\alpha$.

\pro
\label{quadrel}
{\em
a) For every simple root $\alpha\in
I_{\aff}$ we have an exact triangle in
$D^b[\Coh(\gt\tim_{\fg}\gt)]$:
\begin{equation}
\label{KKid}
K_{\tii{s_\alpha}^{-1}}
@>a_\al>>
K_{\tii{s_\alpha}}
@>b_\al>>
\Delta_*(\OO _{\gt\times_\h H_\alpha}),
\end{equation}
 where $\Delta:\gt\to \gt^2$ is the diagonal embedding.

b) For every  $\alpha\in I_{\aff}$, and every $\F\in D=D^b(\Coh(\gt_S))$
 we have a
(canonical) isomorphism in the quotient category $D/\langle
\F\rangle$
$$\tii{s_\alpha}(\F)\cong \tii{s_\alpha}^{-1}(\F)\ \
\mod \langle \F\rangle.$$
Here $\langle
\F\rangle$ denotes the thick triangulated subcategory generated by $\F$,
i.e. the smallest full triangulated subcategory closed under direct summands
and containing $\F$.}

\pf a)  By Lemma \ref{conjugation of h}(a) validity of the claim
for a given  $\alpha\in I_{\aff}$ implies its validity
for any $\beta\in I_{\aff}$ such that $\tii {s_\beta}$ is conjugate to
$\tii {s_\al}$ in $\Baffsc$. Thus in view of Lemma
\ref{aff_root_conj} it suffices to prove the claim for $\al\in I$.
By Lemma \ref{conjugation of h}(b) we have only to check that
the divisor $D_\al\dff
\Delta(\gt\times _\h H_\al)$
in $\Gamma_{s_\al}$ satisfies
\begin{equation}\label{divisor}
\O_{\Gamma_{s_\al}}(-D_\al)
\cong \O_{\Gamma_{s_\al}}\lb-\rho,-\al+\rho\rb.
\end{equation}

It is easy to see that $D_\al$ is the scheme-theoretic intersection of
$\Ga_{s_\al}$ with the diagonal $\Delta_\gt=\Gamma_e$.
Set $Z_\al=\gt\times_{\gt_\al} \gt$, where $\gt_\al=\gt_{\P_\al}$, $\P_\al=G/P_\al$
for a maximal parabolic $P_\al$ of type $\al$. Then $\Ga_{s_\al}$ and $\Delta_\gt$
are irreducible components of $Z_\al$ and \eqref{divisor} follows from
the isomorphism of line bundles on $Z_\al$: $\JJ_{\Delta_\gt}
\cong \O_{Z_\al}\lb-\rho,-\al+\rho\rb$, where
$\JJ$ stands for the ideal sheaf. It suffices to check that the two
line bundles have isomorphic restrictions to $\Delta_\gt$
and to a fiber of the projection $pr_2:Z_\al\to \gt$ (which is isomorphic
to $\Pone$). It is easy to see that in both cases these restrictions are
isomorphic to $\O_{Z_\al}(-\al)$, $\O_{\Pone}(-1)$ respectively. Thus
\eqref{divisor} is verified.

The distinguished triangle in (a) implies that for $\F\in D$ we have
canonical distinguished triangle $\tii {s_\al}^{-1}(\F)\to \tii {s_\al} (\F)
\to \F'$ where $\F'=\O_{D_\al}*\F$. On the other hand, the obvious exact
 sequence of coherent sheaves
$0\to \O_{\Delta_\gt}\overset{\alpha}{\To} \O_{\Delta_\gt}\to \O_{D_\al}
\to 0$ yields a
distinguished triangle
 $\F\to \F\to \F'$ which shows that $\F'\in \langle \F\rangle$.
This implies (b). \qed

\sus{
Uniqueness
}
\label{pruni}
Here we prove the following description of
the
exotic t-structure.

\pro
\label{uni}
{\em Let  $S@>>>\fg$ be an exact affine base change of $\tii\fg@>>>\fg$
(resp. of $\Nt@>>>\fg$).
If  an exotic t-structure
on $\tii S$ (resp. $\tii S'$)
exists then it is unique and given by
$$
D^{\ge 0}=\{\F\ ;\ \RGa (b^{-1}(\F))\in D^{\ge 0} (\Ab)
\ \
\forall b\in \Baff^+
\} ;
$$
$$D^{\le 0}=\{\F\ ;\ \RGa (b( \F))\in D^{\le 0}(\Ab)
\ \
\forall b\in
\Baff^+
\} .
$$
}

\sss{
}
Let $\A=D^{\geq 0}\cap D^{\leq 0}$ be the heart of
an exotic t-structure $\Tau$
and let
$H^i_\Tau:D\to \A$ be the corresponding cohomology functors.
Recall that
we denote for
$w\in \Waff$
by $\tii w\in\Baff$ its canonical lift.

\slem
\label{examo}
 {\em For $M\in \A$ let $\A_M\subset \A$ be the Serre subcategory
generated by $M$. Then for any $\al\in \Iaff$
$$
H^i_\Tau(\tii{s_\alpha}M), \ H^i_\Tau(\tii{s_\alpha}^{-1}M) \in
\A_M\foor i\ne 0, \aand
H^0_\Tau(\tii{s_\alpha}M)\cong H^0_\Tau(\tii{s_\alpha}^{-1}M)
\ \mod \A_M
.$$
}


\pf
 Proposition \ref{quadrel}(b) implies that
$H^i_\Tau(\tii{s_\alpha}M)\cong H^i_\Tau(\tii{s_\alpha}^{-1}M)
\ \mod \A_M$ for all $i$. On the other hand, by
the definition of braid positivity we have $H^i_\Tau(\tii{s_\alpha}M)=0$
for $i>0$ while $H^i_\Tau(\tii{s_\alpha}^{-1}M)=0$ for $i<0$.
\qed

\cor
\label{Ai}
{\em Set $D_{-1}=D$ and
let $D_0\subset D$ be the full subcategory of objects $\F$
such that $\RGa(\F)=0$.
For $i>0$ define inductively
a full triangulated subcategory $D_i\subset D$ by:
$$
D_i=\{\F\in D_{i-1} \ ;\ \tii{s_\alpha}(\F)\in D_{i-1}\ \forall \alpha
\in \Sigma_{\aff} \}
.$$
Set $\A_i=D_i\cap \A,\ i\ge 0$, then we have

a) $\A_i$ is a Serre abelian subcategory in $\A$.

b) Any exotic
t-structure $\Tau$ induces a bounded t-structure on $D_i$,
whose heart is $\A_i$.

c) For $i>0$ and any $\alpha \in \Sigma_{\aff}$,
the composition of $\tii{s_\alpha}$
with the projection to $D/\langle\A_i\rangle$ sends $\A_i$ to
$\A_{i-1}/\A_i\subset D_{i-1}/D_i$;
it induces an exact functor $\A_i\to \A_{i-1}/\A_i$.}


\pf
We prove the statements together by induction.
Assume they are known for  $i\leq i_0$. Validity of statement (c)
for $i\leq i_0$
 implies that for any $\alpha_1, \dots, \alpha_{i_0+1}$
the functor
$$
\F\mapsto \tx{R}^k\Gamma (\tii{s_{\alpha_1}} \cdots
\tii{s_{\alpha_{i_0+1}}}
(\F))
$$
restricted to $\A_{i_0}$ vanishes for $k\ne 0$, and induces an exact functor
 $\A_{i_0}\to Vect$ for $k=0$. The subcategory $\A_{i_0+1}\subset \A_{i_0}$
 is, by definition,
the intersection of the kernels of all such functors; this shows it
is a Serre abelian subcategory in $\A_{i_0}$, hence in $\A$.
 Thus (a) holds for $i=i_0+1$.

Moreover, we see that for $\F\in D_{i_0}$ vanishing of
$\tx{R}^\bu\Gamma (\tii{s_{\alpha_1}} \cdots \tii{s_{\alpha_{i_0+1}}}(\F))$
implies that for an exotic
t-structure $\Tau$
all  $\tx{R}^\bu\Gamma (
\tii{s_{\alpha_1}} \cdots \tii{s_{\alpha_{i_0+1}}}
H^k_\Tau(\F))$ vanish for  $k\in\Z$.
Thus for $\F\in D_{i_0+1}$ we have $H^k_\Tau(\F)
\in D_{i_0+1}$ for all $k$.
This shows that the truncation functors preserves $D_{i_0+1}$, i.e. (b) holds for $i=i_0+1$.

 Finally statement (c)  for $i=i_0+1$  follows from  Lemma \ref{examo}. \qed

\medskip

\sss{
Proof of  Proposition
\ref{uni}
}
Let $\Tau$ be an exotic
t-structure.
Assume that $\F\in D$, $\F\not \in D^{\leq 0}$, and let $i> 0$
be the largest integer
such that $H^i_\Tau(\F)=M\ne 0$. It suffices to show that
$\tx{R}^i\Gamma (b(\F))
\ne 0$ for some $b\in \Baff ^+$.

Lemma \ref{241}(a)
implies that $\capl_i \A_i=\{0\}$, so $M\not \in \A_d$
for some $d$; let $d$ be smallest integer with this property. If $d=0$
then $\tx{R}^i\Gamma(\F)\ne 0$, so we are done. Otherwise let
 $b= \tii{s_{\alpha_1}} \cdots \tii{s_{\alpha_d}}$ be an element such that
$b(M)\not \in D_0$. Then by Corollary \ref{Ai} we have
$\tx{R}^0\Gamma(b(M))\ne 0$. Consider the exact triangle
$$\tau_{<i}^\Tau(\F)\to \F \to M[-i],$$
and apply $b$ to it. Since $\tii{b}$
is $\Tau$-right exact, and $\RGa$
is $\Tau$-exact, we see that $\RGa(b(\tau_{<i}(\F)))\in D^{<i}(Vect)$,
thus we see that $\tx{R}^i\Gamma(b(\F))\iso
\tx{R}^i\Gamma(b(M)[-i])\ne 0$.

This proves the description of $D^{\leq 0}_{\Tau}$, the description
of $D^{>0}_{\Tau}$ is proved similarly.
\qed

\sus{
Reflection functors
$\RR_\al$ for coherent sheaves
}
\label{crefl}
Reflection
functors can be considered as a categorical counterpart of
the
idempotent of the sign representation in a Levi subalgebra
of a Hecke algebra (or the group algebra of the affine
Weyl group).
Reflection
functors on representation categories are usually defined
using translation functors which are direct summands of the functor
of tensoring by a finite dimensional representation.
In this subsection we define geometric reflection functors and show
some favorable properties they share with reflection functors
in representation theory.
In
fact, using the results of \cite{BMR2} it is not hard to check that
these functors are compatible with the usual reflection functors
for modules over the Lie algebra in positive characteristic. We
neither check this in detail nor use in the present paper; however,
the proof of Theorem \ref{RTplu} above is closely related to this fact.

For $\al\in I_{\aff}$, let
$\Xi_\alpha\in D^b[\Coh_{\gt\times_\g \gt}(\gt^2)]
$ denote
 the pull-back of the extension
$K_{\tii{s_\alpha}}^{-1}\to K_{\tii{s_\alpha}}
\to \O_{\Delta\times_\h H_\alpha}$
under the surjection $\O_\Delta \to  \O_{\Delta\times_\h H_\alpha}$,
so we have an extension
$K_{\tii{s_\alpha}}^{-1}\to \Xi_\al\to \O_{\Delta}$.
We define the
reflection functor $\RR_\al$ by
the integral kernel
$\Xi_\al
$.

\sss{
Adjoints of reflection functors
}
We first consider finite simple roots
$\alpha \in I$.
Let $P_\al\subb B$ be a minimal parabolic of type $\al$.
The  canonical projection
$\gt=G\tim_B\fb @>\pi_\al>>\gt_\al$
is generically a ramified two-sheet covering.

\slem
\lab{finite simple roots}
{\em For $\al\in I$,\  $\Xi_\al=\OO_{\gt\tim_{\gt_\al}\gt}$
and the reflection functor
$\RR_\al$
is isomorphic to the functor
$(\pi_\alpha)^*
(\pi_{\alpha})_*$. }

\pf
We will only consider the case of $sl_2$, the general case
follows by considering an associated bundle.
Notice that $\gt\tim_{\gt_\al}\gt$
has two irreducible components
$\De_\gt$ and $S_\al$ which meet transversely along
$\De_\Nt$.
Then
$\OO_{S_\al}(-\rho,-\rho)$ is the ideal of
$\De_\gt\cap S_\al$ inside $S_\al$ and
of $\De_\gt$ inside $\De_\gt\cup S_\al=
\gt\tim_{\gt_\al}\gt$.
So, one has

$$
\CD
0
@>>> \OO_{S_\al}(-\rho,-\rho)
@>>>
\OO_{\gt\tim_{\gt_\al}\gt}
@>>>
\OO_{\De_\gt}
@>>>
0
\\
@.
@V{}VV
@V{}VV
@V{}VV
@.
\\
0
@>>> \OO_{S_\al}(-\rho,-\rho)
@>>>
\OO_{S_\al}
@>>>
\OO_{\De_\Nt}
@>>>
0.
\endCD
$$
The lower line is the construction of
the exact triangle in
the Proposition
\ref{quadrel}.a),
and then the upper line says that
$\Xi_\al$ is $\OO_{\gt\tim_{\gt_\al}\gt}$.
The claim for $\RR_\al$ follows.

\scor
\label{refllem} {\em 
For any simple root $\al\in\Iaff$ :

(a)
We have two
canonical  distinguished triangles
\begin{equation}\label{Xia}
K_{\tii{s_\alpha}}^{-1}\to \Xi_\alpha\to \O_{\Delta_\gt}
\aand
\O_{\Delta_\gt}\to \Xi_\alpha \to K_{\tii{s_\alpha}}.
\end{equation}
(b)
The left and right adjoints of
$\RR_\al$ are both isomorphic to
$\RR_\al$.

(c)
$\RR_\al$ is exact
relative to an exotic t-structure.
}

\pf
The first triangle appears in the definition of $\Xi_\al$. To get the second
one it suffices, in view of
Lemma \ref{conjugation of h}(b), to show that
 \begin{equation}\label{SSS}
\SS(\Xi_\al)[-\dim \g]\cong \Xi_\al
\end{equation}
for all $\al\in \Iaff$, where $\SS$ denotes Grothendieck-Serre duality.
For $\al \in I$ isomorphism \eqref{SSS} follows from
Lemma \ref{finite simple roots}.

Furthermore, Lemma \ref{conjugation of h}(b)
 implies that the conjugation
action of $b\in\Baff$
commutes with
Serre duality.
Thus \eqref{SSS} holds in general by
 Lemma \ref{aff_root_conj}.

This proves (a).
To get (b) we use the following general fact.
 If $X$ is Gorenstein and  $\F\in
 D^b[\Coh(X\times X)]$,
  then it is not hard to show that
 the left adjoint to the functor of (left) convolution
  with $\F$ is given by (left) convolution with
  \begin{equation}\label{adjf}
  \begin{array}{ll}
  \F^{\sv}
  \dff
   \io^*
  \left( \RHHom(\F, \O_{X^2})\otimes pr_2^*K_X \right)
  =\
  \io^*
  \left(
  \SS_{X^2}(\F)\otimes pr_1^*K_X^{-1} \right)
  ;
  \end{array}
  \end{equation}
  where $\io :X^2\to X^2$ is the involution
  $\io(x,y)= (y,x)$ while
$K_X=\Omega_X^{top}[\dim X]$ is the dualizing sheaf on $X$.

Now isomorphism \eqref{SSS}  implies (b) in view of \eqref{adjf}.

Statement (c) follows from (a), since left exactness
of $\RR_\al$ follows from the first distinguished triangle, while
right exactness follows from the second one.
\qed

\srem
Notice that unlike the generators for the affine braid group action,
the geometric reflection functors do not induce a functor
on the derived category of sheaves on $\Nt$ (or varieties obtained
from the latter by base change). This is related to the fact that
the restriction of $\Xi_\alpha$ to the preimage of $0$ under the first
projection to $\h$ is not supported (scheme-theoretically) on
the preimage of $0$ under the second projection.

\sus{
Weak generators for the derived category arising from reflection functors
}
\lab{Weak generators}

The notion of a weak generator was recalled in subsection \ref{vdb}.
In the present subsection we construct
weak generators
for $D[\qcoh(\gt)]$,
which will turn out to be locally free sheaves satisfying the requirements
of Theorem \ref{tstr}.

For a finite sequence
$J=(\alpha_1, \alpha_2, \dots, \alpha_{k})$
 of elements of $\Iaff$,
 we set
$$
\Xi_{J}\dff
\RR_{\alpha_1}\cdots \RR_{\alpha_{k}} (\O_\gt).
$$
For a finite collection
$\JJ$ of  finite  sequences  let
$
\Xi_{\JJ}\dff\
\oplusl_{J\in\JJ} \Xi_{J}
$.
Notice that these are $G\tim\Gm$
equivariant by construction.

\lem\label{241}
{\em a) There exists a finite collection of elements
$b_i\in \Baff^+$, such that each of the two objects
$\oplus b_i(\O_\gt)$ and $\oplus b_i^{-1}(\O_\gt)$ is a weak generator
for
$D[\qcoh(\gt)]$.

b)
There exists a finite collection
$\JJ$ of finite sequences  of $\Sigma_{\aff}$,
 such that $\Xi_\JJ$
is a weak generator for
$D[\qcoh(\gt)]$.
}

\pf
Pick a very ample line bundle
$\O(\la)$ on $\gt$ with $\la\in\La$; thus there exists a locally
closed embedding $\gt\aa{i}\imbed {\mathbb P}^N$,
such that $\O(\la)\cong
\O_{{\mathbb P}^N}(1)\big|_{\gt}$.
It is well-known (and follows from
\cite{BePn}) that the object
$\oplusl_{i=0}^N \O(i)$,
generates
$D^b[\Coh(\Pn)]$ as a Karoubian triangulated category.
Then so does also its $\OO(-N)$ twist  $\oplusl_{i=0}^N \O(-i)$.

If $G\in
D^b[\Coh(\Pn)]$ generates
$D^b[\Coh(\Pn)]$ as a Karoubian triangulated category, we claim that
$i^*G$ is a  weak generator
for $D(\qcoh(\gt))$.
To see this, let
$\FF\in D(\qcoh(\gt))$ be a complex such that $\Ext^\bu(i^*(G),\FF)=0$.
Since
the Karoubian triangulated subcategory of
$D^b[\qcoh(\gt)]$ generated by
the restriction $i^*G$
contains all $i^*\OO_\Pn(k),\ k\in\Z$, we have
$\Ext^\bu(i^*\OO_\Pn(k),\FF)=0$. Assuming that $\FF\ne 0$ we can find
$d$ such that the $d$-th cohomology sheaf of $\FF$ does not vanish.
There exists a coherent subsheaf $\GG$ in the kernel of the differential
$\partial_d: \FF^d\to \FF^{d+1}$ which is not contained in the image
of $\partial_{d+1}$. The sheaf $\G(k)$ is generated by global sections
for  $k\gg 0$, thus we get a nonzero map $i^*\OO_\Pn(-k)[-d]\to \FF$
contradicting the above.

Thus the collection of multiples of $\la$,\
$b_i=i\cd\la\in\
\La^+ \subset \Baff^+,\
 0\le i\le N$,
 satisfies the requirement in (a).
To deduce (b) from (a) it suffices to show that for every $b\in \Baff$
there exists a finite collection of finite sequences $J_i$
in $\Sigma_{\aff}$ such that $b(\O_\gt)$ lies in the triangulated
subcategory generated by $\Xi_{J_i}$.
Let us express  $b$ as  a product $b=
\tii{s_{\alpha_1}}^{\pm 1}\cddd  \tii{s_{\alpha_d}}^{\pm k}$
with $\al_j\in\Iaff$.
 The exact triangles
\eqref{Xia} imply  that $b(\O_\gt)$ lies in the triangulated
category generated by all $\Xi_{J}$ where $J$ runs over  subsequences
of $(\alpha_1,\dots, \alpha_k)$. \qed

\sus{
Existence.
}
\label{prtstr}
\lab{Existence of exotic tilting generators}
This subsection contains a construction of exotic t-structures.

For an object
$\EE\in D^b[\Coh(\gt)]$
and
$(R@>>>\k)\in \FGP$
we denote by
$\EE_\k^{\hat 0}$ the pull-back of $\EE_\k$
to the formal neighborhood $\hatt{\BB_\k}$ of the zero section in
$\gt_\k$.



For $S\in \BB\CC$ or  $S\in \BB\CC'$ (see section
\ref{BCact} for the notation)
we will say that
a tilting generator $\EE\in Coh(\gt_S)$ (respectively, $\EE\in
Coh(\Nt_S)$) is exotic if the t-structure
$\Tau_\EE$
(notation of \ref{Locally free t-structures}) is braid
positive.

\stheo
\lab{Construction of vector bundle E}
{\em Let $\E$ be any  weak generator
$\E=\Xi_\JJ$
from
Lemma \ref{Weak generators} with
$\JJ\ni\emp$. Let $\EE_S$  (respectively, $\EE_{S}'$)
be the sheaf on  $\gt_S$ (respectively,  $\Nt_S$)
obtained from $\E$ by pull-back.

Then for any $S\in \BB\CC$ (respectively, $S\in \BB\CC'$)
$\EE_S$ (respectively, $\EE_{S}'$)
is a locally free exotic tilting generator.
}

\sss{
Reduction to the formal neighborhood of zero section in finite characteristic
}
\lab{Reduction to the formal neighborhood of zero section in finite characteristic}
The proof of the Theorem in subsection
\ref{Existence of an exotic tilting generator}
will proceed by reduction to the case of positive
characteristic. This case is treated by invoking the representation theoretic
picture.
The reduction is achieved in the following Proposition.

\spro
\lab{Reduction}
{\em Let $\EE$ be an object of $D^b(\Coh(\gt))$  containing $\OO$
as a direct summand.

a) If $\EE$ is an exotic locally free tilting generator then
 for any $S\in \BB\CC$ (respectively, $S\in \BB\CC'$) the sheaf
$\EE_S$ on $\gt_S$ (respectively, $\EE_{S}'$ on $\Nt_S$) obtained by
pull-back from $\EE$ is a locally free exotic tilting generator.

b) Assume that $\E$ is $\Gm$ equivariant and   $\E^\hz_\k$
 is an exotic locally free tilting generator
for any $(R\to \k)\in \FGP$.
 Then  $\EE$ itself is an exotic locally free tilting generator.
}

\lem
\lab{Locally free}
{\em a) If $\E \in D^b[\Coh^\Gm(\gt)]$ is such that
for any  $(R@>>>\k)\in \FGP$ the object
$\E^\hz_\k$
is a locally free sheaf,
then $\E$ is a locally free sheaf.

b)  If $\E,\F \in D^b[\Coh^\Gm(\gt)]$ are such that
$\Ext^{>0}_{\Coh(\hatt{\BB_\k})}(\E^\hz_\k, \F^\hz_\k)=0$
for any  $(R@>>>\k)\in \FGP$,
then
$\Ext^{>0}_{\Coh(\gt)}(\E,\F)=0$.
}


\pf\
(a)
Let $U$ be the maximal open subset
 such that $H^0(\E)\big|_U$ is a locally free sheaf and $H^i(\E)\big|_U=0$
for $i\ne 0$. Then $U$
 is $\Gm$-invariant, and the condition on $\E$ shows that this set contains
all closed points of the zero fiber.
Hence $U=\gt$.

 Statement (b) is obtained by applying similar considerations
to the object $\pi_*(\RHHom(\E,\F))$ where $\pi$ is the Grothendieck-Springer
map.
Here we use the formal function theorem
which shows that
 $\Ext^i(
\E^\hz_\k, \F^\hz_\k)$ is the space of sections of the pull-back of
$\tx{R}^i\pi_*
\RHHom(\E_\k,\F_\k)$ to the formal neighborhood of zero.
 \qed

\lem
\lab{positivity criterion}
{\em Let $\E$ be a tilting generator for $D^b(\Coh(X))$
where $X=\tii S$ or $X=\tii S'$ for some $S$ in $\BB\CC$
or in $\BB\CC'$.
Then $\E$ t-structure on $D^b(\Coh(X))$ is exotic iff
$\Hom^i(\E, \tii{s_\al}(\E))=0$ for $i>0$ and $\E$ contains $\O$ as a direct
summand.
}

\pf\
Assume that $\E$ t-structure is  exotic.
Since $\O$ represents the functor $\RGa$ which is exact with respect to
this t-structure, it is a direct summand of any projective generator
of the heart, in particular of $\E$
(a tilting generator is a projective generator of the heart
of the corresponding t-structure).
Also, right exactness of $\tii{s_\al}$
and projectivity of $\E$ show that
$\Hom^i(\E, \tii{s_\al}(\E))=0$ for $i>0$.
Conversely, if $\O$ is a direct summand of $\E$, then the derived
global sections functor is t-exact. Also,  if
$\Hom^{>0}
(\E, \tii{s_\al}(\E))=0$, then $\tii{s_\al}(\E)$ lies in $D^{\leq 0}$
with respect to this t-structure, which implies that
$\tii{s_\al}$ sends
$D^{\leq 0}$ to itself.
\qed

\sss{
Proof of Proposition \ref{Reduction}
}
\lab{Proof of Proposition}
a) It is obvious that the pull-back functor sends a locally free sheaf
containing $\OO$ as a direct summand to a locally free
 sheaf containing $\OO$ as a direct summand,
 while pull-back under an affine morphism sends a weak generator
of
$D(\qcoh)$ to a weak generator of $D(\qcoh)$. In view
of the characterization of tilting generators quoted in
Theorem \ref{gen crit} and criterion for
a t-structure to be exotic
from  Lemma
\ref{positivity criterion}, it remains to see
that $\Ext^{i}(\EE_S,\EE_S)=0$,
$\Ext^i(\EE_S, \tii {s_\al}(\EE_S))=0 $
for $i>0$ (or the similar equalities for $\EE_{S'}$).
The required equalities follow by the base change theorem
which
is applicable due to the Tor vanishing condition.
This proves part (a).

Assume now that $\EE$ satisfies the conditions of (b). Then $\EE$ is locally
free by Lemma \ref{Locally free}(a). By Lemma \ref{Locally free}(b)
it satisfies the above Ext vanishing conditions. Thus
$\EE$ is a tilting generator by Theorem \ref{gen crit}, and it is
an exotic tilting generator in view of
Lemma \ref{positivity criterion}. \qed

\sss{
Proof of theorem \ref{prtstr}
}
\lab{Existence of an exotic tilting generator}
Let $\JJ$ be a finite
collection of finite sequences in  $\Iaff$
such that $\Xi_\JJ $ is a weak generator for $D^b(\Coh(\gt))$
and $\JJ \owns
\emptyset$.
We have to check that $\EE=\Xi_\JJ$  satisfies the properties from the theorem
\ref{prtstr}.
We will first reduce the verification to formal neighborhoods of zero sections
over closed geometric points of positive characteristic, then
the claim will follow from translation to $\fg$-modules.

Recall from  \ref{Weak generators} that
$\E$ is
$G\times \Gm$-equivariant
by construction
and
 $\O$ is a direct summand of $\E$ since $\emptyset
\in \JJ$.
Therefore, by the Proposition
\ref{Reduction}
it  suffices to check that for
any $R@>>>\k$ in $\FGP$
the restriction $\E_\k^{\hat 0}$
of $\E$ to the formal
neighborhood of the zero section in $\gt_\k$ or $\Nt_\k$
is a locally free tilting generator for an exotic
t-structure.

This will follow once we show that
 $\E_\k^{\hat 0}$ is a projective generator for the heart of the RT
t-structure. Indeed, then
$\E_\k^{\hat 0}$
 is locally free
because another projective
generator for the same heart,
namely any splitting vector bundle for
the  Azumaya algebra of differential operators on $\BB_\k$,
 is locally free.
Also
the t-structure given by
 $\E_\k^{\hat 0}$ is the RT t-structure, but
we know that it is exotic from Theorem \ref{RTplu}.

To check that
$\E_\k^{\hat 0}$ is a projective generator for the heart of the RT
t-structure it is enough to check that it is (1) a projective object of
the heart, and (2) is a weak generator.

%
%

Statement (2) is clear since $\E$ is a weak generator for $D^b(\qcoh(\gt))$.
To check (1) it is enough to treat the case of $\gt$, the case of $\Nt$
follows because direct image under closed embeddings $\Nt\imbed
\gt$ is exact relative
to the RT t-structures, this is clear since it corresponds to the full
embedding of
 the categories
of modular representations $mod(U^\la) \to mod(U^{\hatt{\la}})$.
We now consider the case of $\gt$.

The pull-back functor under an exact base change preserves convolutions,
so each summand of $\EE^{\hatt 0}=\Xi_\JJ^{\hatt 0}$ is of the form
$
%
\RR_{\alpha_1,\k}
\cdots  \RR_{\alpha_p,\k} \widehat \O_\k
$
where $\widehat \O_\k$ is the structure sheaf of the formal
neighborhood of the zero section in $\gt_\k$,
and the functor
$\RR_{\al,\k}$ is the convolution with  the base change of
$\Xi_\al$ under $\fg_\k@>>>\fg$.

The structure sheaf $\widehat \OO_\k$
is projective for the
RT structure because  RT structure is exotic by Theorem
\ref{RTplu},  so
the functor $\RHom(\hatt\OO_\k,-)$ is exact (it can be identified with
the direct image to $\fg$).
Thus we will be done if we show that functors
$\RR_{\al,\k}$
preserve the subcategory of
 projective objects in the RT heart.

The latter property is equivalent
to the existence of a right adjoint to $\RR_{\al,\k}$
which is exact relative to the RT t-structure.
According to the  Lemma \ref{refllem}
the right adjoint of $\RR_\al$ is isomorphic to
$\RR_\al$ itself
and $\RR_\al$ can be written as an extension of
the identity functor with either
the action of $\tii {s_\al}$ or $\tii {s_\al}\inv$, the same then holds for
$\RR_{\al,\k}$.
On $\fg$-modules  $\tii {s_\al}$ is right exact and  $\tii {s_\al}\inv$
is left exact.
According to \ref{Baff equivariance of equivalence of reps and cohs}
 the
equivalence of categories of representations
and of coherent sheaves
intertwines the two
actions of $\Baff$,
hence
 we see  that $\RR_{\al,\k}$ is both left and right exact
for the RT t-structure.
\qed

\ex
For $G=SL_2$
vector bundle $\EE$ and therefore also the
algebra $\Abig$ can be described explicitly.
$\EE$ is a sum of positive multiples of
$\OO_\gt$ and $\OO_\gt(1)\dff (\gt@>>> \BB)^*\OO_\Po(1)$.
We know that $\OO_\gt$ is a summand of $\EE$ and
since
$\rho s_\al$ fixes the fundamental alcove according to Remark
\ref{stabilizer Omega of the fundamental alcove},
$(\rho s_\al)^{\sim}\ \OO_\gt=\ \OO_\gt(1)$ is also a summand.
To see that these are all indecomposable summands, it suffices
to see that
$\OO_\gt$ and $\OO_\gt(1)$
are weak generators of  $D^b\Coh(\gt)$.
This is true since
$\OO_\Po,\OO_\Po(1)$ are weak generators of
$D^b\Coh(\Po)$
and
$\gt\inj \Po\tim\fg$.
For further explicit computations of this sort (over algebraically
 closed fields) we refer to:
\cite{BMR1} for
the (sub)regular case,
and the case $e=0$ for $SL_3$; and to
\cite{Anno} for
the case when $e\in sl(2n)$ has two Jordan blocks of equal size.

\sus{Proof of Theorem \ref{tstr_alcove}}
\label{proof_alcove}
In a),
the compatibility with the affine braid group action
(axiom (2)) says that
the choice of a t-structure $\Tau^{S,X}_{\fA_0}$ determines the t-structures
$\Tau^{S,X}_\fA=b_{\fA_0,\fA}(\Tau^{S,X}_{\fA_0})$ for all $\fA$.
We need to check that this collection of t-structures satisfies
the monotonicity property (3) if and only if the
t-structure $\Tau^{S,X}_{\fA_0}$ is exotic.

Notice that for a simple reflection  we have
 $b_{\fA_0,
s_\al\fA_0}=\tii s_\al^{\pm 1}$ where the power is $+1$ exactly when
$s_\al\fA_0$ is above $\fA_0$\ie when $\al$ is not in the finite root system.
So, braid positivity of $\Tau^{S,X}_{\fA_0}$ is implied by monotonicity --
it amounts
to monotonicity applied to  pairs of alcoves $\fA_0, s_\al\fA_0$
where $s_\al$ runs over all simple reflections.
On the other hand, assume $\Tau^{S,X}_{\fA_0}$ is braid positive.
To check property (3) for the collection
 $\Tau^{S,X}_\fA=b_{\fA_0,\fA}(\Tau^{S,X}_{\fA_0})$, it suffices to consider
a pair of neighboring alcoves $\fA_2=\fA_1^{s_\al}$.
Then the automorphism $b_{\fA_0,\fA_1}^{-1}$ sends the pair of
t-structures $(\Tau^{S,X}_{\fA_1}, \Tau^{S,X}_{\fA_2})$ to the pair
$(\Tau^{S,X}_{\fA_0}, \tii s_\al(\Tau^{S,X}_{\fA_0}))$, thus the desired
relation between  $\Tau^{S,X}_{\fA_1}$ and $\Tau^{S,X}_{\fA_2}$ follows
from braid positivity.

Thus part (a) -- and at the same time part (d) --
of the Theorem follows from
existence and uniqueness of exotic t-structure
(Theorem \ref{tstr}).

Part (b) follows from compatibility of the affine braid group
action with base change
(Theorem \ref{Action of braid group on  base changes}).

For (c) notice that for any $S,X$, a projective generator
for the heart of $\Tau^{S,X}_\fA$
can be obtained by pull-back from a projective
generator for the heart of  $\Tau^{\g,\gt}_\fA$ which is of the form
$b_{\fA,\fA_0}(\Xi_{\CJ})$.
So we just need to check that $\E_\fA\dff
b_{\fA_0,\fA}(\Xi_{\CJ})\in D^b(Coh(\gt))$
is a locally free sheaf.

As in the above arguments,
it suffices to check that for any $\k\in\FGP$
the pull-back of $b_{\fA_0,\fA}(\Xi_J)$ to the formal neighborhood of
$\BB_\k$
in $\gt_\k$ is  locally free.
This property of the pull-back follows from compatibility
with localization for $\g_\k$-modules at a Harish-Chandra central character
$\la$ such that
$\frac{\la+\rho}{p}$
 lies in the alcove $\fA$, as explained
in section \ref{RTtstr}.

(e)
It remains to check that $\E_\fA^*$
is a tilting generator for $\Tau^{\g_R,\gt}_{-\fA}$.
We will prove the equivalent claim that
$b_{-\fA,\fA_0}(\E_\fA^*)$ is
a tilting generator for the exotic t-structure.

We will check that:
\
(i)
$(\E_\fA^{S})^*$
is a tilting generator for  $\Tau^{S,\gt}_{-\fA}$ where $S$ is the
 formal neighborhood of the zero section $\BB_\k$
in $\gt_\k$ for $\k\in\FGP$,
and also that
\
(ii)
$\O$ is a direct summand in  $b_{-\fA,\fA_0}(\E_\fA^*)$.
Then  Proposition
\ref{Reduction}(b) implies the desired statement.

Statement  (i) is immediate from the standard isomorphism
$(\DD_X^\LL)^{op}\cong \DD_X^{\LL^{-1}\otimes \Omega_X}$ where
$X$ is a smooth algebraic variety over a field, $\LL$ is a line bundle,
$\Omega_X$ is the line bundle of top degree forms, and
$\DD_X^\LL$
is the sheaf
of $\LL$-twisted differential operators
(cf. \cite{BMR1}). In particular,
we get
$(\DD_\BB^\la)^{op}\cong
\DD_\BB^{-2\rho-\la}$
which shows that
the dual of a splitting vector bundle for $\DD_\BB^\la$
on the formal neighborhood of the zero section is a splitting vector bundle
for $\DD_\BB^{-2\rho-\la}$. It is easy to see that
the choice of the splitting vector bundle for $\DD_\BB^\la$ as in
\cite[Remark 5.2.2.2]{BMR1}
leads to dual vector bundles for $\la$ and $-2\rho-\la$.

To check the claim that $\O$ is a direct summand in
 $b_{-\fA,\fA_0}(\E_\fA^*)$\
recall that
for
$w\in W$, $\la\in Q$ such that $\fA=\la+w\inv\fA_0$
one has
$b_{\fA_0,\fA}=\theta_\la \tii w\inv$
(Lemma
\ref{tstr_alc_sect}.iii).
Since
$-\fA
=-\la+w\inv w_0\fA_0
=-\la+(w_0w)\inv\fA_0
$, this also implies that
$b_{\fA_0,-\fA}= \theta_{-\la}\tii{w_0w}^{-1}$
and so
$b_{-\fA,\fA_0}= \tii{w_0w}\th_{\la}$.

Since  $\O$ is a summand in
$\E_{\fA_0}$,
the sheaf
$\E_\fA=b_{\fA_0,\fA}\E_{\fA_0}=
\theta_\la \tii w\inv\E_{\fA_0}$
has a summand
$\theta_\la \tii w\inv\OO\cong \OO(\la)$
(recall that  $\tii w\O\cong \O$ for $w\in W$ by Lemma
\ref{conjugation of h}(c))
and $\theta_\la(\O)\cong \O(\la)$.
Therefore,
$\O(-\la)$ is a direct summand in $\E_\fA^*$ and
$b_{-\fA,\fA_0}(\E_\fA^*)=\
\tii{w_0w}\th_{\la}(\E_\fA^*)$
has a summand
$
\tii{w_0w}\th_{\la}\ \OO(-\la)\cong
\tii{w_0w}\ \OO
\cong \OO$.
 \qed

\se{
t-structures on $\gt_\PP$ corresponding to alcoves on the wall
}
\label{sect3}

Let $\PP$ be a partial flag variety and consider the space $\gt_\PP$
of pairs of a parabolic subalgebra $\fp\in\PP$ and an element in it.
We have a map $\pi_\PP:\gt\to \gt_\PP$.

Recall that we consider the partition of
$\La_\R=\La\otimes_\Z \RE$ into alcoves,
 which are connected
components of the complement to the hyperplanes $H_{\check\al, n}:=
\{\la\ ;\ \langle \la , 
\check \al \rangle = 
n \}$
parametrized by all $\al\in\De$ and
 $n\in \Z$.
The fundamental alcove $\fA_0\sub \La_\R$ is given by
$
0<   \langle \la , 
\check \al \rangle < 
1$ for all
positive
coroots $\check \al$.

The {\em $\PP$-wall} $\WW_\PP\sub \La_\R$ is given by
$\langle \la 
 , \check \alpha \rangle
 =0$ for roots $\alpha$ in the Levi root subsystem defined by $\PP$.
By a {\em $\PP$-alcove} we will mean a connected
component of the complement in $\WW_\PP$
to those affine coroot hyperplanes which do not contain it.

Let $S\to \g$ be an affine exact base change of
$X@>>>\fg$ where $X=\gt_\PP$ (in particular, $X$ can be
$\gt=\gt_\BB$).
An example would be
a (Slodowy)
slice (the proof for $\gt$ in Lemma \ref{Sl_van}
works also for $\gt_\PP$).
The base change of
$\pi_\PP$
to $S$
is a map $\pi_{\PP}^S: \gt_S \to \gt_{\PP,S}$.

\medskip


\lem {\em 
\label{P_generator}
Let $\varpi:X\to Y$ be  a proper morphism of finite Tor dimension
and assume that $\tx{R}\varpi_*\O_X\cong \O_Y^{\oplus N}$
for some $N$. Let $\Tau_X$, $\Tau_Y$ be  t-structure on $D^b(\Coh(X))$,
$D^b(Coh(Y))$ respectively.

a) If
 $\varpi^*$ is $t$-exact, then $\Tau_Y$ is given by
$$
\FF\in D^{<0}_{\Tau_Y}
\iff
\varpi^*\FF
\in\
D^{<0}_{\Tau_X};\ \ \ \ \ \ \ \ \ \ \ \ \
\FF\in D^{>0}_{\Tau_Y}
\iff
\varpi^*\FF
\in\
D^{>0}_{\Tau_X}.$$

In particular, $\Tau_Y$ is then uniquely determined by $\Tau_X$.

The same applies with $\varpi^*$ replaced by $\varpi^!$.

b)
 If $\E$ is a projective generator for the heart of $\Tau_X$
and  
 $\varpi^!$ is $t$-exact
then  $\varpi_*(\E)$ is a
projective generator for the heart of $\Tau_Y$.

c) 
Assume that $\Tau_X$ corresponds to a tilting generator $\E_X$, and
$\Tau_Y$
corresponds to a
tilting generator $\E_Y$.
Then the functor $\varpi^!$ is $t$-exact iff $\E_Y$ is equiconstituted with
$\varpi_*(\E_X)$.

The functor $\varpi_*$ is $t$-exact iff $\varpi^*(\E_Y)$ is a direct summand
in $\E_X^{\oplus N}$ for some $N$.

}

\pf  The ``$\Rightarrow$''
implication in (a) is immediate from the
$t$-exactness assumption.
 We check the converse for $D^{<0}$, the argument
for $D^{>0}$ is similar.
If $\FF\not \in D^{<0}$ then we have
a nonzero morphism $\phi:\FF\to \GG$, $\GG\in D^{\geq 0}_{\Tau_Y}$.
The projection formula shows that
$\varpi_*\circ \varpi^*\cong Id^{\oplus N}$, thus
the map $\varpi^*(\phi):\varpi^*\FF@>>>\varpi^*\GG$ is not zero.
Since
$\varpi^*(\GG)\in D^{\geq 0}_{\Tau_X}$
this implies that
$\varpi^*(\FF)\not\in D^{<0}_{\Tau_X}$.
This proves the statement
about $\varpi^*$, the proof for $\varpi^!$ is parallel, using the
fact that $\varpi_*\circ \varpi^!\cong Id^{\oplus N}$ (this isomorphism
follows from the one for $\varpi_*\varpi^*$ by Grothendieck-Serre duality).

To prove (b) recall that a functor between abelian categories
 sends projective objects  to projective ones
provided that its right adjoint  is exact.
Thus $\varpi_*$ sends projective objects in the heart
of $\Tau_X$ to projective ones in the heart of $\Tau_Y$
if  $\varpi^!$ is $t$-exact. Also it sends
weak generators of $D(\qcoh(X))$ to weak generators of $D(\qcoh(Y))$
since $\varpi^!$ is conservative (kills no objects), which is clear
from $\varpi_*\circ \varpi^!\cong Id^{\oplus N}$.

The ``only if'' direction in the
 first statement in (c) follows from (b), while the ``if'' part
is clear from the definition of the t-structure corresponding to a tilting
generator and adjointness between $\varpi_*$ and $\varpi^!$.
The second statement
 is clear from the fact that a functor between abelian categories
is exact if and only if its
left adjoint sends a projective generator to a projective object (equivalently,
to a summand of a some power of a projective generator). \qed

\medskip

Recall from
\ref{tstr_alc_sect}
that there is a collection of locally free t-structures
$\Tau_\fA^S$ on $\Sti$ indexed by alcoves $\fA$,
and one can choose the corresponding
tilting generator
$\E_\fA^S$  as a pull-back of a
$G\times \Gm$-equivariant locally free tilting generator $\E_\fA$ on $\gt$.

\theo {\em
a)
There exists a unique
map $\fA_\PP\mm \Tau_{\fA_\PP}^S$ from the set of $\PP$-alcoves $\fA_\PP$
to the set of t-structures on $D^b(Coh(\gt_{\PP,S}))$
such that:

If  $\fA_\PP$ is a $\PP$-alcove in the closure of an alcove $\fA$
then both  of
the functors $(\pi_{\PP}^S)_*$ and
$(\pi_{\PP}^{S})^*\cong (\pi_{\PP}^{S})^!$,
between $D^b(\Coh(\gt_S))$
  and
$D^b(\Coh(\gt_{\PP,S}))$,
are exact for t-structures $\Tau_\fA^S$ and $\Tau_{\fA_\PP}^S$.

b) Each $\Tau_{\fA_\PP}^S$ is  bounded and locally free.
Moreover, for $\fA$, $\fA_\PP$ as above,
any
projective generator $\EE_\fA^S$ for $\Tau_\fA^S$ produces
a locally free  projective generator
$\E_{\fA_\PP}^S=\tx{R}(\pi_{\PP}^S)_*(\EE_\fA^S)$
for $\Tau_{\fA_\PP}^S$.

c) If  $S$ lies  over  a geometric point $\k$ of
$R$ then
$(\pi_{\PP}^S)_*$ sends any irreducible
object in the heart of $\Tau_\fA^S$
 either to zero, or to an irreducible
object in the heart of
$\Tau_{\fA_\PP}^S$.
This gives a bijection of $\Tau_\fA^S$-irreducibles
with  non-zero images and $\Tau_{\fA_\PP}^S$-irreducibles.


}            \pf
Isomorphism  $(\pi_{\PP}^{S})^*\cong (\pi_{\PP}^{S})^!$ follows from
the fact that $\gt$, $\gt_\PP$ are smooth over $R$ of the same dimension.

The direct image
$(\tii\fg@>>>\fg)_*\OO_{\tii\fg}$
 is $\OO_{\fg\tim_{\fh/W}\fh}$ and  this is a
free module of rank  $|W|$ over
$\OO_\fg$ since the same is true for
$\OO(\fh)$ as a module for $\OO(\fh/W)$. Thus Lemma \ref{P_generator} applies.

Uniqueness of a t-structure $\Tau_{\fA_\PP}^S$
for which  $(\pi_{\PP}^{S})^*$ is $t$-exact follows from Lemma
\ref{P_generator}(a).
The remaining part of statements (a,b) is equivalent,
in view of Lemma \ref{P_generator}(c), to the following statement:

$(\bu_S)$\ \ \
$\tx{R}(\pi^S_\PP)_*(\E^S_\fA)$ is a locally free tilting generator.
Moreover, $(\pi^S_\PP)^* \tx{R}(\pi^S_\PP)_*(\E^S_\fA)$ is a direct summand
in $(\E^S_\fA)^{\oplus N}$ for some $N$.

In the special case when  $S$ is the formal neighborhood $\hz_\k$
of zero  in $\fg_{\bbk}$, $\k\in \FGP$, statements (a,b) and hence
$(\bu_S)$, follow from
results of \cite{BMR2}. More precisely, the t-structure arising from
the singular localization theorem and compatible splitting bundle satisfies
the exactness properties because direct and inverse image functors
correspond to translation functors to/from the wall
\cite[Lemma 2.2.5]{BMR2},
which are well known to be exact.
A projective generator for the heart
of the t-structure can in this case be chosen to be a
splitting
bundle
on $\hz_\k$
for an Azumaya algebra $\tii\DD$ on $\gt\tw_\PP$
\cite[Remark 1.3.5]{BMR2}.
This projective generator (and hence any)  is clearly locally free.

The general case of $(\bu_S)$ follows from the above special case
by the reasoning of section \ref{sect2}.

 First, $\tx{R}(\pi^S_\PP)_*(\E^S_\fA)$
is a weak generator because $\E^S_\fA$ is a weak generator and the right
adjoint functor $(\pi^S_\PP)^!$ is conservative.

Now (b) would follow once we verify that
 $\tx{R}(\pi^S_\PP)_*(\E^S_\fA)$ is  locally free and satisfies the
Ext vanishing condition in Theorem \ref{gen crit}, while
 $\Hom(\E^S_\fA, (\pi_\PP^S)^*\tx{R}(\pi_\PP^S)_*(\E^S_\fA))$ is projective
as a module over $End(\E^S_\fA)$ (notice that
vanishing of $\Ext^i(\E^S_\fA, (\pi_\PP^S)^*
\tx{R}(\pi_\PP^S)_*(\E^S_\fA))$ for $i\ne 0$
follows from the tilting property of
$\tx{R}(\pi_\PP^S)_*(\E^S_\fA))$).
It suffices to do it in the ``absolute'' case $S=\g$. Then
local freeness and Ext vanishing
 follow from the above special case in view of Lemma \ref{Locally free}.
Similarly, projectivity holds since it holds after base change to any
 algebraically
closed field of positive characteristic and completing by the grading topology.

Finally (c) follows from the next Proposition. \qed

We keep the notations of the Theorem, fix (and drop from
notations) the alcove $\fA$, and
$S=\g$. Set $\AA=End(\EE)^{op}$, $\AA_\PP=End(\EE_\PP)^{op}$,
$\MM=\Hom(\EE,\pi_\PP^*(\EE_\PP))$. Thus
 $\MM$ is an $\AA -\AA_\PP$ bimodule.

\pro {\em
a) Under the above equivalences $D^b(Coh(\gt))\cong
D^b(mod^{fg}(\AA))$,  $D^b(Coh(\gt_\PP))\cong
D^b(mod^{fg}(\AA_\PP))$ the functor
$\tx{R}\pi_*$ is identified with the functor
$N\mapsto \MM \Lotimes _{\AA} N$.

b) The natural map $\O(\h^*)\otimes _{\OO(\h^*)^{W_L}}\AA_\PP^{op}\To
End_{\AA}(\MM)$ is an isomorphism.

}            \pf a) is obvious from the definitions. It suffices to prove
 that (b) becomes true
after
 base change to the formal neighborhood
of zero in $\g^*_\k$, $\k\in \FGP$. It is clear that validity
of the statement is independent on the choice of tilting generators
$\EE$, $\EE_\PP$ for the hearts of a given t-structure. An appropriate choice
of $\EE$, $\EE_\PP$ yields $\AA=\Gamma (\tii \DD)$, $\AA_\PP=
\Ga (\tii \DD_\PP)$, while $\MM$ is the space of sections of the bimodule
$\BB_\la^\mu$
providing the equivalence between the Azumaya algebras $\tii \DD^{\hat \la}$,
$\tii \DD^{\hat\mu}$ on $FN(\Nt)_{\gt}$
\cite[Remark 1.3.5]{BMR2}.
Then the statement follows from
\begin{multline*}
End_{\AA}(\MM)
=
\Gamma (\EEnd_{\tii \DD^{\hat \la}}(\BB_\la^\mu))
=
\Gamma (\tii \DD^{\hat\mu})^{op}
\\
\cong \Gamma (\pi_\PP^*(\tii \DD_\PP^{\hat\mu})^{op})
=
\O(\h^*)\otimes _{\O(\h^*/W_L)}\Gamma ((\tii \DD_\PP^{\hat\mu})^{op})
=
\O(\h^*)\otimes _{\O(\h^*/W_L)}\AA_\PP.
\end{multline*}

\se{
t-structures on $T^*(G/P)$ corresponding to alcoves on the wall
}
\label{sect4}

Let $P\subset G$ be a parabolic, $\PP=G/P$, $L\sub P$  be a Levi subgroup.
The projection $\BB\to \PP$ yields the maps
 $ T^*(\PP)@<pr_1<< T^*(\PP)\times_{\PP}\BB\aa{i}\imbed T^*\BB$.

Let $\La^P\subset \La$ be the sublattice of weights orthogonal to
coroots in $L$.
We have embeddings of finite index
$\La^P\subset Pic(\PP)$,
 $\La\subset Pic(\BB)$ compatible with the pull-back under projection map
$\BB\to \PP$.

(If $G$ is simply connected then both embeddings are isomorphisms .)


\sus{
t-structures in positive characteristic
}
In the next Proposition we work over a field $\bbk=\bar{\bbk}$ of
 characteristic $p>h$.
As in
\cite[1.10]{BMR2},
$\fD^\la_\PP$
 denotes the  sheaf of
crystalline
differential operators on $\PP$ with a twist $\la\in \La^P$.

\spro {\em
\label{TstarP}
For a  weight
$\la\in \La^P$
such that the element ``$(\la+\rho)\ \mod\ p$'' of $\fh^*$
is regular, consider the functor
$$
\pi^\star:\
D^b(\Coh(T^*\PP))
@>>>
D^b(\Coh(T^*\BB))
,\ \ \
\pi^\star\F= i_{*}pr_{1}^*(\F)\otimes \O(\rho)
$$
and its left adjoint
$\pi_\star\GG=(pr_1)_{*}i^* [\GG(-\rho)]
$.

a) If $\tx{R}^{>0}\Ga (\fD ^\la_\PP)=0$
then
there exists a unique t-structure on $D^b(\Coh(T^*\PP))$, such that
the functor $\pi^\star$ is $t$-exact,
 where the target is equipped
with the t-structure $\Tau_\fA$
corresponding to the alcove $\fA$ containing $\la$.

This t-structure is locally free and a $\Gm$-equivariant projective generator
of its heart is given by
$
\EE_{\PP}\dff
\pi_\star \EE_\fA
$
for any $\Gm$-equivariant
projective generator
$\EE_\fA$ of the heart of $\Tau_\fA$.

b)
If the map $U\g \to \Gamma (\fD ^\la_P)$ is surjective, then
the $t$-exact functor  $\pi^\star$
 sends irreducible
objects to irreducible ones, and it is  injective on  isomorphism classes
of irreducibles.



}            \pf
The functor $
\pi^\star\F=i_{*}pr_{1}^*(\F)\otimes \O(\rho)$ is conservative, thus its
left adjoint
$\pi_\star
$
sends a generator to a generator. Thus
for a t-structure on $D^b(Coh(T^*\PP))$
the functor
$\pi^\star$ is $t$-exact iff
$\EE_\PP=\pi_\star[ \EE_\fA]$
is a projective
 generator of its heart.
This shows uniqueness in (a)
and reduces the rest of statement (a) to showing that
 $\EE_\PP$ is a tilting vector bundle.

 As above, it suffices to check that the restriction of
$\EE_\PP$ to the formal neighborhood
 of the zero section of $T^*\PP$, is a tilting vector bundle.
Reversing the argument of the previous paragraph
 we see that it suffices to show the existence of a locally free t-structure on $D^b(Coh_\PP(T^*\PP))$
 such that  $\pi^\star$ is $t$-exact.


It is shown in
\cite[1.10]{BMR2} that under
the assumptions of (a) the derived global sections functor
$\RGa : D^b [mod^{fg}(\fD_\PP^\la)]\to
D^b[mod^{fg}(\Ga(\fD_\PP^\la))]$ is an equivalence.
Furthermore, $\fD_\PP^\la$ is an Azumaya algebra over
$T^*\PP\tw$ which is
split on the formal neighborhood of the fibers of the moment
map ${\bf \mu}_\PP:T^*\PP\tw\to (\g^*)\tw$, in particular
on the formal neighborhood of the zero fiber. Thus we get
an equivalence between
$D^b[Coh_\PP(T^*\PP)]$ and the derived category of
modules
over a certain algebra \cite[Corollary 1.0.4]{BMR2}.
In view of
\cite[Proposition 1.10.7]{BMR2}
we see that the resulting
t-structure on
$D^b[Coh_\PP(T^*\PP)]$
satisfies the desired exactness property, thus (a) is proved.

The same Proposition also implies claim (b), since it shows that the functor
between the abelian hearts  induced by the functor
$D^b(Coh_{{\bf \mu}_\PP^{-1}(e)}(T^*\PP))
@>\pi^\star>>
D^b(Coh_{\pi^{-1}(e)}(T^*\BB))
$,
can be identified with the pull-back
functor between the categories of modules
(with a fixed generalized central character)
corresponding to the ring homomorphism
$U^\la(\g) \to \Ga (\fD^\la_\PP)$.
If this ring homomorphism is surjective then the pull-back functor
sends irreducible modules to irreducible ones
and distinguishes the isomorphism classes of irreducibles.  \qed

\sus{
Lifting to characteristic zero
}
We now return to considerations over an arbitrary base.
Let $S\to \g$ be a base change exact for both
$T^*(\BB)\to\fg$
an
$T^*(\PP)\to\fg$.
Again, any Slodowy slice is an example of such $S$.
Making base change to $S$ we get
maps $i_S$, $pr_{1,S}$
and the functor
$\pi^\star_S:\
D^b[Coh(\tii S_\PP)]
@>>>
D^b[Coh(\tii S)]
$, where
$\tii S_\PP=\ S\tim_{\fg}T^*\PP$.

\stheo {\em 
There exists an integer $N>0$ (depending on the type of $G$ only), such
that the following is true provided that $N$ is invertible on $S$.

(a) Fix an alcove $\fA$. Assume that there exists a  weight
$\la\in \La^P$, such that $\frac{\la +\rho}{p}\in \fA$.

Then there exists a unique t-structure
$\Tau_\fA^{\tii S_\PP}$
on $D^b(\Coh(\Sti_P))$, such that
the functor
$\pi^\star_S:\
D^b[Coh(\tii S_\PP)]
@>>>
D^b[Coh(\tii S)]
$\
is $t$-exact,
 where the target is equipped
with the t-structure $\Tau_\fA^S$ corresponding to $\fA$.

The t-structure
$\Tau_\fA^{\tii S_\PP}$
is locally free over $S$, a projective generator
of its heart is given by
$\EE_{\PP}^S= pr_{1,S_*}i_S^* [ \EE_\fA(-\rho)]$
for any projective generator $\EE_\fA$
of the heart of $\Tau_\fA^S$.

(b) Let $\bbk$ be a geometric point of $R$
such that  the pull-back map $\OO(\g_\k) \to \Gamma (\OO(T^*(\PP_\k)))$
is surjective.
Then
the $t$-exact functor  $\pi^\star_S$
 sends irreducible
objects to irreducible ones and is
injective on isomorphism classes of irreducibles.

}            \pf
It is well known \cite{Bro}
  that for fields $\k$ of characteristic zero,
hence also for $\k$ of a sufficiently large positive
characteristic, we have
$
H^{>0}(T^*\PP_\k,\OO)=0
$.
It follows that
 the cohomology vanishing condition of Proposition \ref{TstarP}(a)
holds over such a field. Thus Proposition \ref{TstarP} shows that
statement (a) is true when $S=\g_\k^*$ where $\k$ is an algebraically closed field
of sufficiently
large positive characteristic. Then the general case follows, as in the proof of
 Theorem \ref{tstr},
by Proposition
\ref{Reduction to the formal neighborhood of zero section in finite characteristic}.
  This proves part (a).

The general case of statement (b) follows by a standard argument
from the case when $\k$ has (large) positive characteristic.
(Notice that
if the map
$\OO(\g_\k) \to \OO(T^*\PP_\k)$ is surjective for
$\k$ of characteristic
zero,
then it is surjective for $\k$ of large positive characteristic).
 In the latter case the statement follows from
  Proposition
\ref{TstarP}(b), since
surjectivity of the map
$\OO(\g_\k)  \to \OO(T^*\PP_\k)$ implies surjectivity of the map
 $U(\g_\k) \to \Gamma(\fD^\la_{\PP_\k})$. \qed

\rem
a) It is well known
that conditions of part (b) hold when
$G$ is of type $A_n$ and
$\bbk$ is of characteristic zero.

b)
The twist by $\rho$ appearing in the last Theorem
is caused by  normalizations in the definition of the braid group action.
Removing this twist would produce a twist in the preceding Theorem.

c) Notice that control on the set of primes for which the result is valid is weaker here than
in other similar results of \cite{BMR2}
 and sections \ref{sect1}, \ref{sect3}. The only reason for
this is that higher cohomology vanishing for the
sheaf $\OO(T^*\PP)$ has not been
established in general (to our knowledge)
 in positive characteristic $p$ except for indefinitely
large $p$.

\rem\label{KosRem}
It can be deduced from the results of
\cite{Simonchik} that the
construction of the present subsection is related to that
of the preceding one by  Koszul duality (see {\em loc. cit.} for details).

The matching of  combinatorial parameters is explained in the next subsection.


\sus{
Shifting the alcoves
}
The set of t-structures constructed in this subsection is indexed
by the set of alcoves $\fA$ such that $\frac{\la+\rho}{p}\in \fA$
for some  integral $\la\in \La^P$.
Let us denote this set
by $\Alc_\PP$.

\lem {\em 
\lab{AlcPP}
The set $\Alc_\PP$ is in a canonical bijection with
the set of $\PP$-alcoves. The bijection sends $\fA\in\Alc_\PP$ to the interior
of $\bar \fA\cap\ \WW_\PP$. 

}            \pf
By a {\em $\PP$-chamber} we will mean a connected component of the complement
in the $\PP$-wall
$\WW_\PP$,
of intersections with
coroot hyperplanes not containing the wall. The $\PP$-chambers
are in bijection with parabolic subalgebras in the
Langlands dual\fttt{
Of course,
these are also in bijection with
parabolic subalgebras in $\g$ with the fixed Levi. However, the argument uses coroots rather than roots,
hence the appearance of $\ch G$.
}
 algebra $\ch\g$ with a Levi $\ch\fl\sub\ch\g$ whose
semisimple
part
is given by coroots orthogonal to $\la$.

A Weyl chamber
is a connected component of the complement to coroot hyperplanes in $\h^*_\RE$.
We will say that a Weyl chamber is {\em near the $\PP$ wall} if it contains a
$\PP$-chamber
in its closure.
The set of Weyl chambers is in bijection with the set of Borel subalgebras with a fixed Cartan.
So, a Weyl chamber is near the $\PP$-wall
iff the corresponding Borel
subalgebra
$\ch\fb$ is contained
in a parabolic subalgebra with Levi
$\ch\fl$, in other words if the subspace $\ch\fl+\ch\fb\subset \ch\g$
is a subalgebra.

It suffices to show that for every  weight
$\la\in\La^P$ such that $\la+\rho$ is regular,
$\la+\rho$ lies in a Weyl chamber which is near the $\PP$ wall.
This amounts to showing that
if
$\ch\al, \ch\beta$ and $\ch\al+\ch\beta$ are coroots and
$\langle \ch\al, \la+\rho \rangle >0$,
$\langle \ch\beta, \la \rangle =0$, then
either  $\langle \ch\al+\ch\beta, \la+\rho \rangle >0$ or
$\langle \ch\al+\ch\beta, \la \rangle =0$.

It is enough to check that  $\langle \ch\al+\ch\beta, \la+\rho \rangle >0$
assuming $\ch\beta$ is a simple
negative coroot in $\ch\fl$.
For such $\ch\be$ we have   $\langle
\ch\beta, \la+\rho \rangle = \langle \ch\beta, \rho \rangle =-1$, so
 $\langle \ch\al+\ch\beta, \la+\rho \rangle \geq 0$.
However, since we have assumed that $\la+\rho$ is regular
 and $\ch\al+\ch\beta$ is  a coroot, this implies that
$\langle \ch\al+\ch\beta, \la+\rho \rangle \ne 0$.
\qed

%
%



\np
\se{
Applications to Representation Theory
}\label{sect5}

The results of section \ref{sect2} imply
 independence of the numerics of $\fg_\k$ modules
on  the characteristic $p$ of the base field $\k$ for $p\gg 0$.
This  is spelled out in
subsection \ref{Independence of p}.

In \ref{Gradings and bases in K-theory}, \ref{54}
we briefly recall Lusztig's conjectural
description of the numerical structure of the theory and
 reduce its verification
to a certain positivity
property of a grading on the category of representations.

Given the results of sections
\ref{sect1}, \ref{sect2} it is
easy to construct a family
of gradings on the category.
However, showing that this family contains
a grading which satisfies  the positivity property is more subtle.
This is established
in section \ref{sect6} by making use of the results of \cite{AB}. This
reduction relies on  $G$-equivariant versions of some of the above
constructions. This technical variation is presented in
\ref{Equivariant versions and Slodowy slices}.

To simplify notations we only treat the case of a regular block and a nilpotent
Frobenius central character, generalization to any block and an arbitrary $p$-central character
is straightforward.

\sus{
Generic
independence of  $p$
}
\label{Independence of p}
In this section  $\E$ is a  vector bundle on $\gt$ that satisfies
Theorem \ref{Vector bundle E}(b) (so up to equiconstitutedness
it is  unique and, furthermore, has the form $\Xi_\JJ$ from Theorem
\ref{Construction of vector bundle E});
while
$\Abig$ is the $R$-algebra $\End(\EE)^{op}$.

\theo {\em
\label{Moreq}
\label{flatness of A over g}
The algebra $A$ satisfies the following.
For any $\k \in \FGP$
and any
$e \in \N_\k$,
there are canonical Morita equivalences
compatible
with the action of    $\O_{\g\times _{\h/W}\h}$ :
$$
U^{\hat 0}_{\k,\hat  e}
\ \sim\
\Abig^{\hat 0}_{\k, \hat  e}
\aand
U^{0}_{\k,\hat  e}
\ \sim\
\Abig^{0}_{\k, \hat  e}
.$$
Here,
$U_{\k,\hat  e}^{\hat 0}$,
$\Abig^{\hat 0}_{\k, \hat  e}$
are completions
of
$U_\k$ and $\Abig_\k=\Abig\otimes _R \k$,
at the central ideals corresponding
to $ e$ and to $0\in \fh$
(resp.,
$\Abig^{0}_{\k, \hat  e}$,
$U_{\k,\hat  e}^{0}$
are completions
of
$U_\k^0$ and $\Abig_\k^0=\Abig_\k\ten_{\OO(\fh)}\k_0$
at the central ideals corresponding
to $ e$).

}            \pf
Recall that
$\Az\dff \Abig\ten_{\OO(\fh)}\OO_0$ is
the algebra
$\End(\E|_\Nt)^{op}$ (Lemma \ref{Algebras Abig and Az}),
and that there are canonical
equivalences
$$
mod^{fg}(U^{\hatt 0}_{\hatt e})
\cong\Cohe(\hBBke),
mod^{fg}(\Abig)\cong \Cohe(\gt),
$$
$$
mod^{fg}(U^{0}_{\hatt e})
\cong\Cohe(\hBBkep),$$
$$mod^{fg}(\Az) \cong \Cohe(\Nt)$$
compatible with the pull-back functors. 
Both times,  the first equivalence is from
Theorem \ref{RTplu},
and the second one is
the fact that $\E$ (resp. $\EE|_\Nt$) is a  tilting generator
and that the locally free t-structure it produces is the exotic t-structure.
The following
induced equivalences of categories
give  the
desired Morita equivalences of algebras:
$$
mod^{fg}(U^{\hatt 0}_{\hatt e})
\cong\Cohe(\hBBke)
\cong
mod^{fg}(\Abig_{\k,\hatt e}^{\hatt 0})
\aand
mod^{fg}(U^{0}_{\hatt e})
\cong
\Cohe(\hBBkep)
\cong
mod^{fg}(\Az_{\k,\hatt e})
.$$

We can define an
action of
$\O_{(\g\times _{\h/W} \h/W)_\k}$
on the leftmost terms using the center of $
U_\k
$.
Recall that it is isomorphic to
$
\OO_{ \fg_\k\tim_{\fh_\k/W}\fh_\k/W }
$, where map $\fh_\k@>>>\fh_\k$ is the
Artin-Schreier map (see  \ref{The center of Ug}).
However, for a {\em regular} $\la\in\La$ (say $\la=0$)
and any $e\in\N_\k$,
the completion of this center at the point
$(e,W\la)$
is canonically
isomorphic to the completion
of
$
\O_{\g_\k\times _{\h_\k/W} \h_\k}
$ (where this time $\h_\k\to\h_\k/W$ is just the quotient map)
at the
point
$(e,\la)$.
\qed


\rems (1)
Let $ e\in \NN(R)$.
Since Morita invariance is inherited by central reductions
we get also $R$-algebras
$\Abig^0_{e}$, $\Abig^\hz_{e}$, $\Abig^0_{\hat e}$
whose base change to any $\k\in FGP$
is Morita equivalent to the corresponding central reductions
$U^0_{\k,e},U^\hz_{\k,e},U^0_{\k,\hatt e}$
of the enveloping algebra.
Here $U^\hz_{\k, e}$ is the most popular version
-- its category of modules
is the principal block of the category of
$U_{\k,e}$-modules.

(2)
A similar result for representations of algebraic groups (rather than
Lie algebras) has been established by
 Andersen, Jantzen and Soergel
\cite{AJS}.
This is equivalent to the case $ e=0$ of the theorem.
\medskip



\sss{
Cartan matrices
}
Recall that the set of nilpotent conjugacy classes in $\g_\k$ does not depend
on the algebraically closed field $\k$ provided its characteristic is a good prime.

\scor {\em
There exists a finite set of primes $\Pi$
such that
for any $\k, \k'\in \FGP$
with  characteristics outside $\Pi$,
the following holds.
If
the
conjugacy classes  of $ e\in \N_\k$, $ e'\in \N_{\k'}$
correspond to each other,
then the Cartan matrices\fttt{By this we mean the matrix whose entries are
multiplicities of irreducible modules in indecomposable projective ones.}
 of $U_ e^{\hat 0}$, $U_{ e'}^{\hat 0}$
 coincide.

The same applies to $U_ e^{ 0}$,
 $U_{ e'}^{ 0}$.

}            \pf
For a finite rank $R$-algebra $\AA$ independence of the Cartan matrix
of $\AA_\k$ on $\k$,
for an algebraically closed field $\k$ of sufficiently large
 characteristic,
is well known (cf. also Lemma \ref{Cartan_matr_general}
below). Thus the claim follows from part (1) of the the previous Remark. \qed

\sss{
Lifting irreducible and projective objects
}
\lab{Lifting irreducible and projective objects}
In the remainder of the subsection we will
enhance the last Corollary
to a geometric statement.


A  {\rm quasifinite domain} over a commutative ring $R$ is
a
commutative ring $\RR$ over $R$
for which there exists $R'\subb R$ such that
$R'$ is a finite domain over $R$ and $\RR$ is a finite localization of $R'$.

\spro {\em
\label{EiLi}
For any $e\in\NN(R)$
there exists a quasifinite
domain $\RR$ over $R$,
such that
there exist
\bi
a finite set $I_ e$,
\i
a collection
$\E_i$, $i\in I_ e$, of locally free sheaves on
$\hatt{\BB_{\RR, e}}$,
and
\i
a collection
of complexes of coherent sheaves
$\LL_i \in  D^b[
\Coh(\widehat {\BB_{\RR, e}}')
]$, $i\in I_ e$;
\ei
with the following properties.

Let $\la\in\La$ be such that
$\fra{\la+\rho}{p}$ is in the fundamental alcove.
Then for every
 characteristic $p$ geometric point
$\k$ of $\RR$,
the set of isomorphism classes of irreducible
 $U_{\k, e}^\la$
 modules
is
canonically
parametrized
by  $I_ e$, we denote this $I_ e\ni i\mapsto L_{\k,i}\in
Irr(U_{\k, e}^\la)$.
This parametrization is such that
the equivalence of \cite{BMR1}
sends:
\bi
(I) irreducible $L_{\k,i}$ to
$(\LL_i)_\k\dff\LL_i\Lten_R\k\in\
D^b[\Coh(\hBBkep)]
$;


\i
(P$_{1}$)
the projective cover of $L_{\k,i}$ over $U_{\k,\hat  e}^{\hat \la}$
to  $(\E_i)_\k
\dff\E_i\ten_R\k
\in\
D^b[\Coh(\hBBke)]
$;

\i
(P$_{2}$) the projective cover of $L_{\k,i}$ over
$U_{\k,\hat  e}^{\la}$
to
$
(\E_i)_\k\big|_\BBkep
\in\
D^b[\Coh(\hBBkep)]
$;
\i
(P$_{3}$)
the projective cover of $L_{\k,i}$ over $U_{\k,  e}^{\hat \la}$
to  $\E_i\Lten _{\O_\g} \k_ e
\in\
D^b[\Coh(\hBBke)]
$.
\ei

}            \pf The first three claims  are immediate from the
Theorem
\ref{Moreq},
 together with the next general Lemma. Part
(P$_{3}$)  follows from
part (P$_{1}$)  because
the localization equivalence is compatible with derived tensor product
over the center, while the enveloping algebra $U$
is flat over its Frobenius center.
\qed

\lem {\em 
\label{Cartan_matr_general}
a) Let
$\sR$ be a quasifinite domain over $\Z$
and $\AA$ be a finite rank flat
algebra over
$\sR$.
Then there exist  a quasifinite domain
$\sR'$ over $\sR$
and collections
$L_i$, $P_i$ of $\AA_{\sR'}$-modules indexed by
$I=Irr(\AA_{\barr\Q})$,
such that
for any geometric point $\k$ of $\sR'$ the corresponding $\AA_\k$ modules
$(L_{i})_\k$, $(P_{i})_\k,\ i\in I$,
provide complete nonrepeating lists of
irreducible
$\AA_\k$-modules
and their
projective covers.

b)
The same holds for topological algebras of the form $\varprojlim \AA_i$ where
$\AA_i$ are as in (a) and
$\AA_i\to \AA_{i+1}$ is a
{\em square zero extension}\ie
a surjective homomorphism with zero multiplication on the kernel.

}

{\em Sketch of proof.}
Part (b) follows from (a) since
irreducibles for $(\AA_i)_\k$ do not depend on $i$
and any projective module admits a lifting to a projective module
over a square zero extension.


(a)
Let $L_i^o,\ \ii$,
be the complete list of
irreducible modules for $\AA_{\Qubar}$
and let
$P_i^o$ be the indecomposable projective
cover of $L_i^o$.
These modules lift to
some finite extension $\sR_1$ of $\sR$\ie
there are
$\AA_{\sR_1}$-modules
$P_i$, $L_i$
 such that
$(P_i)_{\Qubar}\cong P_i^o$, $(L_i)_{\Qubar} \cong L_i^o$.
(For that choose
a presentation for a module as a cokernel of a map between finitely
generated free modules and multiply the matrix of the map by an element in $\sR$ to make
its entries integral over $\Z$. Then $\sR_1$ is obtained by adjoining the new entries to $\sR$.)

After a finite localization
$\sR_2$ of $\sR_1$
we can also achieve that each $P_i$ is projective.
The reason is that for some $d_i>0$,
the sum $\oplus_I\ (P_i)_{\Q}^{\oplus d_i}$ is a free
$\AA_\Q$-modules, so $\oplus_I\ P_i^{\oplus d_i}$ becomes a free
module after a finite localization $\sR_2$.

After further replacing $\sR_2$ by a finite localization $\sR_3$
we can achieve
that $(P_i)_\k$ are pairwise non-isomorphic indecomposable modules.
For this consider two $\AA_{\sR_2}$-modules $M,N$ and a
commutative $\sR_2$ algebra $\fR$.
The   existence of  an isomorphism
$M_{\fR}\cong N_{\fR}$ and 
a nontrivial idempotent in $End(M_{\fR})$;\ amounts
to existence of an $\fR$ point of a certain affine algebraic variety over
$\sR_2$.
Such an algebraic variety has a $\Qubar$ point iff it has a
$\k$ point for all finite geometric
points $\k$ of $Spec(\sR_2)$ of  almost any prime characteristic.
Here, $\sR_3$ is obtained from $\sR_2$ by inverting finitely many primes.

Finally,
$\Hom(P_i,L_j)$ is a finite $\sR_3$-module such that
$\Hom(P_i,L_j)\otimes \Qubar \cong \Qubar^{\delta_{ij}}$. Thus after
a
finite
localization $\sR'$ we can assume
that
$\Hom(P_i,L_j)\cong (\sR')^{\delta_{ij}}$.
\qed

\rem
Our proof of Lusztig conjectures
about numerics of modular representations
gives the result for $p\gg 0$ only,
because
we rely
on the (very general) Lemma
\ref{Cartan_matr_general},
which is not constructive in the sense that it
provides
no information on the set of primes one needs to invert in $\sR$ to get
$\sR'$ with needed
properties.
The same problem appears in \cite{AJS}
which contains results that are
essentially equivalent to
the case $e=0$ of our results.
Notice that
Fiebig  has found an explicit (and very large) bound $M$
and showed that the argument of \cite{AJS} works for $p>M$
\cite{Fiebig bound}.

\rem
\lab{A and E}
The following is just the extra information one finds in the proofs of
Theorem \ref{Moreq}
and Proposition \ref{EiLi}.


\scor {\em
Let $\EE$ be a vector bundle on $\gt$ as in \ref{tstr},
then the algebra $\Abig$ in Theorem \ref{Moreq}
can be chosen as $\Abig=\End(\E)^{op}$.
For any $e\in \NN(R)$, the collection $\E_i$ in Proposition
\ref{EiLi}, \
can be chosen as
representatives of isomorphism classes
of
indecomposable constituents of the restriction
of $\EE$ to
$\hatt {\BB_{\RR, e}}$.


$I_ e$ is the set of isomorphism classes of irreducible modules of
$\Abig|_{\hatt e}\ \ten_R\C$ (i.e., of
$\Abig_e\ten_R\C$ for
$\Abig_e= \Abig\ten_{\OO(\fg)}\OO_e$).

}

\sss{
Numerical consequences
}
For $\k\in FGP$
and $e\in\NN(\k)$ recall  the isomorphisms of Grothendieck
groups
\begin{equation}\label{KKKK}
\begin{array}{lll}
K(U_{\k, e}^0)\cong K(\BBke),\\
K(\BBke)_\Qu\cong K(\BB_{\Qubar, e})_\Qu,\\
K(\BB_{\Qubar, e})\cong K((\Abig_\Qubar)_ e^0)
.\end{array}
\end{equation}
The first isomorphism is \cite{BMR1} Corollary 5.4.3,
  the second one follows from \cite{BMR1} Lemma 7.2.1 and
Proposition 7.1.7, and
the last one is immediate  from Theorem \ref{tstr} above.
By ``homotopy invariance of the Grothendieck group'' we have
$K^0[mod^{fg}(U_{\k,e}^0)]\
=K^0[mod^{fg}(U_{\k,e}^{\hat 0})]
=K^0[mod^{fg}(U_{\k,\hat  e}^0)]$,
and
the
same for $U$ replaced by $\Abig$.

\scor {\em
\label{multbigp}
 a) For almost all characteristics $p$ the composed isomorphism
$$
K[mod^{fg}(U_{\k, e}^ 0)]_\Qu
\cong
K[mod^{fg}((\Abig_\Qubar)_ e^0)]_\Qu
$$
sends classes of irreducibles to classes of
irreducibles and classes of indecomposable
projectives to classes of indecomposable projectives;
the
same applies to $U^{\hat 0}_{\k, e}$ and $U^{0}_{\k,\hat  e}$.

(b) For all $p> h$
the composition  sends the class of the dual of baby Verma $U_\k$-module
of ``highest'' weight zero,
to the class
of a structure sheaf of a point.\fttt{
A baby Verma module for $U_{\k,e}^0$ involves the data
of a Borel subalgebra
$\fb\in\BBke$.
Notice that it is easy to see that
the class $[\OO_\fb]\in K(\BBke)$ is independent of the choice of
$\fb$.
}

(c)
Multiplicities of irreducible modules in baby Verma modules
are independent of $p$ for large $p$.


}

\sus{
Equivariant versions and Slodowy slices
}
\label{Equivariant versions and Slodowy slices}

From now on
we denote by $\fG$ a copy of the group
$G\tim \Gm$, supplied
with a map
$\fG@>>>G\tim\Gm$ by $(g,t)\mm (g,t^2)$.
Since
$G\times \Gm$ acts on $\g$
by $(g,t): x\mapsto t\cd ^gx$ the group
$\fG$
acts on $\g$
by $(g,t): x\mapsto \ t^2\ ^gx$,

\sss{
Equivariant tilting t-structures
}

\spro {\em \label{equiva}
Let $R'$ be a Noetherian $R$-algebra.
Let $S\to \g_{R'}$ be a map of affine schemes over $R'$
which satisfies the Tor vanishing condition
from \ref{BCact} relative to $\gt_{R'}$.

Let $H$ be a flat  affine algebraic group over $R'$ endowed with
a morphism
$\phi:H\to \fG_{R'}$. We assume that
$H$ acts on $S$ and that the map $S\to \g_{R'}$ \
(respectively $S\to \N_{R'}$)
is $H$-equivariant.
Then we have a natural equivalence
$$
D^b(Coh^H(\tii S))
\cong
D^b(mod^{H,fg}(\Abig_S))
,\ \ \ \tx{respectively}\ \ \
D^b(Coh^H(\tii S'))
\cong
D^b(mod^{H,fg}(\Az_S))
.$$

}            \pf We construct the first equivalence,
the second one works the same.
The tilting bundle $\EE$ constructed in section \ref{sect2}
is manifestly $G\times \Gm$ equivariant.
Thus the vector bundle $\EE_S$
and the algebra $\Abig_S$
 carry natural $H$-equivariant
structures.  Therefore,
for $\F\in Coh^H(\tii S)$ and $ M\in mod^{H,fg}(\Abig_S)$,
the $\Abig_S$-modules
$\Hom(\E_S, \F)$ and
$M\Lotimes _{\Abig_S}\E_S\in\Coh^H(\tii S)$
carry
natural $H$-equivariant structures.
Passing to the derived functors we get two adjoint functors
$\CC^H,\II^H$
between
$
D^b(Coh^H(\tii S))
$ and
$
D^b(mod^{H,fg}(\Abig_S))
$.
A standard argument shows that these functors are compatible with
the pair of adjoint functors
$\CC,\II$ between $D^b(Coh(\tii S))$
and $D^b(mod^{fg}(\Abig_S))$,
which are given by the same formulas.
Moreover, this compatibility extends to the
adjunction morphisms between the identity functors
and compositions $\II^H\CC^H,\CC^H\II^H$
and
$\II\CC,\CC\II$.
Since the adjunction
morphisms are isomorphisms in the non-equivariant setting,
they are also isomorphisms in the equivariant one.
\qed

\sss{
Slodowy slices
}
\lab{Slodowy slices I}
\label{Slodowy slice}
Let $\k$ be  a geometric point
of $R$.
Fix a nilpotent $e\in \N(\k)$ such that
there is
a
homomorphism $\varphi:SL(2)\to G$ with
$d\varphi\smat 0&1\\0&0\esmat= e$.\fttt{
This is always possible if $p>3h-3$
\cite{Hu}.
}
The corresponding     $sl_2$ triple    $e,h,f$ defines a Slodowy slice
$S_{\k,e}\dff e+Z_{\fg_\k}(f)$ transversal to the conjugacy class of $e$.

Let $C$ be a maximal torus
in the centralizer of the image of $\varphi$, it is also
a maximal torus in the centralizer $G_e$ of
$e$.
Let $\phi:\Gm\to G$ by $\phi(t)\dff \varphi\smat t&0\\ 0&t\inv\esmat$.
We denote by
$\fGm$
a copy of  the group $\Gm$
and by
$\Ctil$ the group $C\tim\fGm$ supplied with a morphism
$\bbi:\Ctil\to \fG$ by
$\bbi(c,t)\dff (c\phi(t),t\inv)$.
The action of $\Ctil$ on $\fg$
 is by $(c,t)\ x=\ t^{-2}\ ^{c\phi(t)}x$, it
preserves the Slodowy
slice $S_{\k,e}$ and the action of $t\in\fGm$ contracts it
to
$e$ for $t\to\yy$.


Now we use the fact that
$C$ contains
$\phi(-1)\bz$ for some
$\bz\in  Z(G)$,
see Appendix \ref{Eric}.
Notice that
the element
$\bbb m=(\phi(-1)\bz,-1)\in\Ctil$
is sent to
$(\phi(-1)\bz \phi(-1),-1)=\ (\bz,-1)\in\ \fG$.
So, the action of $\Ctil$ on $\fg$
factors through the quotient
 by the  subgroup
generated by
$\bbb m$.

\sss{
Equivariant lifts of irreducibles and indecomposable projectives
}
\lab{Equivariant lifts of irreducibles and indecomposable projectives}

Let $\k$ be a geometric point of $R$ and let
$e\in\NN(\k)$.
Let $H$ be a torus endowed with  a map into  the stabilizer of $e$
in $\fG_\k$.
We say that for an object $\F$ in the derived category of coherent sheaves
its equivariant lifting
is an object in the equivariant derived category whose image under
forgetting the equivariance functor is isomorphic to $\F$.

\spro {\em
\label{gra}
Let $H$ be a $\k$-torus mapping to $\fG_\k$ and fixing
$e,h,f$, so that in particular  it
preserves $\Ske$.

a) Every irreducible exotic sheaf $\LL$ on either of the spaces:
$\Stilke$, $\Stilkep$, $\hBBke$, $\hBBkep$,
$\gt_\k$, $\Nt_\k$,
 whose (set-theoretic) support is contained in  $\BBke$
admits an $H$-equivariant lift $\tii \LL$. Any other lift is isomorphic to
a twist of $\tii \LL$
by a character of $H$.

b)
Every
projective exotic sheaf
$\WW$ on either of the spaces:
$\Stilke$, $\Stilkep$,
$\hBBke$, $\hBBkep$  admits an $H$-equivariant lifting
$\tii \WW$.
If $\WW$ is indecomposable then
every  equivariant lifting of
$\WW$ is isomorphic to a
twist of $\tii \WW$ by a character
of $H$.

c) Assume that $char(\k)=0$.
Then there exists a quasifinite $R$-domain
$R'$, such that
the following holds.
\bi
i) The nilpotent $e$,
torus $H$ and the homomorphism $H\to \fG$ are
defined over $R'$.
\i
ii) For each $\LL\in D^b(Coh_\BBke(\Stilkep))$, $\WW\in \Coh(\Stilkep)$
and equivariant lifts $\tii \LL$, $\tii \WW$ as above there exist
 $\LL_{R'},\ \WW_{R'}\in D^b(Coh^{H_{R'}}(\tii S_{R',e}'))$
and their equivariant lifts
 $\Ltil_{R'},\ \tii\WW_{R'}$, such that their
base change to $\k$ is isomorphic to $\LL,\WW,\tii \LL,\tii\WW$.
\i
iii) For every geometric point $\k'$ of  $R'$ the base change
to $\k'$ of
$ \LL_{R'}$, $\WW_{R'}$
and
$ \tii\LL_{R'}$, $\tii\WW_{R'}$
are,
respectively, irreducible and indecomposable projective equivariant
exotic sheaves.
Every equivariant
irreducible or indecomposable projective exotic sheaf
on
$\tii S_{\k',e}$,
$\tii S_{\k',e}'$,
$\hatt{\BB_{k',e}}$,
$\hatt{\BB_{k',e}}'$
arises in this way.
\ei

}            \pf
It is clear that the direct image
of an irreducible exotic sheaf
$\LL\in D^b (Coh_{\BBkep}(\Stilkep))$
 under each of
the closed embeddings
$\Stilkep\imbed \Nt_\k$, $\Stilkep\imbed \Stilke \imbed \gt_\k$
is again an irreducible
exotic sheaf supported on $\BBkep$, and all irreducible
exotic sheaves on these spaces supported on $\BBke$
are obtained this way. The resulting  sheaves on $\gt_\k$, $\Nt_\k$
can be also
thought of as sheaves on $\hBBke$, $\hBBkep$ respectively. Thus in (a)
 it suffices
to consider the case of $\Stilkep$ only.

 Applying Proposition \ref{equiva}
we get an equivalence $D^b(Coh^H(\Stilkep))\cong
D^b(mod^{H,fg} (\Az_\Ske))$. Thus, statement (a) reduces to showing that every
  irreducible $\Az_S$ module with central character
$e$ admits an $H$ equivariant structure.
 The torus $H$ acts trivially on the finite set of irreducible
$\Az_S$-modules with central character $e$.
 It follows that $H$ acts projectively on such a representation.
 Since every cocharacter of $PGL(n)_\k$
 admits a lifting to $GL(n)_\k$, we see that the representation admits an $H$ equivariant structure.
 This proves (a).

 Similarly, in order to check (b) it
suffices to equip indecomposable projective modules over the respective
 algebras with an $H$ equivariant structure. By a standard argument,
projective cover of an irreducible
 module in the category of graded (equivalently, $H$-equivariant) modules is also a projective cover
 in the category of non-graded modules, which yields (b).

To check (c), it suffices to consider
equivariant projective modules on $\Stilkep$
(then the rest follows as in \ref{EiLi}).
Equivariant
indecomposable projectives over $\k$ are direct summands of
$\E_{\Stilkep}$.
We can find
 a quasifinite $R$-domain $R'$
such that the corresponding idempotents are defined and orthogonal over $R'$.
 \qed

\sss{
Equivariant localization
}
\label{BMRequivav}
\lab{Equivariant localization}

In this subsection we link $C$-equivariant exotic sheaves to
representations graded by weights of $C$ by proving Theorem \ref{163equivav}.
In this argument it will be important to distinguish between a variety
and its Frobenius twist, so we bring the twist back into the notations.
We concentrate on the first equivalence, the second one is similar.

 The torus $C\tw$ acts on  $\hBBke\tw $, composing this action with the
Frobenius morphism $C\to C\tw$ we get an action of $C$ on $\hBBke\tw$.
Consider the category of equivariant coherent sheaves
$Coh^C(\hBBke\tw)$. The finite group scheme $C_1=Ker(C\overset{Fr}{\To}C\tw)$
maps to automorphisms of the identity functor in this category;
since the category of $C_1$-modules is semisimple with simple objects
indexed by $\fc^*(\Fp)=X^*(C)/p$,
the category splits into a direct sum
\begin{equation}\label{cohspl}
Coh^C(\hBBke\tw)=\oplusl
_{\eta\in X^*(C)/p} Coh^C_\eta
\end{equation}
 where  $Coh^C_\eta$ consists of
such equivariant sheaves that $C_1$ acts on each fiber by the character $\eta$.
Notice that $Coh^C_0(\hBBke\tw)\cong Coh^{C\tw}(\hBBke\tw)$ canonically
and for every
 $\tii \eta\in X^*(C)$ the functor of twisting
by $\tii \eta$ provides an equivalence $Coh^C_\theta(\hBBke\tw)\cong
Coh^C_{\theta+\eta}(\hBBke\tw)$ where $\eta=\tii \eta\,\mod p  X^*(C)$.

The sheaf of algebras $\tii \DD$
is equivariant with respect to the $G$
action, hence $\tii \DD|_{\hBBke\tw}$ is equivariant with respect to the
action of $C$ on $\hBBke\tw$. Consider the category
$
mod^{C,fg}(\tii \DD|_{\hBBke\tw})
$
of $C$-equivariant coherent sheaves of modules
over this $C$-equivariant sheaf of algebras; here "coherent" refers to coherence as a sheaf
of $\OO$-modules over the scheme $\hBBke\tw$.
A sheaf $\F\in
mod^{C,fg}(\tii \DD|_{\hBBke\tw})$ carries two commuting actions
of $\fc$, $\al_C$ and $\al_\g$ (see \ref{163equivav}) whose difference
commutes with the action of $\tii \DD$.



It is easy to see (cf. \ref{163equivav})
that every $\F\in
mod^{C,fg}(\tii \DD|_{\hBBke\tw})$ splits as a direct sum
$\F=\oplusl_{\eta\in \fc^*(\Fp)} \F_\eta$ where
$\al_\g(x)-\al_C(x)-
\langle
x,\eta
\rangle Id
$ induces a pro-nilpotent endomorphism of $\F_\eta$.
Here by a pro-nilpotent endomorphism we mean one which becomes nilpotent
when restricted to any finite nilpotent neighborhood of $\BB_e\tw$.


Thus we get a decomposition of the category

\begin{equation}\label{modDDspl}
mod^{C,fg}(\tii \DD|_{\hBBke\tw})=\oplusl_{\eta\in \fc^*(\Fp)}
mod^{C,fg}_\eta(\tii \DD|_{\hBBke\tw}).
\end{equation}

Let
$\EE$ be a
splitting bundle for the Azumaya algebra $\tii \DD$ on  $\hBBke\tw $.
We claim that $\EE$ admits a $C$-equivariant structure
compatible with the equivariant structure on $\tii \DD|_{\hBBke\tw}=
\EEnd(\EE)$.
For $\la,\mu\in \h^*(\Fp)$ the bimodule providing Morita equivalence between the restrictions
of $\tii \DD$ to the formal neighborhoods of the preimages of $\la$ and $\mu$
under the projections from the spectrum of the center of $\tii \DD$
to $\h^*$ (see \cite[2.3]{BMR1}) is manifestly $G$-equivariant, thus it suffices
to consider the case $\la=-\rho$. Then $\tii \DD|_{\hBBke\tw}$ is identified
with the pull-back of a $C$-equivariant
 Azumaya algebra on the formal neighborhood of
$e\in \g^*\tw$. We now construct compatible
$C$-equivariant splitting bundles on the $n$-th
infinitesimal neighborhood of $e$ for all $n$ by induction in $n$.
The base of induction follows from the fact that every extension
of $C$ by $\Gm$ splits, while the induction step follows from splitting
of an extension of $C$ by an additive group. Thus existence of a $C$-equivariant
splitting bundle is established.

We fix such a $C$-equivariant structure
on $\EE$; we can and will assume that the restriction of the resulting
equivariant $\tii \DD$ module
to a finite order neighborhood of $\BBke$ belongs to
 $mod^{C,fg}(\tii \DD|_{\hBBke\tw})$ (since this restriction is
an indecomposable sheaf of $\tii \DD$-modules, this can be achieved
by twisting an arbitrarily chosen $C$-equivariant lift of $\EE$ by a character
of $C$).
Then we get a functor $\F \mapsto \F \otimes _\OO \EE$
from  $\Coh^{C\tw}(\hBBke\tw)$
to the category
$
mod^{C,fg}(\tii \DD|_{\hBBke\tw})
$
of $C$-equivariant sheaves of modules
over a $C$-equivariant sheaf of algebras.
We will compose this functor
with  the global sections functor
$D^b[mod^{C,fg}(\tii \DD|_{\hBBke\tw})]
\to
D^b[mod^{C,fg}(U_{\hatt{e}}^{\hatt{\la}})]$.
We claim that
the composition lands in the full subcategory
$D^b[mod^{fg}(U_{\hatt{e}}^{\hatt{\la}}, C)]$ and provides the desired
equivalence. This follows from the next

\slem {\em
a) The derived global  sections functor provides an equivalence
$$
\RGa:
D^b[mod^{C,fg}(\tii \DD|_{\hBBke\tw})]
\to D^b[mod^{C,fg}(U_{\hatt{e}}^{\hatt{\la}})]
.$$

b)
The functor $\F\mapsto \F\otimes \EE$ provides an equivalence
$\Coh^{C}(\hBBke\tw)\to
mod^{C,fg}(\tii \DD|_{\hBBke\tw})$.

c)
Both of these  equivalences are compatible
with the canonical  $\fc^*({\Bbb F}_p)$-decompositions of categories.
In particular for $0\in\fc^*({\Bbb F}_p)$, we get equivalences
$$
mod^{fg}(U_{\hatt{e}}^{\hatt{\la}},C)
\ \cong\
\Cohe^{C\tw}(\hBBke\tw)
.$$

}            \pf (a) follows by the argument of Proposition
\ref{equiva} (the two
adjoint functors
commute with forgetting the equivariance).

(b) is  just the observation  that once $\EE$ is equivariant,
the standard equivalence between coherent sheaves
on $\gt$ and modules over the sheaf of algebras  $\AA=\EEnd(\EE)$
extends to  the equivariant setting.

(c) follows from the definition of the decompositions. \epf

\sus{
Gradings and bases in K-theory
}
\label{Gradings and bases in K-theory}

In this subsection
we work over
a geometric point $\k$ of $R$.
We reduce  the conjectures
of \cite{Kth2} which
motivated this project
to a certain property $\bst$ of exotic
sheaves.
All substantial proofs of claims in this subsection are
postponed
to subsection \ref{Proof of Theorem and Proposition}.
In the next section \ref{sect6}
we will see that property $\bst$
(thus a proof of Lusztig's conjectures)
follows ({\em for large $p$})  from the results of \cite{AB}.

\sss{Construction of gradings}
\lab{Construction of gradings}
On the formal neighborhood
$\hBBke$ of a Springer fiber $\BBke$
in $\gt_\k$  there is a canonical (up to an isomorphism)
vector bundle
$\pl_{i\in \II_e}\ \E_i$
which is a minimal
projective generator for the heart $\Cohe(\hBBke)$
of the exotic
t-structure.
We just take $\E_i$'s to be
representatives for
isomorphism classes of
indecomposable summands in the
pull back
$\E|_\hBBke$ of
any vector bundle $\E$ from Theorem \ref{tstr}.


By Proposition \ref{gra},
vector bundles $\E_i$ admit a $\fGm$ equivariant structure.
We temporarily fix such a structure in an arbitrary way and
let $\Etil_i\in \Coh^\fGm(\hBBke)$
denote the resulting equivariant bundle.
In view of
Lemma \ref{gra}, any other choice yields
an equivariant vector
bundle isomorphic to
a twist $\Etil_i(d)$ by some character  $d\in\Z$ of $\fGm$.

Consider the restriction
of $\Etil_i $
to the formal neighborhood of
$\BBke$ in $\Stilkep$.
Since $\fGm\subset \Ctil$ acts on $\Stilkep$ contracting it to
the projective variety $\BBke$,
there exists a unique (up to a unique
isomorphism),
$\fGm$ equivariant
vector bundle on $\Stilkep$ whose pull-back to the formal
neighborhood of
$\BBke$ is identified with
$\Etil_i|$. 
We denote
this vector bundle by
 $\Etil_i^S\in\Coh^\fGm(\Stilkep)$.

To summarize,   $\Etil_i^S$
is a set  of representatives for equivalence classes
of  indecomposable projective $\fGm$-equivariant
exotic sheaves on $\Stilkep$ modulo $\fGm$-shifts.

\sss{
Property $\bst$
}
\lab{Property bst}
The following statement will be the key to the proof of Lusztig's conjectures.
$$
\bst
\Ca
\tx{\em
There exists a choice
of $\fGm$-equivariant lifts
$\Etil_i,\ i\in
\II_e$,
which is
}
\ \
\ \ \
\ \ \
\ \ \
\ \ \
\ \ \
\ \ \
\ \ \
\ \ \
\\
\tx{\em
invariant under
the action
of the centralizer of $e,h$ in $G$ and such that:
}
\ \ \
\ \ \
\ \ \
\ \ \
\ \ \
\\
\tx{(i$_\st$)}
\ \ \
\ \ \
\Hom_{\Coh^\fGm(\Stilkep)}(\Etil_i^S, \Etil_j^S(d))=0\ \ \ {\rm for}\ \ \  d>0,
\ \ \
\ \ \
\ \ \
\ \ \
\ \ \
\ \ \
\ \ \
\ \ \
\ \ \
\ \ \
\ \ \
\ \ \
\ \ \
\ \ \
\ \ \
\ \ \
\ \ \
\\
\tx{(ii$_\st$)}
\ \ \
\ \ \
\Hom_{\Coh^\fGm(\Stilkep)}(\Etil_i^S,\Etil_j^S)\ =\ \bbk^{\delta_{ij}}.
\
\ \ \
\ \ \
\ \ \
\ \ \
\ \ \
\ \ \
\ \ \
\ \ \
\ \ \
\ \ \
\ \ \
\ \ \
\ \ \
\ \ \
\ \ \
\ \ \
\ \ \
\ \ \
\ \ \
\ \ \
\ \ \
\ \ \
\\
\tx{\em Such a choice is unique up to twisting all $\Etil_i$
by the same character of $\fGm$.}
\ \ \
\ \ \
\ \ \
\ \ \
\Eca
$$
In other words
	$\Hom_{\Coh(\Stilkep)}(\Etil_i^S, \Etil_j^S)$ has no negative
$\fGm$ weights and the zero weight spaces are spanned by
the  identity maps.

\medskip


\sss{Normalizations}
\label{Normalizations}
Assuming that $\bst$ holds
for $\k$ and $e$,
we will reserve notation $\Etil_i$ for
$\fGm$-equivariant vector bundles
on $\hBBke$,
satisfying the above conditions.
Recall that the vector bundle
which is a tilting generator for the exotic t-structure
has $\OO$ as a direct summand,
this implies that $\E_{i_0}\cong \O_{\hBBkep}$ for some
$i_0\in\II_e$.
We will assume that
$\Etil_{i_0}=\O_{\hBBke}(2\dim \BB_e)$;
since the collection $\Etil_i$ is  unique up to a simultaneous twist
by some character of $\fGm$,
this   fixes the set of isomorphism classes of $\tii\E_i$ uniquely.



According to  Proposition
\ref{Slodowy slices I}
we can equip the vector bundle $\Etil_i$
with some $\Ctil$ equivariant structure compatible with the
$\fGm$ equivariant structure fixed above.
Let $\Etil_{i,0}$ denote
the resulting $\Ctil$-equivariant vector bundle and set
$\Etil_{i,\la}=\Etil_{i,0}(\la)$ for $\la\in X^*(C)$.
As above, from $\Etil_{i,\la}$ we obtain
a  $\Ctil$-equivariant vector bundle
$\Etil_{i,\la}^S$ on $\Stilkep$.
As we vary $i\in\II_e$ and $\la\in X^*(C)$,\
vector bundles $\Etil_{i,\la}$ on $\hBBke$
(resp., $\Etil_{i,\la}^S$
on $\Stilkep$) form a
complete list of representatives  modulo  $\fGm$-shifts, of
indecomposable exotic projectives in
$D^b[\Coh^\Ctil(\hBBke)]$
(resp.  $D^b[\Coh^\Ctil(\Stilkep)]$).
Also, we
define irreducible exotic objects
$\Ltil_{i,\la}$
of
$D^b(\Coh^\Ctil(\hBBke))$ (respectively, $\Ltil_{i}$
of
$D^b(\Coh^\fGm(\hBBke))$)
such that
$\Ext^\bu_{D^b\Coh^{\Ctil}(\hBBke)}
(\Etil_{i,\la}, \Ltil_{j,\mu})=\
\k^{\delta^i_j \delta^\la_\mu}
$, $\Ext^\bu_{D^b\Coh^{\fGm}(\hBBke)}
(\Etil_{i}, \Ltil_{j})=\
\k^{\delta^i_j}
$.


\sss{Uniformity of K-groups}
\lab{Uniformity of K-groups}
Although not strictly necessary for the
the proof of Conjectures, the following result allows a neater formulation
of the next Theorem and clarifies the picture

Fix two geometric points of $Spec(R)$:\ $\k_0$ of characteristic
zero and $\k$ of characteristic $p>h$. Recall that for a flat
Noetherian scheme
$X$ over $R$ one has
{\em  specialization map}
$Sp_X: K(X_{\k_0})\to K(X_{\k})$.

\spro {\em
\label{ptoCe}
The maps
$Sp_{\BBe}$, $Sp_{\Stile}$,
$Sp_{\Stilep}$ are isomorphisms.

}

\pf
Tensoring the maps with $\Q$
we get isomorphisms because
the modified Chern character map identifies
$K(\BBe)\otimes \Qlb$ with
the dual
of $l$-adic
cohomology:
in \cite{BMR1} it was shown in Lemmas 7.4.2 and 7.4.1
 that the map is injective,
 however the dimensions are the same by Theorem 7.1.1 and Lemma 7.4.3.
The
 independence of $l$-adic cohomology on the base field
was established by Lusztig
(\cite{Lu2} section 24, in particular theorem 24.8 and
subsection 24.10).

On the other hand,
using the above equivalences of categories we see that
over the field $\k$ the classes of
irreducible (respectively, projective) objects, form bases in respective
Grothendieck groups.
In particular, each of the K-groups is a free
abelian group of finite rank and the Ext pairing between
$K(\BBke)$ and $K(\Stilke)$ is perfect.
The corresponding statements for $\k_0$ were proved by Lusztig.
It is clear that the specialization
map is compatible with this pairing.
Thus the pair of maps
$Sp_{\BBe}$, $Sp_{\Stile}$ is an example of the
following situation.
We are given free abelian groups $A$, $A'$, $B$, $B'$,
all of the same finite rank, maps $F:A\to B$, $F':A'\to B'$
and perfect pairings $A\times A'\to \Zet$, $B\times B'\to \Zet$ such that
$\langle F(x), F'(y)\rangle = \langle x, y\rangle$ for all $x\in
A$, $y\in A'$. It is clear that in this situation $F$, $F'$ have to be
isomorphisms. \epf

\sss{
Reduction of  Lusztig's conjectures
to $\bst$
}
Lusztig's conjectures from \cite{Kth2}
will be recalled in detail  in
sections \ref{Recollection},
\ref{Lusztig's formulation of conjectures and  how it follows} and
\ref{Three Involutions}.
In \ref{Recollection}
we will state our precise results and then we will see
in
\ref{Lusztig's formulation of conjectures and  how it follows}
that they imply the following Theorem.

\stheo {\em 
\label{LuLu}
\label{Lusztig Conjectures}
(1)
If $\bst$ holds for some  geometric point $\k$ of $R'$, then
Conjectures
 5.12, 5.16 of \cite{Kth2} hold
(existence of certain signed bases of K-groups of
complex schemes $\BBCe$
	and
$\StilCep$).
The two bases discussed in Conjecture 5.12
are given by the
classes of
$\Ltil_{i,\la}^\k(-2\dim\BB_e)$,
$\Etil_{i,\la}^\k$, $i\in\II_e,\ \la\in X^*(C)$;
and the two bases discussed in Conjecture
5.16
are given by the
classes of
$\Ltil_i^\k$, $\Etil_i^\k$.\fttt{
Super index $\k$ means that we use the
sheaves defined over $\k$. A priori these define elements of
K-groups for $\k$-schemes, however by  Proposition
\ref{Uniformity of K-groups} these K-groups are canonically identified for
all  $\k$.
}

(2) If $\bst$ holds for
some $\k\in FGP$,
then Conjecture 17.2 of {\em loc. cit.}
 (relation to modular
representations over $\k$;
omit the last paragraph in 17.2 on the quantum
version\fttt{The quantum version is
closely related to Conjecture \ref{conj_Uq} above.})
holds for $\k$.

}

\sss{
A reformulation of Lusztig's conjectures
}
\lab{Recollection}
Here we
formulate
a list of  properties
that naturally appear from the present point of view.
Then in
\ref{Lusztig's formulation of conjectures and  how it follows}
we will recall Lusztig's conjectures and
show that they follow from these properties.
We will omit \cite[Conjecture 5.16]{Kth2}, 
since it is similar to {\em loc. cit.} Conjecture 5.12;
the only
difference is
that
5.12 deals with coherent sheaves equivariant with respect
to the torus $\Ctil$,
while
5.16 is about sheaves equivariant with
respect to a one-parameter subgroup $\fGm\subset \Ctil$.
So the existence of  bases with properties from 5.12 implies the same for
5.16  and Lusztig's  uniqueness argument (recalled in footnote
(\ref{standard argument})) applies equally to  both conjectures.

The K-group of a torus $T$ is the group
algebra of its character lattice
$R_T\dff\Z[X^*(T)]$, it contains
a {\em subsemiring}
$R_T^+=\Z_+[X^*(T)]$.
\label{RpC}
Let $\fK_T\ \dff Frac(R_{T})$
be the fraction field of $R_T$.
Denote by
$\del:R_\Ctil\to R_\fGm$
the constant coefficient map
$\sum_{\nu\in X^*(C)}\ p_\nu[\nu]\mm p_0$
(\cite[5.9]{Kth2}).

Recall from \cite[Theorem 1.14.c]{Kth2}
that the direct image map
gives an embedding
$K^{\Ctil}(\BBke)\inj K^{\Ctil}(\Stilkep)$.
Lusztig's conjectures
involve certain involutions
$\bee$ on $K^{\Ctil}(\BBke)$
and
$\beS$ on $ K^{\Ctil}(\Stilkep)$
(denoted $\tii\beta$, $\beta$ in \cite{Kth2}),
a
certain
pairing $(\  \Vert \ )$
on  $K(\Stilkep)$ with values in
the fraction field
$\fK_\Ctil$,
and
 a certain element $\nabla_e$ of $\fK_\Ctil$.
We denote by $P\mapsto
P^\sv
$ the involution of $\fK_\Ctil$
corresponding to inversion on $\Ctil$.

The following proposition will
be verified
in
\ref{Proof of geometric formulation of Conjectures}.

\spro {\em
\lab{geometric formulation of Conjectures}
Let $\k$ be a geometric point of $R$ such that $\bst$ holds for $\k$.
Define $\Etil_i^S,\ \Etil_{i,\la}^S$, $\Ltil_i,\ \Ltil_{i,\la}$
as in
\ref{Normalizations}.

\sff{A}
The following subsets are bases over the ring $R_\fGm=\Z[v^\pmo]$
(we will often omit the super index $\k$):
$$
\bfBS^\k\ \dff\{[\Etil_{i,\la}^S];
\ i\in\II_e,\ \la\in X^*(C)\}
\sub K^\Ctil(\Stilkep)
,$$
$$
\bfBe^\k\dff\{v^{-2 \dim \BB_e} [\Ltil_{i,\la}];
\ i\in\II_e,\ \la\in X^*(C)\}
\}\sub K^\Ctil(\BBke)
.$$

Elements of
$\bfBS$ are fixed by $\beS$ and elements of $\bfBe$
by $\bee$:

\begin{equation}
\label{tiibeta}
\bee(v^{-2 \dim \BB_e}[\Ltil_{i,\la}])=v^{-2 \dim \BB_e}[\Ltil_{i,\la}],
\end{equation}
\begin{equation}
\label{beta}
\label{betaS}
\beS([\Etil_{i,\la}^S])=[\Etil_{i,\la}^S].
\end{equation}
Both bases satisfy the condition of
{\em
asymptotic orthonormality}:
\fttt{
To make sense of it
we use the embedding
$R_{\Ctil}\subset \RC((v^{-1}))$ and the induced
embedding of fraction fields.
}
\begin{equation}
\label{asymptotic orthonormality}
(b_1 \Vert b_2)\in
\Ca
\ \ \ \ \ \ v^{-1}\RC[[v^{-1}]]
&
\ \ \tx{if}
&
\ \ \ \ \ \ \ \tx{$b_2\not\in X^*(C)b_1$},
\\
1+v^{-1}\RC[[v^{-1}]]
&
\ \ \tx{if}
&
\tx{$b_2=b_1$}.
\Eca
\end{equation}

\sff{B}
The
two bases are dual for the pairing
$(\ \Vert\ )$.

\sff{C}
$
\lab{nabla}
(\bfBS \Vert \bfBS)\sub\ \fra{1}{\nabla_e}\RCp[v^{-1}]
\cap \RCp[[v^{-1}]]
$
and
$
(\bfBe \Vert \bfBe)\sub\ \RC[v^{-1}]
$.
\fttt{
For basis $\bfBe$ the corresponding positivity statement
$(\bfBe \Vert \bfBe)\sub\ \RCp[v^{-1}]$
follows from the result of section
\ref{Conjecture on Koszul property} below.
}


\sff{D}
For
$b_1,b_2\in \bfBS$, each of
 the coefficient polynomials $c^\nu_{b_1,b_2}\in \Z_+[v^\inv],\ \nu\in
X^*(C)$,
in the expansion
$\nabla_e(b_1\Vert b_2)=\
\sum_{\nu\in X^*(C)}\ c^\nu_{b_1,b_2}[\nu]\in \RCp[v^{-1}]$,
is either even or odd.

\sff{E}
Basis
$\bfBS$
is the unique $R_\fGm$ basis of
$K^\Ctil(\Stilkep)$
which is pointwise fixed by $\beS$,
satisfies asymptotic orthonormality,
a normalization
property
$\bfBS\ni\ v^{2\dim(\BB_e)}[\O_{\Stilkep}]$
and either of positivity properties:\ \
$
(\bfBS \Vert \bfBS)
\in\ \RCp[[v^{-1}]]
$
or
$
(\bfBS \Vert \bfBS)
\in\ \fra{1}{\nabla_e}\RCp[v^{-1}]
$.


\sff{F}
If $p=char(\k)>0$ then there is a canonical   isomorphisms
$
K^C
(\BBke)
\ @>\io_e>\cong>\
K^0[mod^{fg}(U^\hz_{\k,e},C)]
$
which sends
$\bfBe$ to classes of irreducible modules
and
$\nab_e^\sv\bfBS$
to classes of indecomposable projective
modules.

\sff{G}
If $p=char(\k)>0$ then
for $b_i\in\bfBS$, the evaluation  of the polynomial
$\del[\nabla_e (b_1 \Vert b_2)] $
 at  $1\in  \fGm$
is
equal to
the corresponding entry of the Cartan matrix of
$mod^{fg}(U^\hz_{\k, e},C)$, i.e., the dimension
of the $\Hom$ space between the corresponding indecomposable
projective objects.
The dimension of the $\Hom$ space
in category $mod^{fg}(U^\hz_{\k, e})$
is given by the evaluation  of the polynomial
$\nabla_e (b_1 \Vert b_2) $
at the point $(1,1)\in
C\tim\fGm$.

}

\sss{
Proof of Theorem \ref{LuLu} modulo
Proposition \ref{geometric formulation of Conjectures}
}
\lab{Lusztig's formulation of conjectures and  how it follows}
Regarding K-groups,
Lusztig considers the case $\k=\C$ and defines
$
\bfB_{\BB_e}^\pm
\sub K^\Ctil(\BBke)^\bee$ and
$
\bfB_{\La_e}^\pm
\sub K^\Ctil(\Stilkep)^\beS$
by the condition of
{\em asymptotic norm one}:\
$(b\Vert b)\in\ 1+v\inv\Z[[v\inv]]$.
Conjectures \cite[5.12]{Kth2} (a,b)
assert that these are
{\em signed} $R_\fGm$-bases.
Parts (c,d) of the Conjecture 5.12 say that
signed bases $\bfB_{\BB_e}^\pm,
\bfB_{\La_e}^\pm$ are asymptotically orthonormal
and
parts (e,f) of the conjecture say that the
two bases are dual.

Since the identifications of K-groups
in Proposition \ref{ptoCe} are easily
shown to be compatible
with involutions $\bee$, $\beS$ and pairing $(\, \Vert\, )$,
we see that if we define
$\bfB_{\BB_e}^\pm,\bfB_{\La_e}^\pm$ in the same way for all
$\k$, what we get will be independent of $\k$
and the same will  hold for validity of conjectures (a-f).
However,
if $\bst$ is known for some $\k$ then conjectures (a-f)
follow from   \sff{A} and \sff{B}.
The point is that
\sff{A}  implies that for this $\k$ one has
$\bfB_{\BB_e}^\pm=\bfBe^\k\sqcup\ -\bfBe^\k$
and
$\bfB_{\La_e}^\pm=\bfBS^\k\sqcup\ -\bfBS^\k$, so these are
indeed signed bases.\fttt{
\lab{standard argument}
This is a standard argument.
Let $\xi\in K^\Ctil(\Stilkep)$ with $(\xi\Vert \xi)\in 1+v\inv\Z[[v\inv]]$.
Write  $\xi$ as $\sum_{b\in\bfBS,\ n\in\Z}\
c_b^nv^n\ b$. If $N$ is the highest power of
$v$ that appears then $\sum_b\ (c_b^N)^2\ge 1$
and
asymptotic orthonormality of the basis
implies that $\sum_b\ (c_b^N)^2= 1$ and $N=0$.
If also $\beS \xi=\xi$ then $c_b^{-n}=c_b^n$ since
$\beS(vx)=v\inv\beS x$. Therefore,
$c^n_b\ne 0$ implies $n=0$ and such $b$ is unique.
}

Part (g) of the Conjecture says
 that $\del(\nabla_e b_1 \Vert b_2) \in \pm
\Z_{\geq 0}[(-v)^{-1}]$ for
$b_1,b_2\in \bfB^\pm_{\La_e}$\ie for
$b_1,b_2\in \bfBS^\k$.
Using the $X^*(C)$-action on
$\bfBS^\k$ we see that this is  equivalent
to the claim that for $b_i\in\bfBS^\k$ all coefficient polynomials
$c^\nu_{b_1,b_2},\ \nu\in
X^*(C)$\ of $\nab_e(b_1\Vert b_2)$
(see \sff{D} for notation),
are in $\pm \Zet_{\geq 0}[(-v)^{-1}]
$.
This follows from
the first claim in \sff{C} and the ``parity vanishing'' statement
\sff{D}.

Part (h) of the conjecture says that
$\bfBS^\k$ satisfies the normalization property from \sff{E}.

Thus part (1) of the Theorem is established.
For part (2) recall that
Conjecture 17.2 claims that
if a subcategory $\fM\sub
mod^{C,fg}(U_{\k,e})$
is a generic block then
there exists a bijection $
Irr(\fM)\con \bfB^\pm_{\La_e}/\{\pmo\}$
which is compatible with the action of $X^*(C)$
and identifies
the Cartan matrix of $mod^{C,fg}(U_{\k,e})$ with
the matrix
$|\del[\nab_e(b_1\Vert b_2)_{(1,-1)}]|,\  b_1,b_2\in\bfB^\pm_{\La_e}$,
of absolute values of evaluations at $(1,-1)\in C\tim\fGm=\Ctil$.

We know that the subcategory  $\fM=mod^{C,fg}(U^\hz_{\k,e})$
is a generic block
(see \cite{BG}). Also,
Proposition \ref{geometric formulation of Conjectures}
together with the established part (1) of the Theorem
yields
 bijections
$\bfB^\pm_{\La_e}/\{\pmo\}\cong
\bfBS\cong Irr(\fM)$.
Now the difference between
Lusztig's formulation and the second sentence in
\sff{G}
is that the former uses evaluation at
$-1\in\fGm$ and absolute value, while the latter uses evaluation
at $1\in\fGm$. This  is
accounted for by the
parity vanishing property $\sff{D}$. \epf

\sss{
Pairing
$(-\Vert-)$
and
Poincare series
of sheaves on $\Stilkep$
}
The next Lemma explains the categorical  meaning of the pairing
$(-\Vert-)$.
To present it we need another  notation.

Let $Rep^+(\Ctil)$
be the category of representations $U$ of $\Ctil$
with finite multiplicities and
with $\fGm$-isotypic components $U_d,\ d\in\Z$, vanishing for $d<<0$.
We denote by $[U]$ its image in the K-group
$K^0[Rep^+(\Ctil)]\cong\ \RC((v))$ where $v$ is the image of the
standard representation of $\fGm$ in the K-group.
This extends to a map $U\mm [U]$ from
$D^b[Rep^+(\Ctil)]$ to
$K^0[Rep^+(\Ctil)]$.



Now, for $\F\in  D^b\Coh^\Ctil (\Stilkep)$
we have
$\RGa(\F)\in D^b[Rep^+(\Ctil)]$
and it is easy to show
(see also \cite{Kth2}),
that $[\RGa(\F)]\in\RC((v))$ is Laurent series
of a rational function\ie it lies in
$\fK_\Ctil\subset \RC((v))$.
Of course, if $\GG\in  D^b\Coh^\fGm (\Stilkep)$
then the same applies to
$\RHom(\FF,\GG)=\ \RGa[\RHHom(\FF,\GG)]$.

Recall that  $P\mapsto
P^\sv
$ denotes the involution of $\fK_\Ctil$
corresponding to inversion on $\Ctil$.

\slem {\em
\label{spariva}
Let
$\F,\G\in  D^b(Coh^\Ctil(\Stilkep))$.

 a) If $\G$ is set theoretically supported
on $\BBkep$
and the class $[\G]$ is invariant
under $\bee$, then  $(\F\Vert \G)=\ [\RHom(\G,\F)]$.

b)   If  $[\G]$ is invariant
under $\beS$ then $(\F\Vert \G)=\ [\RHom(\F,\G)]^\sv$.

\lab{a sign}
c) Let
 $\sF,\sG\in D^b\Coh^\Ctil(\hBBke)$ be such
that the restrictions of $\sF$ (respectively, $\sG$)
to the formal neighborhood of
$\BBkep$ in
$\Stilkep$
is isomorphic
to the restriction of $\F$ (respectively, $\G$).

If  $[\G]$ is invariant
under $\beS$ then
$$
\nabla_e(\F\Vert\G)
\ =\
[\RHom(\sF,\sG)\Lotimes_{\O(\hatt e)}\k_e]^\sv
.$$

}

\medskip

The Lemma will be proven in section \ref{Proof of the Pairing Lemma.}.

\sus{
Proofs for subsection \ref{Gradings and bases in K-theory}
}
\lab{Proof of Theorem and Proposition}
\label{54}
In \ref{Three Involutions}
we recall $\bee,\ \beS,\ (\ \Vert\ )$
and in
\ref{Proof of the Pairing Lemma.}
we check formulas for $(b\Vert c)$ when $c$ is fixed by $\bee$ or $\beS$.
Then we prove in
\ref{Proof of ``fixed by involutions''}
that
$\bee$ fixes the K-class of (a shift of)
$\Ltil_{i,\la}$ and $\beS$ fixes the class of $\Etil_{i,\la}$.
This is all the preparation we need for the proof of Proposition
\ref{geometric formulation of Conjectures}
in
\ref{Proof of geometric formulation of Conjectures}.

\sss{
Involutions
$\bee,\beS$ and $\Upsilon$
}
\lab{Three Involutions}
Involutions
$\beS$ on $K^\Ctil(\Stilkep)$
and
$\bee$ on $K^\Ctil(\BBke)=K^\Ctil(\BBkep)$
are defined in \cite[section 5.11, page 304]{Kth2}
by
$$
\beS\dff (-v)^{-\dim\BB+2\dim\BB_e}  \Upsilon T_{w_0}^{-1} \D
\aand
\bee\dff (-v)^{-\dim \BB} \Upsilon\ci T_{w_0}^{-1}\ci \D
.$$
Here $\Ups$ is a certain involution,
$T_{w_0}$ is an element of
a
standard basis for the affine Hecke algebra
corresponding to the long element $w_0\in W$ and
$\D$ is the
Grothendieck duality functor.
Since the direct image for the closed embedding
$i:\BBkep\to \Stilkep$
intertwines Grothendieck duality functors, we have
$i_*\bee = v^{-2\dim \BB_e} \beS i_*
$, since $i_*$ is an embedding  we write this as
$
\beS=\ v^{2\dim\BB_e} \bee
$.

Actually,
$(-v)^{-\dim \BB} T_{w_0}^{-1}
$
is the effect on the K-group
of the action of $\tii{ w_0}\in\B$
on $D^b(\Coh^\Ctil(\Stilkep))$
(see  Theorem \ref{Action of braid group on  base changes}.b).
Therefore,
$
\bee=\ \Upsilon \tii{w_0} \D $.

The only information about $\Upsilon$
(defined in \cite[5.7]{Kth2})
that we will use
is as follows:
\begin{equation}
\label{Ups1}
\Upsilon = \sum_{s=1}^l
a_{s} g_s^*,
\ \tx{with}\
g_s\in A(\Ctil ,\k\cd e),  {\rm ord}(g_s)<\infty
\ \tx{and}\
a_s\in \Qu,
\sum_s\ a_s=1;
\end{equation}
\begin{equation}\label{Ups2}
\overline{\Upsilon}=\overline{T_{w_0}\circ \D}
.\end{equation}
Here $A(\Ctil,\k \cd e)$ is the group of automorphisms of
$G$ normalizing the line $\k \cdot e$ and $\Ctil$;
bar denotes the induced action on the Grothendieck
group $K^0(Coh^C(\Stilkep))$
and $T_{w_0}$ is the action of
$w_0\in W$ on $K(\Stilkep)$ from \cite{Kth2}.

Claim \eqref{Ups1} is immediate from the definition of $\Upsilon$
in \cite{Kth2} and \eqref{Ups2} will be shown in
the Appendix
\ref{Involutions on homology of Springer fibers}.

Lusztig defines pairing $(\ :\ )$ on $K^\Ctil(\Stilkep)$
by $([\F]\, : \,[\G])\dff
[\RGa(\F\Lotimes \G)]
$ (see \cite[2.6]{Kth2}),
and uses it to define the pairing
$(\ \Vert\ )
$
by
$$([\FF]\ \Vert\ [\GG])
\ \dff \
(-v)^{\dim\BB-2\dim\BB_e} ([\FF] : T_{w_0}\Ups [\GG])
=\
v^{-2\dim\BB_e} ([\FF] : \tii{w_0}\inv\Ups [\GG])
$$
(see \cite[5.8]{Kth2}).
Since
$\bee=\bee\inv=\ \D \tii{w_0}\inv\Upsilon $ gives
$\D\bee=\ \tii{w_0}\inv\Upsilon $, we have
$$
([\FF]\Vert [\GG])
=\
v^{-2\dim\BB_e} ([\FF] : \D\bee [\GG])
=\
([\FF] : \D\beS [\GG])
.$$

These  pairings on $K^\Ctil(\Stilkep)$
descend to
pairings on
$K^C(\Stilkep)$ which we denote the same way.
We will denote
$([\FF]\ \Vert\ [\GG])$ simply by $(\FF\Vert \GG)$.


\srem	
The involutions $\bee,\beS$ are K-group
avatars of dualities that would fix irreducibles
(resp. projectives)
corresponding to the fundamental alcove $\fA_0$
The point is that (if one neglects $\fGm$ equivariance),
the
duality $\RHHom(-,\OO)$ takes projectives for
$\fA_0$ to
projectives for  $-\fA_0$, and then $\tii{w_0}$
returns them to projectives for $\fA_0$.\fttt{
Notice that because of the difference between dimensions of supports,
the analogous procedure for irreducibles
would use  $\RHHom(-,\OO)[2\dim\BB_e]$ instead of
$\RHHom(-,\OO)$.
}
This composition creates a permutation of indecomposable projectives
or irreducibles for
$\fA_0$.
In order to undo this permutation
Lusztig uses  the centralizer action to
describe a $\Z[v,v\inv]$ {\em linear} involution  $\Ups$
on the K-group which induces the same permutation.
This is a generalization of the relation of the Chevalley involution
to duality for irreducible representations of a reductive group.


\sss{
Proof of Lemma \ref{spariva}
}
\lab{Proof of the Pairing Lemma.}
a)
$\bee$-invariance of $[\GG]$
gives
$$
(\FF\Vert \GG)
\ =\
([\FF]: (v)^{-2\dim\BB_e}\D\bee[\GG] )
\ =\
([\FF]: (v)^{-2\dim\BB_e}\D[\GG] )
.$$
However,
$$
\D[\G]
\ =\ \
[\RHHom(\G, \Om_{\Stilkep})[\dim\Stilkep]\ ]
\ =\ \
v^{2\dim\BB_e}\
[\RHHom(\G, \O_{\Stilkep})]
.$$
For the second equality recall that for the
standard symplectic form $\omega$ on $\tilde \N$,
restriction $\omega|_{\Stilkep}$
is again symplectic, so its top
wedge power $\omega^{\dim(\BB_e)}|_{\Stilkep}$
is a non-vanishing section of
the canonical
line bundle $\Omega_{\Stilkep}$.
	Now the claim follows since
$\om$ is invariant under the action
of $G$ and transforms by the tautological character under the action of
$\Gm$ by dilations, while 	
$\fGm$ acts by a combination of $G$ and the square of dilations.

Thus we see that
$$
(\FF\Vert \GG)
\ =\
(\FF:\RHHom(\G,\O) )
\ =\
[\RGa(\F\Lotimes \RHHom(\G,\O))]
\ =\
[\RHom(\G,\F)]
.$$

b)
$\beS$-invariance of $[\GG]$
gives
$$
(\FF\Vert \GG)
\ =\
([\FF]: \D\beS[\GG] )
\ =\
([\FF]: \D[\GG] )
\ =\
[\RGa(\F\Lotimes\D\GG)]
\ =\
[\RGa (\D\RHHom(\FF,\GG))]
.$$
It is easy to show that if $\Ctil$ acts linearly
on a vector space $V$ and $\fGm$ contracts it to the  origin
for $\fGm\ni t\to\yy$,
then for any
$\KK\in D^b[\Coh^\Ctil(V)]$
one has
$[\RGa(\D\KK)]=\ [\RGa(\KK)]^\sv$.
Applying this to
the sheaf $\KK=
(\Stilkep\to\Ske)_*\RHHom(\GG,\FF)$ on
the space $\Ske\cong Z_\fg(f)$
we get the result.

c)
For a finite dimensional $\Ctil$-module
$V$ Lusztig denotes
by $
\boxed{V}\dff
[\we^\bu V]\in R_\Ctil
$ the image of the super-module $\we^\bu V$
in $K^\Ctil$.
So,
for
$\KK\in D^b[Coh^\Ctil(V)]$ and $i:0\inj V$
a use of Koszul complex gives
$$
[i_*i^*\KK]=\ [\KK \Lotimes _{\O(V)}\ \O(V)\ten_\k\ \we^\bu V^*]
=\
\boxed{V^*}\
[\KK]
=\
\boxed{V}^\sv\
[\KK]
.$$
By definition
$\nabla_e\dff
\boxed{Z_\g(f)}\ \boxed{\fh}^{-1}
$ (\cite[3.1]{Kth2}), where $\Ctil$ acts on $\fh$ by $(c,t)h=c^{-2}h$.

As in
\ref{Construction of gradings},
there is a unique
$\Gc\in D^b[\Coh^\Ctil(\Stilke)]$
which agrees with $\sG$ on $\Stilke\cap \hBBke$.
Notice that because  the restrictions to
$\Stilkep\cap \hBBke$
agree for $\Gc,\ \sG$ and $\GG$, we also have
$\Gc|_\Stilkep\cong\GG$.
Now, in order to calculate the K-class of
$
\RHom(\sF,\sG)\ \ten_{\OO(\hatt e)}\ \k_e
\cong\
\RHom(\sF,\sG\ \ten_{\OO(\hatt e)}\ \k_e)
$,\
observe that by the definition of $\Gc$ we have
$
\sG\ \Lten_{ \OO(\hatt e) }\ \OO(\Ske\cap\hatt e)
\cong\
\Gc\ \Lten_{ \OO(\Ske)    }\ \OO(\Ske\cap\hatt e)
$,
and this gives
$
\sG\ \Lten_{\OO(\hatt e)}\ \k_e
\cong\
\Gc\ \Lten_{\OO(\Ske)}\ \k_e
$.
Similarly, $\sF$ gives $\Fc$ with analogous properties.
Therefore in $K^\Ctil=R_\Ctil$ one has
$$
[\RHom_\hBBke(\sF,\sG)\ \Lten_{\OO(\hatt e)}\ \k_e]
=\
[\RHom_\hBBke(\sF,\Gc\ \Lten_{\OO(\Ske)}\ \k_e)]
=\
[\RHom_{\hBBke\cap\Stilke}(\sF|_\Stilke,\Gc\ \Lten_{\OO(\Ske)}\ \k_e)]
.$$
When we replace
$\sF|_{\Stilke}$
with  $\Fc|_\hBBke$ we can view this as
$$
[\RHom_{\Stilke}(\Fc,\Gc\ \Lten_{\OO(\Ske)}\ \k_e)]
=\
\boxed{Z_\fg(f)^*}\ [\RHom_\Stilke(\Fc,\Gc)]
=\
\fra{\boxed{Z_\fg(f)^*} }{ \boxed{\fh^*} }\
[\RHom_{\Stilke}(\Fc,\Gc\ \Lten_{\OO(\fh)}\ \k_0)
.$$
Now,
$\Gc\ \Lten_{\OO(\fh)}\ \k_0
\cong\
\Gc|_\Stilkep\cong \GG
$.
The same observation for
$\Fc$ and  adjunction give
$$
=\
\nabla_e^\sv\
[\RHom_{\Stilke}(\Fc,\GG)]
=\
\nabla_e^\sv\
[\RHom_{\Stilkep}(\FF,\GG)]
.$$
So, the claim follows from b).

\sss{
Proof of invariance of the bases under the involutions}
\lab{Proof of ``fixed by involutions''}
{\bf ($\bbb\aleph$)
Proof of  \eqref{beta}.
}

{\bf
($\bbb\aleph$.i)
Reduction to:\ $\beS$ preserves
$\pl_{i,\la}\ \Qu[\Etil_{i,\la}^S]$.
}
According to \cite{Kth2},
the restriction of equivariance map induces
an isomorphism
$K^\Ctil(\Stilkep) /(v-1) \con K^C(\Stilkep)$.
The $\Qu$-vector subspace
$\pl_{i,\la}\ \Qu[\Etil_{i,\la}^S]$
in $K^\Ctil(\Stilkep)_\Qu$
maps isomorphically to $K^C(\Stilkep)_\Q$.
In view of \eqref{Ups2},  the action of $\beS$ on
$K^C(\Stilkep)_\Q$ is trivial, so it suffices
to see that the vector subspace
$\pl_{i,\la}\ \Qu[\Etil_{i,\la}^S]$
is invariant under
$
\beS
$.

We will factor $\beS$ into $\Ups\DD$,
a functor
$\DD\F\dff \tii{w_0}\big(\RHHom(\F,\O)\big)(4\dim\BB_e)$
and
$\Ups$ which is only defined on the K-group.
Indeed,
$\beS=v^{2\dim\BB_e}\bee
=
v^{2\dim\BB_e}\Ups \tii{w_0}\D
$
and
$\Om_\Stilkep=\OO_\Stilkep(2\dim\BB_e)$
(see the beginning of
\ref{Proof of the Pairing Lemma.}), so that
$\D\FF=\RHHom(\FF,\OO_\Stilkep)(2\dim\BB_e)[2\dim\BB_e]
$.
We will actually show that
$\pl_{i,\la}\ \Qu[\Etil_{i,\la}^S]$
is invariant under both
$\Upsilon$ and $
\DD$.

{\bf
($\bbb\aleph$.ii)
Invariance under $\Upsilon$.}
By \eqref{Ups1},
it suffices to show that for any
finite order element
$g$
in $A(\Ctil,\k \cd e)$,
the pull back $g^*$ permutes
$[\Etil_{i,\la}^S]$'s. Since $g$ is
 an automorphism
commuting with  the multiplicative
group
$\fGm$ and fixing
the line of $e$,
we see that $\{ g^*(\Etil_i)\}$  is a set of
$\fGm$-equivariant
vector bundles on $\hBBke$ satisfying
the properties of $\Etil_i$ from $\bst$,
thus
 uniqueness part of
$\bst$
implies that for each $i$ there
exists some $i^g$ such that
$g^*(\Etil_i)\cong \Etil_{i^g}(d_g)$, where integer $d_g$
does not depend on $i$.
Obviously $d_{g^n}=
nd_g$, and therefore
 $d_g=0$ since $g$ is assumed to have finite order.
The
isomorphism
$g^*(\Etil_i)\cong \Etil_{i^g}$ implies
$g^*(\Etil_i^S)\cong \Etil_{i^g}^S$
for the corresponding $\fGm$-equivariant bundles on
$\Stilkep$.
Since $g$ fixes
$\Ske$,\ $g^*$ fixes
$\fGm$-equivariant vector bundle
$\Etil_{i_0}=\OO_\Stilkep(2\dim\BB_e)$ and
by uniqueness in \ref{Normalizations}
this implies
that
$g^*$ permutes the collection of
$\Etil_{i}^S$'s and then also the collection of
all $\Etil_{i,\la}^S$'s.

Also note that
$\Upsilon$  fixes the K-class of
$\Etil_{i_0,0}=\OO_\Stilkep(2\dim\BB_e)$
since this is true for all relevant $g^*$
and
in \eqref{Ups1} we have
$\sum_s\ a_s=1$.

{\bf
($\bbb\aleph$.iii)
$\DD$
permutes $[\Etil_{i,\la}^S]$'s.
}
$\DD$ factors to $D^b[\Coh(\hBBkep)]$
as $\DD=\tii{w_0}\RHHom(-,\OO)$.
Part (e)
of the Theorem \ref{tstr_alcove}
shows that the dual vector bundles $\E_i^*$
are exactly all indecomposable
projectives in the heart of the t-structure
$\Tau^{\hBBkep,\Nt_\k}_{-\fA_0}$
on $D^b(\Coh(\hBBkep))$,
Since,
$b_{-\fA_0,\fA_0}=\tii{w_0}$
by example
\ref{braid group elements giving transition between alcoves},
part (a.2)
of the Theorem \ref{tstr_alcove} now  shows that
the sheaves $\DD\Etil_i=\ \tii{w_0}(\E_i^*)$
are all
indecomposable projectives in the heart of  $\Tau^{\hBBkep,\Nt_\k}_{\fA_0}$.
Thus
we have
$\DD\Etil_i=\ \tii{w_0}(\E_i^*)\cong\ \E_{\check{i}}$
for some permutation $i\mapsto \check{i}$ of the indexing set.\fttt{
This could also be deduced from \cite[Corollary 3.0.11]{BMR2}.
}

Let us now add  $\fG_m$-equivariance.
Since a $\fGm$-equivariant structure on $\Etil_i$
is unique up to a twist
(Lemma
\ref{Equivariant lifts of irreducibles and indecomposable projectives}.b),
we have
 $\DD\Etil_i\cong\
\Etil_{\check{i}}(d_i)
$
for some integers $d_i$.
 The  uniqueness statement in $\bst$ implies that
$d_i=d_j$ for all $i,j$.
On the other hand, it follows from \cite[5.14]{Kth2} that $\beS$
sends the class of
$\Etil_{i_0}^S\cong
\O_{\Stilkep}(2\dim \BB_e)$
to itself.
Since we have already checked that
$\Upsilon$ fixes $\O_{\Stilkep}(2\dim \BB_e)$ (the last remark in
($\bbb\aleph$.ii)),
we find
the same is true for $\DD$, therefore
$d_i=0$ for $i=i_0$ and then the same holds for all $i$'s.

We can transport
$\DD[\Etil_i]
=\
[\Etil_{\check{i}}]
$
to $\Stilkep$
to get
$
\DD[\Etil_i^S]
=
[\Etil_{\check i}^S]
$.
Similarly, uniqueness of a torus equivariant structure (up to a twist)
gives
$\DD[\Etil_{i,\la}^S]
\ \cong\
[\Etil_{\check{i},\nu(i,\la)}^S]
$
for some $\nu(i,\la)\in \La$.
We will write this as
$\DD[\Etil_{i,\la}^S]
\ \cong\
[\Etil_{(i,\la)^\sv}^S]
$

{\bf
($\bbb{\aleph\aleph}$)
Proof of \eqref{tiibeta}.
}
Recall from
($\bbb\aleph$.i)
that
$\bee
= v^{-2\dim \BB_e} \beS
=\
\Upsilon v^{-2\dim \BB_e} \DD
$.
In particular,
 $\bee$ acts
on
$
K^C(\Stilkep)$
the same as $\beS$\ie trivially.
Therefore, as in the proof of \eqref{beta}
we only need that
$\pl_{i,\la}\  v^{-2\dim \BB_e} \Q[\Ltil_{i,\la}]\ \sub K^\Ctil(\hBBke)$
be invariant under $\bee$,
and this will follow from more detailed information:\
$ v^{-2\dim \BB_e}[\Ltil_{i,\la}]
$
are permuted
by\
(i) finite order elements of $A(\Ctil,\k e)$
and
 (ii) $v^{- 2\dim \BB_e} \DD=\tii{w_0}\D$.

 Since we have checked that
 finite order elements of $A(\Ctil,\k e)$
permute
$\Etil_{i,\la}$'s it follows that
they also permute
$\Ltil_{i,\la}$'s, hence also $ v^{-2\dim \BB_e}[\Ltil_{i,\la}]$.
On the other hand,
$$
\RHom(\DD\tii\LL_{i,\la},\tii\EE_{j,\mu})
\cong\
\RHom(\DD\tii\EE_{j,\mu},\tii\LL_{i,\la})
\cong
\k^{ \de^{(i,\la)} _{(j,\mu)^\sv} }
\,$$
gives
$$
\k^{\de^{(i,\la)^\sv}_{(j,\mu)}}
\cong
\RHom(\DD\tii\LL_{i,\la},\tii\EE_{j,\mu})^*
\cong\
\RHom(\tii\EE_{j,\mu},\DD\tii\LL_{i,\la}\otimes \Omega_{\Stilkep}[2 \dim \BB_e])
,$$
hence
$
\DD(\tii\LL_{i,\la})
=
\tii\LL_{(i,\la)^\sv}\otimes \Omega_{\tii S_e'}^{-1}[-2 \dim \BB_e]
$.
Thus
$$
(v^{-2 \dim \BB_e}\DD)\big
(v^{-2 \dim \BB_e}  [\tii\LL_{i,\la}])
=\
\DD
[\tii\LL_{i,\la}]
=\
v^{-2 \dim \BB_e} [\tii\LL_{(i,\la)^\sv}]
.$$

\sss{
Proof of Proposition
\ref{geometric formulation of Conjectures}
}
\lab{Proof of geometric formulation of Conjectures}
\sff{A}
We know that
$\bfBS
$
and
$\bfBe$
are sets of representatives -- modulo $\fGm$ shifts --
of isomorphism classes of
respectively,
indecomposable projective objects in $\Cohe^\Ctil(\Stilkep)$
and
of irreducible objects in $\Cohe^\Ctil(\hBBkep)$.
Since the exotic t-structure is bounded,
they  form
bases in the respective
Grothendieck groups
$K^\Ctil(\Stilkep)$
and $K^\Ctil(\hBBkep)$ over the ring $R_\fGm$.
Pointwise invariance of $\bfBS
$
and
$\bfBe$
under the involutions $\beS$ and $\bee$
has been proved in the previous subsection
\ref{Proof of ``fixed by involutions''}.

Lemma  \ref{spariva}.b) implies that for
$b_i\in\bfBS$  one has
$(b_1\Vert b_2)=[\RHom(b_1,b_2)]^\sv$\
(because $\beS$ fixes
$b_2$).
Since $b_i$ are projective objects this is really
$[\Hom(b_1,b_2)]^\sv$, so it lies in $\RCp[v^\pmo]$.
Now, property
$\bst$
says that the algebra $\AA=End_{Coh(\Stilkep)} (\bigoplus\limits_i \Etil_i^S)$
equipped with
the
grading coming from the $\fGm$ equivariant structure
on $\Etil_i^S$
has no components of negative degree and the
component of degree zero is spanned by
identity endomorphisms of $\Etil_i^S$'s.
This is the same
as saying
that
if $X^*(C)b_1\ne X^*(C)b_2$ then
$[\RHom(b_1,b_2)]\in v\RCp[[v]]$ and
if $b_1=b_2$ then
$[\RHom(b_1,b_2)]\in  1+\RCp[[v]]$.
So, we have established
for the basis $\bfBS$
the
asymptotic orthonormality property and also a positivity
property
$(\bfBS||\bfBS)\sub \RCp[[v\inv]]$.

Similarly,
Lemma  \ref{spariva}.a) implies that for
$b_i\in\bfBe$  one has
$(b_1\Vert b_2)=\ [\RHom(b_1,b_2)]$,
because $\bee$ fixes
$b_2$.
The properties of the $\fGm$-grading of
$\AA$ imply that
the $\fGm$-grading on
$\Ext^\bu_\AA[\bigoplus\limits_i \Ltil_i^S,\bigoplus\limits_i \Ltil_i^S)]$
has no positive $\fGm$-degrees and the  component
 of degree zero is
spanned by
identity maps. This is the same as
as saying
that
if $X^*(C)b_1\ne X^*(C)b_2$ then
$[\RHom(b_1,b_2)]\in v\inv\RC[[v\inv]]$
and
if $b_1=b_2$ then
$[\RHom(b_1,b_2)]\in 1+v\inv\RC[[v\inv]]$.

\sff{B}
Since  $\beS$ fixes $\Etil_{j,\mu}^S$,
Lemma \ref{spariva}.b) and  Calabi-Yau property of $\gt$ give
$$
(\Ltil_{i,\la}\Vert \Etil_{j,\mu}^S)
=\ [RHom(\Ltil_{i,\la},\Etil_{j,\mu}^S)]
=\
[\RHom_\k
\left(\RHom(\Etil_{j,\mu}^S,\Ltil_{i,\la}\ten\Om_\Stilkep[\dim\Stilkep]),\k\right)]
$$
$$
=\
\big(
v^{2\dim\BB_e}
[ \RHom( \Etil_{j,\mu}^S,\Ltil_{i,\la})]
\big)^\sv
\ =\
v^{-2\dim\BB_e}
[\k^{\de^i_j\de^\la_\mu}]
\ =\
\de^i_j\de^\la_\mu
\ v^{-2\dim\BB_e}
.$$

\sff{C}
In \sff{A} we have already checked
that $(\bfBS\Vert \bfBS)\ \sub
\RCp[[v\inv]]$.
Recall that
$\Etil_{i,\la}^S$ was constructed so that on
$\Stilkep\cap\hBBke$ it coincides with  a certain projective
exotic object
$\Etil_{i,\la}\in \Cohe^\Ctil(\hBBke)$
(see \ref{Normalizations}).
So,
because $\beS$ fixes $\Etil_{j,\mu}^S$,
Lemma \ref{spariva}.c) gives
$$
\nabla_e(\Etil_{i,\la}^S\Vert \Etil_{j,\mu}^S)
\ =\
[\RHom(\Etil_{i,\la},\Etil_{j,\mu})\Lten_{\OO(\hatt e)}\k_e]
\ =\
[\RHom
(\Etil_{i,\la},
\Etil_{j,\mu} \Lten_{\OO(\hatt e)} \k_e)
]
.$$
Here, $\Etil_{j,\mu} \Lten_{\OO(\hatt e)} \k_e$
is exotic\ie
under the equivalence $D^b(Coh(\gt))\cong D^b(mod^{fg}(\Abig))$
(restricted to $\hBBke$),
 the object
$\Etil_{j,\mu}\Lotimes_{\O(\hatt e)}\k_e$
corresponds to a module
rather than a complex of modules.
The reason is that
the  algebra $\Abig$ is flat over $\O(\g)$
(see Lemma \ref{Algebras Abig and Az}), hence the same is true for its
projective modules.
 Therefore,
the result is just
$[\Hom (\Etil_{i,\la}, \Etil_{j,\mu} \Lten_{\OO(\hatt e)} \k_e)]
$ which  lies in $\RCp[v^\pmo]$.
However, as
$\nab_e\in 1+v\inv\RC[v\inv]$
(\cite[Lemma 3.2]{Kth2}), from
$(\bfBS\Vert \bfBS)\ \sub
\RCp[[v\inv]]$
we now get
$\nab_e(\bfBS\Vert \bfBS)\ \sub
\RCp[v\inv]$.

The second claim follows from asymptotic orthonormality
from \sff{A} and
the fact
$(K^\Ctil(\BBkep)\Vert K^\Ctil(\BBkep))\ \sub R_\Ctil$
which is checked  in \cite{Kth2}.

\sff{D}
Recall
from
\ref{Slodowy slice}
that
for a certain $\bz\in Z(G)$
the element $\bbb m=(\phi(-1)\bz,-1)$ of $\Ctil$
acts trivially on $\gt$.
This implies that it acts on any
$\Etil_{i,\la}$ by a scalar $\ep_{i,\la}$.\fttt{
If $G$ coincides with its adjoint quotient $\barr G$ then
$\bbb m^2=1$, so since
$\Etil_{i,\la}$ is indecomposable the claim is true and
$\ep_{i,\la}=\pmo$.
By definitions in
\ref{Construction of gradings}-\ref{Normalizations}
if $G$ is replaced by its adjoint quotient $\barr G$
the collections
$\EE_i,\ \Etil_i,\ i\in\II_e$,
do not change, and if  $\barr C$ is the image of $C$ in
$\barr G$ then
$B^{\barr G}\dff
\{\Etil_{i,\la};\ i\in\II_e,\ \la\in X^*(\barr C)\}
$ and the corresponding object for $G$ are related by
$B^G\cong B^{\barr G}\tim_{X^*(\barr C)}X^*(C)$.
This implies the general case.
We also see that
$\ep_{i,\la}=\ep_{i,0}\ \lb\la,\barr m\rb$.
}
For $\nu\in X^*(C)$ and $d\in\Z$
the coefficient of $v^d$ in
$c^\nu_{b_1,b_2}\in \Z[v^\pmo]$
is the dimension of
$\Hom_\Ctil[\Etil_{i,\la},\Etil_{j,\mu+\nu}(d)\Lten_{\OO(\fg)}\k_e]$.
If this is not zero then
$(-1)^d=\ep_{j,\mu+\nu}\ep_{i,\la}^{-1}$
since $\bbb m\in\Ctil$ acts on
$\Hom[\E_{i,\la},\Etil_{j,\mu+\nu}(d)\Lten_{\OO(\fg)}\k_e]$ by
$(-1)^d\ep_{j,\mu+\nu}\ep_{i,\la}^{-1}$.

\sff{E}
The normalization property
is a part of the definition of $\bfBS$ in \ref{Normalizations},
and we have already checked  that $\bfBS$ satisfies
all  other properties.
Any  $R_\fGm$-basis $\bfB$ of
$K^\Ctil(\Stilkep)$
which is pointwise fixed by $\beS$ and
satisfies asymptotic orthonormality is
of
the form
$\vap_b b,\ b\in\bfBS$,
for some $\vap\in\{\pmo\}$, this much was established
immediately after the statement of Proposition
\ref{geometric formulation of Conjectures}.
If $\bfB$ satisfies
normalization property
then
$\vap_{i_0}=1$.
Now either of positivity properties for the pairing
$(\ \Vert\ )$ implies
$\vap=1$. The reason is that
the equivalence relation
$\sim$
on $\bfBS$ generated
by $b_1\sim b_2$ if $(b_1 \Vert b_2 )\ne 0$, is transitive since
$(\Etil_{i,\la}^S \Vert \Etil_{j,\mu}^S )
\ne 0$ is equivalent to
$\Hom_{D^b[\Coh^\Ctil(\Stilkep)]}
(\Etil_{i,\la}^S \Vert \Etil_{j,\mu}^S )\ne 0$
and the category
$\Cohe^\C(\Stilkep)$ is indecomposable
(because $D^b\Coh(X)$
is indecomposable for a connected variety $X$
and
$D^b\Coh(\Stilkep)\cong D^b[mod^{fg}(A_S)]\cong D^b[\Cohe^C(\Stilkep)]$).
Notice also that the last claim is equivalent to indecomposability of
$\Cohe^C(\hBBke)$
and
then
the corresponding statement in representation theory is well known
(see \cite{BG}).


\sff{F}
The
equivalence
$
\Cohe^C
(\hBBke)
@>\io_e>>
mod^{fg}
(U^\hz_{\k,\hatt e},C)
$
from Theorem \ref{opportunity}(c) 
provides compatible bijections
(of isomorphism classes) of irreducibles
and indecomposable projectives
and
an isomorphism
$
K^C
(\hBBke)
@>\io_e>\cong>
K^0[mod^{fg}
(U_{\hatt e}^{\hatt 0},C)]
$
\
which we can view as
$
K ^C
(\BBke)
@>\io_e>\cong>
K^0[mod ^{fg}(U_e^\hz,C)]
$.
The list of irreducibles and their projective covers
in
$
\Cohe^C
(\hBBke)
$
is given by
images $\LL_{i,\la}
,\EE_{i,\la}
$
of
the corresponding objects
$\Ltil_{i,\la}
,\Etil_{i,\la}
$
of
$
\Cohe^\Ctil
(\hBBke)
$,
and we denote by
$L_{i,\la}
,E_{i,\la}
$
their images in
$mod^{fg}
(U^\hz_{\k,\hatt e},C)$.
The projective cover
of
$L_{i,\la}
$
in
$
mod^{fg}
(U^\hz_{\k,e},C)
$
is the restriction
$E_{i,\la}
\Lten_{\OO(\hatt e)}\k_e$.
So, it remains to notice that
the K-class of the restriction
$\Etil_{i,\la}
\ten_{\OO(\hatt e)}\k_e$
is
$\nabla_e^\sv[\Etil_{i,\la}
^S]$.
This calculation we repeat from
part (c) of \ref{Proof of the Pairing Lemma.}.
We use an intermediate object
$(\Etil_{i,\la})^o
\in \Cohe(\Stilke)$,
by its definition
$
\Etil_{i,\la}
|_{\Ske}
=\
(\Etil_{i,\la})
^o|_{\hBBke}
$, so one gets
$\Etil_{i,\la}
|_e
=\
(\Etil_{i,\la}
)^o|_e
$, hence
$[\Etil_{i,\la}
|_e\,]
=\
[(\Etil_{i,\la}
)^o|_e\,]
=\
\boxed{Z_\fg(f)^*}\,[(\Etil_{i,\la}
)^o]
$.
Also,
$
\Etil_{i,\la}
^S
=\
(\Etil_{i,\la})^o
\ten_{\OO(\fh)}\k_0
$
gives
$[\Etil_{i,\la}
^S]
=\
\boxed{\fh^*}\,[(\Etil_{i,\la})^o]
$, hence
$[\Etil_{i,\la}
|_e\,]
=\
\boxed{Z_\fg(f)^*}\ \boxed{\fh^*}\inv[\Etil_{i,\la}
^S]
=\
\nabla_e^\sv[\Etil_{i,\la}
^S]
$.

\sff{G}
In order  to avoid the dg-setting,\fttt{
The above localization of
$\Az_{\k,\hatt e}$-modules to coherent sheaves
on $\hBBkep$ specializes to a localization
of the category of
$\Az_{\k,e}$-modules on the Springer fiber, however one is forced to
use the dg-version
of the Springer fiber $\BBkep$
\cite{Simonchik}.
}
we will
pass here from
exotic sheaves to $\Abig$-modules
by means of
the equivalence
$
D^b[\Coh^T(\hBBke)]
@>\bAA>\cong>
D^b[mod^{T,fg}(\Abig_{\k,\hatt e})],\ \bAA\dff
\RHom(\EE|_{\hBBke},-)
$, where $T$ could be $\{1\}$, $C$ or $\Ctil$.

Let us start with the non-equivariant statement\ie
$\Hom$ in $U^\hz_{\k,\hatt e}$-modules.
We are  interested in
the composition of
equivalences
$$mod^{fg}(U^\hz_{\k,\hatt e})
\cong \Cohe(\hBBke)
@>\bF>\cong>
mod^{fg}(\Abig_{\k,\hatt e}).
$$
Due to compatibility with the action of
$\OO(\fg\tim_{\fh/W}\fh)$ it restricts to
an equivalence
$
mod^{fg}(U^\hz_{\k,e})
\cong\
mod^{fg}(\Abig_{\k,e})$.

We will start as in
\sff{C},  so
$\Etil_{j,\mu}\Lten_{\OO(\fg)}\k_e$ is an exotic sheaf and
$$
\nabla_e(\Etil_{i,\la}^S\Vert \Etil_{j,\mu}^S)
\ =\
[\RHom_{\hBBke}(\Etil_{i,\la},\Etil_{j,\mu}\Lten_{\OO(\fg)}\k_e)]
\ =\
[ \RHom_{\Abig_{\k,\hatt e}}
(\bF\Etil_{i,\la},\bF\Etil_{j,\mu}\Lten_{\OO(\fg)}\k_e)]
.$$
By adjunction in sheaves of $\Abig$-modules,
$$
\nabla_e(\Etil_{i,\la}^S\Vert \Etil_{j,\mu}^S)
\ =\
\RHom_{\Abig_{\k,e}}
(\bF\Etil_{i,\la}\Lten_{\OO(\fg)} \k_e,
\bF\Etil_{j,\mu} \Lten_{\OO(\fg)} \k_e)
]
.$$
So, the evaluation
$\nabla_e (\Etil_{i,\la}^S\Vert \Etil_{j,\mu}^S)(1_C,1_\fGm)
$
is the image of
$
\Hom_{\Abig_{\k,e}}
(\bF\Etil_{i,\la}\Lten_{\OO(\fg)} \k_e,
\bF\Etil_{j,\mu} \Lten_{\OO(\fg)} \k_e)
$
in $K^0(mod^{fg}(\k))$\ie
the dimension
of this vector space.

It remains
to notice that
$
\bF
\Etil_{i,\la} \ten_{\OO(\fg)} \k_e
$
is a projective cover of
$
\bF\Ltil_{i,\la}
$
in
$
mod^{fg}(\Abig_{\k,e})
$.
Since
$
\bF\Etil_{i,\la}
$
is a projective cover of
$\bF\Ltil_{i,\la} $
in
$mod^{fg}{\Abig_{\k,\hatt e}}$, it is projective over $\k$, hence
$$
\bF(\Etil_{i,\la} \Lten_{\OO(\fg)} \k_e)
\cong
\bF\Etil_{i,\la} \Lten_{\OO(\fg)} \k_e
\cong
\bF\Etil_{i,\la} \ten_{\OO(\fg)} \k_e
.$$
Since
$\AA\Ltil_{i,\la} $
is
irreducible
in
$mod(\Abig|_{\hatt e})$,
it is supported scheme theoretically on $e$, therefore
we find by adjunction that
$\bF\Etil_{i,\la} \ten_{\OO(\fg)} \k_e$
is a projective cover of
$\bF\Ltil_{i,\la} $
in
$
mod^{fg}{\Abig_{\k,e}}
$.

If one is interested in
maps in
$mod^{fg}(U^\hz_{e},C)$ only,
one uses
equivariant equivalences
$mod^{fg}(U^\hz_{\k,\hatt e},C)
\cong \Cohe^C(\hBBke)
@>\bF>\cong>
mod^{C,fg}(\Abig_{\k,\hatt e})
$ and one also
needs to
take  $C$-invariants in the above calculation.
This has the effect of  applying $\del$
to
$\nabla_e(\Etil_{i,\la}^S\Vert \Etil_{j,\mu}^S)$.
\qed

\subsection{
Koszul property
}
\label{Kos_conj}
\lab{Conjecture on Koszul property}
This subsection is not used in the rest of the text.
Set $\AA_e=
\End(\oplus \E_i^S)$ where the vector bundles $\E_i^S$ on $\Stilkep$
are as above.
The $\fGm$-equivariant structure $\tii \E_i^S$ on $\E_i^S$
introduced in \ref{Normalizations} equips $\AA_e$ with a grading.

\spro {\em
Properties
i$_\st$, ii$_\st$
of \ref{Property bst}
imply that the graded algebra
$\AA_e$ is a Koszul quadratic algebra.

}

\pf
 \fttt{The proof is due to Dmitry Kaledin.}
 For two graded modules $M$, $N$ over $\AA_e$ let $Ext^i_j(M,N)$
denote the component of inner degree $j$ in $Ext^i_{\AA_e}(M,N)$.
Then   i$_\st$, ii$_\st$ imply that $Ext^i_j(\tii \LL_1,\tii \LL_2)=0$ for $j<i$
where $\tii \LL_1$, $\tii \LL_2$ are irreducible
graded $\AA_e$-modules concentrated in graded degree zero.

The canonical line bundle of $\Stil_e'$ admits a trivialization which
transforms under the action of $\fGm$
by  the $2d_e$-th power of the tautological character.
So, Serre duality shows that for finite dimensional
graded $\AA_e$ modules we have
$$Ext^i_j(M,N)=Ext^{2d_e-i}_{2d_e-j}(N,M)^*.$$
Thus we see that $Ext^i_j(\tii \LL_1,\tii \LL_2)=0$
 for $j\ne i$, which is one of characterizations of Koszul algebras. \epf

\rem
For $e=0$ the work of  S.~Riche
\cite{Simonchik} provides a representation theoretic interpretation
of the algebra $\kappa(\AA_e)$
 which is Koszul dual to $\AA_e$. It would be interesting to
generalize this to nonzero nilpotents.

When $e$ is of principal Levi type, the relation between
the parabolic semi-infinite module over the affine Hecke algebra
and $K(Coh^{\tii C}(\Stil_e'))$ (see \cite{Kth2}, sections 9, 10) suggests that
the category of $\kappa(\AA_e)$-modules can be identified
with the category of perverse sheaves on the
parabolic semi-infinite flag variety
of the Langlands dual group. For $e=0$ this follows from the result of
\cite{5} compared with \cite{Simonchik}.

\se{
Grading that satisfies property
$\bst$
}
\label{sect6}

In  subsection \ref{subsCetop}
we reduce verification of property $\bst$
(see \ref{Property bst}),
to the case of a characteristic zero base field.
From then on until the end of
the  section we work over the field $\k=\Qlb$ of characteristic zero.

Our goal is to construct a $\fGm$-equivariant structure
on projective exotic sheaves that satisfies
property $\bst$.
For this we
use a derived equivalence between the category of
$G$-equivariant coherent sheaves on
$\Nt$ and of  certain perverse constructible sheaves on
the affine flag variety $\Fl$. In the new setting the
$\fGm$-structure is related to Frobenius (Weil) structure on $l$-adic
sheaves,
which we choose  to
be pure of weight zero.
In
\ref{Perverse t-structures on $\Az$-modules}
we compare the exotic t-structure on coherent
sheaves to the standard t-structure on perverse constructible sheaves on
$\Fl$, this involves the notion of
{\em perversely exotic} $G$-equivariant coherent sheaves.
In
\ref{Reduction to a $G_e$-equivariant setting}
we reduce $\bst$ to a property $\bsttt$ which is stated in terms of
$G$-equivariant sheaves.
Finally, in
\ref{End of the proof} we verify $\bsttt$.

\sus{
Lusztig's conjectures for $p\gg 0$
}\label{subsCetop}


\pro {\em 
If $\bst$
holds in characteristic zero,
it holds for almost all positive characteristics.

}            \pf
By Proposition \ref{gra}(c)
the choice of a graded lift of indecomposable projectives and
irreducibles in characteristic zero defines
such a choice in almost all prime characteristics.
We claim that the required properties are inherited from characteristic zero to almost all prime characteristics.
Indeed, the fact that the given choice of graded lifts satisfies the  positivity requirement amounts to vanishing of the components of negative
degree in the Hom space between indecomposable projective modules.
Since the sum of these components is easily seen to be a finite $R'$ module (here we use the fact that this Hom space is a finite module over the center
$\O(\Stilep_{R'})$), it vanishes after a finite localization provided that its base change to a characteristic zero field vanishes.
Invariance of the graded lifts under the action of the centralizer clearly holds in large positive characteristic if it holds in characteristic zero (notice that
the centralizer acts on the set of isomorphism classes of (graded) modules through its group of components, which is the same in almost all characteristics).

Uniqueness of the graded lifts with required properties amounts to non-vanishing of components of degree minus one in $\Ext^1$ between certain pairs of irreducibles
(see the proof of 6.2.1 below). After possibly replacing $R'$ with its localization we can assume that
$\Ext^1$ between the "extended irreducible" modules over $R'$
are flat over $R'$, thus dimensions of each graded component in
$\Ext^1$ between the corresponding irreducibles over every geometric
point of $R'$ is the same.
\qed

\sss{
The final form of the results
}
Since $\bst$ for $\k=\C$ will be established
in the remainder of this section, the proposition implies
that there
exists a quasifinite $R$-domain $R'$  such that
for all geometric points $\k$ of $R'$ property
$\bst$ holds
and therefore so do all claims
\sff{A}-\sff{G} from Proposition
\ref{geometric formulation of Conjectures}.
In particular this  establishes
the following version of Lusztig conjectures.

\stheo {\em 
\label{Lusztig Conjectures Theorem}
(1)
Conjectures  5.12, 5.16 of \cite{Kth2}
(existence of certain signed bases of K-groups of
complex schemes $\BBCe$
	and
$\StilCep$) hold.

(2)
The part of Conjecture 17.2 of {\em loc. cit.} concerning modular
representations holds for all
finite characteristic  geometric points $\k$ of $R'$.

}

\sus{
Perverse t-structures on $\Az$-modules
}
\lab{Perverse t-structures on $\Az$-modules}
Recall that the triangulated category
$D^b[Coh^G(\N)]$
carries a certain
t-structure called
{\em perverse coherent
t-structure of middle perversity}
\cite{CohPerverseNilpotent} (see also \cite{DimaRoma}
for the general theory of such t-structures).
As above let $\Az$ be $\End(\EE|_\Nt)$ for the vector bundle $\E$ from
Theorem \ref{tstr}.
This is an $\O(\N)$-algebra equipped with a
$G\times \fGm$-action.
This allows us to define
a perverse coherent
t-structure
$\TT^{\bG}_{pc}(\Az)$ of middle perversity on
$D^b[mod^{\bG,fg}(\Az)]$
where $\bG$ is one of the groups
$G,G\tim\Gm$ or $\fG$.
These are characterized by the requirement that
the forgetful functor to $D^b[Coh^G(\N)]$ is $t$-exact when
the target category is equipped with the
perverse coherent t-structure of middle perversity.

Recall the equivalence of derived categories of coherent and constructible
sheaves
\begin{equation}
\label{ABeq}
\PPhi:\
D^b(Coh^G(\Nt))
\iso
D^b(Perv_\Fl)
,\end{equation}
constructed in \cite{AB}, where $D^b(Perv_{\Fl})$
is the derived category of {\em
anti-spherical} perverse sheaves
on the affine flag variety
of the dual group.

\theo {\em
\label{ABtoM}
The composed equivalence $\PPhi_{\Az}$
$$
\PPhi_{\Az}\ \dff\ \
[\ D^b(mod^{G,fg}(\Az))
\con D^b(Coh^G(\Nt)) @>\PPhi>\cong>
D^b(Perv_{\Fl})
\ ]
.$$
sends the  perverse coherent t-structure of middle perversity
$\TT^{G}_{pc}(\Az)$
to the tautological t-structure on $D^b(Perv_{\Fl})$.

}            \pf
In  Lemma
\ref{Existence of Perversely exotic t-structures}
below, we show that the t-structure on
$D^b(Coh^G(\Nt))$
that comes
from
$D^b(Perv_{\Fl})$ satisfies a certain property
and in
Lemma
\ref{Uniqueness of Perversely exotic t-structures}
we show that the only t-structure on
$D^b(Coh^G(\Nt))$ that could
satisfy this property is the one coming from
the t-structure
$\TT^{G}_{pc}(\Az)$ on
$
D^b[mod^{G,\fg}(\Az)]
$.

\sss{
Perversely exotic t-structures
}
We will
say that a t-structure
on
a triangulated category $\CC$
is
{\em compatible with
a thick triangulated subcategory}
$\CC'$ if there
exist t-structures on $\CC'$, $\CC/\CC'$ such that the embedding
and projection functors are $t$-exact (cf. \cite{BBD}).
Inductively one extends
this definition to the definition of a t-structure compatible
with a filtration by thick
triangulated subcategories.

By a {\em support filtration}
on $D^b(Coh^G(\Nt))$ we will mean the filtration by
full subcategories
of complexes supported
(set theoretically) on the preimage of the closure of a given
$G$ orbit in $\N$
(we fix a complete order on the set of orbits compatible with the adjunction
partial order).

Finally, we say that
a t-structure on
$D^b(Coh^G(\Nt))$, is
{\em perversely exotic} if it is
\ben
compatible with the support filtration;

\i braid positive
(see \ref{Braid positive and exotic t-structures on T*(G/B)});
\i such that the functor $\pi_*$ is $t$-exact
when
the target category $D^b(Coh^G(\N))$ is equipped
with perverse coherent t-structure of middle perversity.
\een
Uniqueness of such t-structure follows from:

\lem {\em 
\lab{Uniqueness of Perversely exotic t-structures}
A perversely exotic
t-structure $\TT$ on
$D^b(Coh^G(\Nt))$ corresponds
under the equivalence $D^b(Coh^G(\Nt))\cong
D^b[mod^{G,fg}(\Az)]$
to
$\TT^{G}_{pc}(\Az)$,
the
perverse coherent t-structure of middle perversity.

}            \pf
It is a standard fact that for a triangulated category $\CC$, a thick subcategory $\CC'$
and t-structures $\TT'$ on $\CC'$, $\TT''$ on $\CC/\CC'$ a t-structure $\TT$ on $\CC$ compatible
with $\TT'$, $\TT''$ is unique if it exists. Thus uniqueness of an exotic t-structure implies uniqueness of perversely exotic t-structure.

On the other hand, the t-structure
 corresponding to $\TT^{G}_{pc}(\Az)$ is perversely exotic
as is clear from the fact that the t-structure
 corresponding to the tautological
one on $D^b(mod^{G,fg}(A))$ is exotic (the last fact is the definition of
$A$ as endomorphism of the exotic tilting generator
$\EE$). \qed

\lem {\em 
\lab{Existence of Perversely exotic t-structures}
The t-structure on $D^b(Coh^G(\Nt))$ which
under the equivalence  \cite{AB}
corresponds to the perverse t-structure
on $D^b(Perv_\Fl)$,\
is perversely exotic.

}            \pf
Properties (1) and   (3) are satisfied by \cite{AB}, Theorem 4(a) and Theorem 2
respectively.
We now deduce property (2)  from the results of
\cite{CohTiltMod}.
In {\em loc. cit.} it is shown that the t-structure corresponding to the one of $Perv_\Fl$
can be characterized
as follows:\
$$
D^{\ge 0}=\langle \Delta^\la[-d]\rangle_{d\ge 0,\ \la\in\La}
\aand
D^{\le 0}=\ \langle \nabla_\la[d]\rangle_{d\ge 0,\ \la\in\La}
,$$
where
$\langle\ ,\ \rangle$ denotes the full subcategory generated
under extensions and  $\Delta^\la$, $\nabla_\lambda$, $\la \in \La$,
 are certain explicitly defined objects in $D^b(Coh^G(\Nt))$.


Furthermore, we claim 
 that
\begin{equation}\label{nala}
\nab_\la=\tii w_\la (\OO),\ \ \  \Delta_\la =( \tii {w_\la^{-1}})^{-1}(\OO)
\end{equation}
where
 $w_\la$ is any 
 representative of the coset $\la W\subset W_{aff}$
 and $\tii w$ denotes
the canonical representative in $\Baff$ of $w\in W_{aff}$.
For $\la$ dominant we have $\nab_\la\cong \OO(\la)$, $\Delta_{-\la}
\cong \O(-\la)$, so in this case \eqref{nala} is clear
from the description of the action of $\theta_\la$ in Theorem \ref{132}(a,ii).
The general case follows from \cite[Proposition 7(b)]{CohTiltMod}:
in {\em loc. cit.} one finds a distiguished triangle 
$$\nab_\la\to \nab_{s_\al(\la)}\to F_\al'(\nab_{s_\alpha(\la)})$$
for a certain functor $F_\al'$, where we assume that $s_\al(\la)\preceq \la$. 
The functor $F_\al'$ is readily identified as convolution 
with $\F_\al[1]$, where $\F_\al = 
Ker(\OO_{\Ga_{s_\al}'}\to \OO_\Nt)\in Coh^\LG(\Nt^2)$
(notations of \ref{111}); here the arrow is the restriction to diagonal
map. The arrow  $\nab_{s_\al(\la)}\to F_\al'(\nab_{s_\alpha(\la)})$
appearing in \cite[Proposition 7]{CohTiltMod} coincides with 
the one coming from the map $\O_\Nt\to \F_\al[1]$ of
the distiguished triangle $\F_\al\to \OO_{\Ga_{s_\al}'}\to \OO_\Nt\to \F_\al[1]$.
This shows that $\tii s_\al(\nab_{s_\al(\la)})\cong \nab_\la$, which
implies \eqref{nala}. 

Now Proposition
\ref{quadrel}(a)
yields  exact triangles available for any $\F\in D^b(Coh^G(\Nt))$:
$$
\tii{ s_\al}^{-1} \F\to
\tii{s_\al}\F \to \F \oplus \F[1]
.$$
Thus if $\ell(s_\al w_\la)
<  \ell(w_\la)$
(where $\ell$ is the length function on $W_{aff}$),
then $\tii s_\al\inv\nab_\la=\nab_{s_\al(\la)}$, so
we have
an exact triangle:
$$
\nabla_{s_\al(\la)} \to \tii{s_\al}\nabla _\la\to
\nabla_\la \oplus \nabla_\la[1]
,$$
which shows that
$\tii{s_\al}\nabla _\la\in D^{\leq 0}$.
Also, if $\ell(s_\al w_\la)> \ell(w_\la)$, then
$\tii {s_\al}\nabla_\la\cong \nabla _{s_\al(\la)}$.
Thus $\tii s_\al: D^{\leq 0}
\to D^{\leq 0}$ which implies
braid positivity property (2).\qed

\rem
A more conceptual proof of
braid positivity property (2)
in the last lemma follows
from the preprint \cite{BHeck} (see also announcement in \cite{ICM}).
 It permits to relate the $\Baff$ action described above to a standard action on the
 category of constructible sheaves on the affine flag space $\Fl$. In the latter case
 the generator $\tii s_\al$ acts by convolution with a constructible sheaf $j_{s_\al*}$
 (in the notations of, say, \cite{AB}), i.e.
 the $*$ extension of the constant sheaf shifted by 1 on the  Iwahori orbit corresponding to $s_\al$.
 It is well
 known that convolution with such a sheaf is right exact with respect to the
 perverse t-structure (see e.g.
\cite{BeilinsonBernstein on Intertwining functors}).

 In fact, these considerations have led us to the notion of a braid positive t-structure,
 which was introduced as a way to relate modular representations to
 perverse sheaves on the affine flag space. We have chosen to present the above ad hoc argument
 in an attempt to keep the present paper self-contained.



\sus{
Reduction to a $G_e$-equivariant setting
}
\lab{Reduction to a $G_e$-equivariant setting}
We consider the algebra $\Az_e\dff \Az\otimes _{\O_\g} \k_e$.
It is graded by means of the action of $\fGm\sub\Ctil$
as above.


\sss{Reduction to a property of
$\Az_e$-modules
}

\slem {\em
Property
$\bst$ follows from:

$\bstt$
there exists a $G_e$-invariant
choice of a
graded lifting $\tii L$ for every
irreducible representation
$L$ of $\Az_e$,  such  that:

\ben
Components of
nonnegative
weight in
$\Ext^1_{mod(\Az_e)}(\tii L_1, \tii L_2)$
vanish for $L_i\in\ Irr(\Az_e)$.

\i
Consider the preorder on the quotient
of the set of irreducible representations of $\Az_e$
by the action of $G_e$,
generated by: $\al_1 \leq \al_2$
if for some representatives $L_i$ of $\al_i$
the component of degree\ --$1$  in
$\Ext^1_{mod(\Az_e)}(\tii L_1, \tii L_2)$ does not
vanish.
This preorder is actually a
transitive equivalence relation, i.e. $\al_1\leq \al_2$ for all
 $\al_i\in Irr(\Az_e)/G_e$.
\een

}            \pf
{\bf (i) Existence.}
Set $\Az_S=\Az\otimes _{\O_\g} \O(\Ske)$. Then we have
$D^b[mod^{fg}(\Az_S)]\cong D^b(Coh(\Stilke))$
and the same holds with equivariance under $\fGm$ or $\Ctil$.
Recall that  $\bst$
(see  \ref{Property bst},
involves a   $G_{e,h}$-invariant choice
of graded liftings $\Etil_i^S$ of exotic sheaves $\E_i^S$\ie
a $G_{e,h}$-invariant choice
of graded liftings of
indecomposable projective $\Az_S$-modules.
This is
equivalent to a
$G_{e,h}$-invariant choice of graded liftings $\tii L_i$
   of irreducible modules
$L_i$ supported at $e$.

Since $\bstt$ provides a choice with stronger  equivariance,
it remains to  check  that the choice of $\tii L_i$
satisfying the vanishing property (1)  of
$\bstt $ yields a choice of $\Etil_i^S$ satisfying the vanishing
requirements of
$\bst$\ie
$\pl_{i,j}\ \Hom_{\Coh(\Stilkep)}(\Etil_i^S, \Etil_j^S)$ has no negative
$\fGm$ weights and zero weights are spanned by identity maps.
By a standard argument this property from
$\bst$ is equivalent to saying that
$$
\Ext^1_{\Az_{S_e}}(\tii L_i,\tii L_j(d))=0 \ \  {\mathrm{for}} \  \
d\ge 0
.$$

If $L_i\not\cong L_j$ then
    any $\Az_S$-module which is an extension of $L_i$ by $L_j$
    is actually an $\Az_e$ module, because the action of a regular function on  $\Ske$ vanishing at $e$ on such an extension factors through a map $L_i\to L_j$, such map is necessarily zero.
    On the other hand, if $L_i\cong L_j$ and an extension $0\to L_j\to M\to L_i\to 0$ is such that $M$ does not factor through $\Az_e$, then some function as above induces  a nonzero map $L_i\to L_j$. Since $\fGm$ acts on the ideal in $\O(\Ske)$ by positive weights, we see that the class of the extension has negative weight.


{\bf (ii) Uniqueness.}
Finally, we will see that the uniqueness statement in
$\bst$ follows
from
property (2).
Any
graded lifting of $L_i$'s
is of the form  $\tii L_i(d_i)$ for some integers $d_i$.
If it  satisfies
the requirements, then $d_i$ is clearly monotone with respect
to our preorder.
  Thus property (2) implies that $d_i=d_j$ for all $i,j$. \qed


\sss{
Category $mod^{G_e,fg}(\Az_e)$ of $G_e$-equivariant $\Az_e$-modules.
}
Notice that  $mod^{G_e,fg}(\Az_e)\cong mod^{G,fg}(\Az|_{O_e})$, where
$O_e$ is the $G$-orbit of $e$ and the category in the right hand side is the
category of equivariant quasicoherent sheaves of modules over the sheaf
of algebras $\Az|_{O_e}$.
This category has  a graded version
$mod^{\fG,fg}(\Az|_{O_e})$, compatible with
the graded version
$mod^{\fG,fg}(\Az)\cong D^b[\Coh^\fG(\Nt)]$
considered above.
In terms of the stabilizer $\fG_e$ of $e$ in
$\fG$ this is
$mod^{\fG_e,fg}(\Az_e)$.

Tensor category $Rep(G_e)$
clearly acts on the category
$mod^{G_e,fg}(\Az_e)$ where for $V\in Rep(G_e)$ and
$M\in mod^{G_e,fg}(\Az_e)$
one equips
the
tensor product $V\otimes M$
with the diagonal action
of $G_e$.
We will now see that
a tensor subcategory
$Rep^{ss}(G_e)$
of semisimple representations  of $G_e$
acts on
$
mod^{\fG_e,fg}(\Az_e)
$.

We  use
morphisms $SL_2@>\varphi>> G$
and
$\fGm@>\bbi>>\fG,\ \bbi(t)= (\phi(t),t\inv)$,
chosen in
\ref{Slodowy slice}.
Notice  that $(g,t)\in\fG=G\tim\Gm$ lies in
$\fG_e$ iff
$e= t^2\cd {^g}e=\ ^{g\phi(t)}e$\ie $g\phi(t)\in G_e$.
So,  $\fG_e$ contains
$G_e\cd\bbi(\fG_m)$
and this is equality since $\tii g=g\phi(t)\in G_e$ implies that
$(g,t)=(\tii g,1)\cd \bbi(t\inv)$. We have exact sequence
$0\to\fG_m\sub \fG_e@>p>> G_e\to 0$ for $p(g,t)=g\phi(t)$ and
maximal reductive subgroups of $G_e$ and
$\fG_e$
can be chosen as the stabilizer $G_\varphi=Z_G(Im(\varphi))$ of $\varphi$
in $G$ and  $G_\varphi\cd\bbi(\fGm)$.
Now  $G_\varphi\cd \bbi(\fGm)@>p>>G_\varphi$
gives a tensor functor
$$
Rep^{ss}(G_e)
\cong\
Rep(G_\varphi)
@>p^*>>
Rep[G_\varphi\cd \bbi(\fGm)]
\conn
Rep^{ss}(\fG_e)
.$$

\slem {\em
(a)
For any $M\in \ Irr^{G_e}(\Az_e)$, restriction to $\Az_e$ is a multiple of a sum
over some $G_e$-orbit $\OO_M$ in $Irr(\Az_e)$.

(b)
For  $M_1, M_2\in\ Irr^{G_e}(\Az_e)$, the space
$\Hom_{\Az_e}(M_1,M_2)$ is a semi-simple $G_e$-module.

}            \pf
(a)
For an  irreducible $\Az_e$-module $L$  denote by
$G_{e,_L}\subset G_e$
the stabilizer of  the isomorphism class of $L$. Then
$G_{e,_L}$
is a finite
index subgroup in $G_e$, and $L$ can be equipped with a compatible projective
action of
$G_{e,_L}$\ie
an action of a central extension
$0\to \AA_L\to \GG_{e,_L}\to G_{e,_L}\to 0$
by a finite abelian group
$\AA_L$.\fttt{
V.~Ostrik has informed us that
he can prove that in fact this extension can be assumed to be
trivial provided $G$ is simply connected.
(An equivalent statement is that the set of
"centrally extended points" appearing in
\cite{BO} is actually a plain finite set with a $G_e$ action.)
We neither prove nor use this fact here.
}

The subgroup
$\AA_L$
acts
on
$L$ by a character $\chi_L$.
For any
irreducible
representation $\rho$
of $\GG_{e,L}$,
with
$\AA_L$
acting by $\chi_L^{-1}$,
 we get an irreducible
object $\rho\otimes L$ in $mod^{G_{e,L}}(\Az_e)$.

Then
$
Ind_{G_{e,L}}^{G_e}(\rho\otimes L)
$
is an irreducible object
of $mod^{G_e}(\Az_e)$, and all irreducible objects $M$ arise this way.

Now for
$M=Ind_{G_{e,L}}^{G_e}(\rho\otimes L)$,
we have
$$
M|_{\Az_e}
=\oplusl_{g\in G_e/G_{e,L}}\ ^g\rho\ten\ ^gL|_{\Az_e}
\cong\ \k^{\pl\dim(\rho)}\ten \oplusl_{K\in G_e\cd L}\ K
$$
To check (b) observe that
$\Hom_{A_e^0}(M_1,M_2)=0$ if $\OO_{M_1}\ne \OO_{M_2}$.
Thus assume that $\OO_{M_1}=\OO_{M_2}$ is the orbit of $L\in Irr(A_e^0))$, then we get:
$$
\oplusl_{g,h\in G_e/G_{e,L}} \Hom_\k[^g\rho,^h\rho]\ten
\Hom_{\Az_e}(^gL,^hL)
= \oplusl_{g\in G_e/G_{e,L}} \Hom_\k[\rho,\rho]
=\ Ind^{G_e}_{G_{e,L}} \End_\k[\rho]
.$$
Since $\rho$ is a semisimple $\GG_{e,L}$-module,
$\End_\k(\rho)$ is a semisimple
$G_{e,L} $ module.  As
$G_{e,L}$ has finite index in $G_e$,
it follows that
$
Ind^{G_e}_{G_{e,L}}\ \End_\k[\rho]
$ is also semisimple.


\sss{
Reduction to a property of
$G_e$ equivariant $\Az_e$-modules
}
\lab{Reformulation in Ge equivariant Ae-modules}
\slem {\em
$\bstt$ follows from the following.

$\bsttt$  there exists a choice of a graded lifting
$\tii L
\in mod^{\fG_e}(\Az_e) $
for every irreducible object
$L$ of $mod^{G_e}(\Az_e)$, such that:

	\bi
$\zsttt$ For
$L\in Irr^{G_e}(\Az_e)$ and
any irreducible representation $V$ of $G_e$ we have
$V\otimes \tii L\cong \sum \tii L_i$
for some
$L_i\in Irr^{G_e}(\Az_e)$.

	\i
$\osttt$
components of nonnegative weight in
$\Ext^1_{\Az_e}(\tii L_1, \tii L_2)$
vanish
for
$L_1,\, L_2\in Irr^{G_e}(\Az_e)$.

	\i
$\tsttt$
Consider the preorder $\leq$
on the set of irreducible objects in
$mod^{G_e}(\Az_e)$
generated by: $L_1\leq L_2$ if
component of degree $(- 1)$ in $\Ext^1_{\Az_e}(\tii L_1, \tii L_2)$
does not vanish.
This partial preorder is actually a
transitive equivalence relation, i.e. $L_1\leq L_2$ for all $(L_1, L_2)$.
	\ei

}            \pf
Property $\zsttt$ is equivalent to saying that
for $L_i\in Irr^{G_e}(\Az_e)$, the multiplicative group
$\fGm$ acts trivially on
$\Hom_{\Az_e}(\tii L_1,\tii L_2)$.
Indeed, if
$\zsttt$
 holds then
 $\fGm$ acts trivially
on
$\Hom_{\Az_e}(\tii L_1,\tii L_2)$,
because
for
any
$\si\in Irr(G_e)$,
$\Hom_{G_e}[\si,\Hom_{\Az_e}(L_1,L_2)]=\
\Hom_{mod^{G_e}(\Az_e)}(\sigma \otimes L_1,L_2)$
and $\sigma \ten \tii L_1$ is a sum of $\tii L$'s.

Consequently,
a choice of graded lifts of
$G_e$ equivariant irreducibles defines uniquely a choice of graded lifts
of irreducible $\Az_e$ modules, such that the forgetful functor sends the graded lift of
an equivariant irreducible $\tii M$ to
a sum of graded lifts of non-equivariant irreducibles $\tii L_i$.

It is clear that $G_e$ permutes the isomorphism classes of those $\tii L_i$.

Suppose that $\Ext^1(\tii L_1, \tii L_2(d))\ne 0$, for some $L_1,L_2$, $d>0$. Fix
$M_1, M_2$ such that $L_1$, $L_2$ are direct summands in $M_1$, $M_2$ considered
as $\Az_e$ modules. Then
 $\Ext^1_{\Az_e}(M_1, M_2(d))\ne 0$.
 The space  $\Ext^1_{\Az_e}(M_1, M_2)$ carries a (not necessarily semi-simple) $G_e$ action,
 and for an irreducible representation $\rho$ of $G_e$ we have an embedding of graded
 spaces $\Hom(\rho, \Ext^1_{\Az_e}(M_1, M_2))\to
  \Ext^1_{mod^{G_e}(\Az_e)}(\rho\otimes
 M_1,M_2)$. The latter embedding can be obtained as follows:
 given an element in the source space we get an extension of $\Az_e$
 modules $0\to  M_2 \to M \to \rho \otimes M_1\to 0$. Twisting this module
 by $z\in G_e$ we obtain an isomorphic extension; since
$\Hom_{\Az_e}(M_1,M_2)=0$
 we actually get  a unique isomorphism $M^z\cong M$ compatible with the given
 equivariant structures on  $\rho\otimes M_1$, $M_2$. Thus we get a $G_e$ equivariant
 structure on $M$.

 Since $\fGm$ acts on the Lie algebra of the unipotent  radical of $G_e$ by positive weights,
 the $G_e$ submodules generated by the degree $d$ components in $\Ext^1_{\Az_e}
 (M_1,M_2)$ is concentrated in positive degrees. This subspace has an irreducible
 subrepresentation $\rho$, which produces a nonzero $\Ext^1(\rho\otimes M_1,M_2)$
 of positive degree contradicting properties
$\zsttt$,$\osttt$.
 To prove property
$\tsttt$
it is enough to show that if $M_1$, $M_2$ are irreducible
 objects in $mod^{G_e}(\Az_e)$ such that
$\Hom_{\Az_e}(M_1, M_2)=0$, and
 $\Ext^1_{mod^{ \fG}(\Az_e)}
 (\tii M_1,\tii M_2(1))\ne 0$, then
$\Ext^1_{mod^{\fGm}(\Az_e)}(\tii M_1,\tii M_2(1))\ne 0$.
 It suffices to check that applying the forgetful functor
$mod^{G_e}(\Az_e)
 \to mod(\Az_e)$ to a nontrivial  extension
$0\to M_2\to M\to M_1\to 0$ we get a nontrivial extension. However, if there exists
an $\Az_e$ invariant splitting $M_1\to M$, then its image has to be invariant under $G_e$,
since $\Hom_{\Az_e}(M_1,M_2)=0$ and the isomorphism class of the
$\Az_e$-module $M_1$
is $G_e$ invariant. Thus existence of a non-equivariant splitting implies
the existence of an equivariant splitting. \qed

\sus{
End of the proof
}
\lab{End of the proof}
Here we prove $\bsttt$, the proof is based on the equivalence
$\PPhi:
D^b[\Coh^G(\Nt)]\iso
D^b(Perv_\Fl))$
from \cite{AB},
which we use in the form
$\PPhi_{\Az}:
D^b[mod^{G,fg}(\Az)]\iso
D^b(Perv_\Fl)$
(see  Theorem \ref{ABtoM}).

\sss{
The choice of grading
}
\lab{Translation between perverse coherent and constructible}

Equivalence $\PPhi_{\Az}$ makes category
$mod^{G_e,fg}(\Az_e)=Coh^G(\Az|_{O_e})$
a full subcategory
in a  Serre quotient category of $Perv_\Fl$.
We will now
show that
property
$\bsttt$ holds when the
graded lifting  $\tii L$
of irreducibles $L$ in $mod^{G_e}(\Az_e)$ is chosen so that it
corresponds to pure Weil structure of weight zero.
What is meant by this  is the following.

First,
it is shown
in \cite{AB} that the
Frobenius functor corresponding to a finite field $\Fq$
on the perverse sheaves category, corresponds
to the functor $\GG\mapsto {\bf q}^*(\GG)$ on coherent sheaves,
where
${\bf q}:\Nt@>>>\Nt$ by ${\bf q}(\b,x)=(\b, qx)$.
The same then applies to $\F\in D^b[mod^{G,fg}(\Az)]$ with
${\bf q}:\N@>>>\N$ by ${\bf q}(x)=q\cd x$.

Thus, for a perverse coherent sheaf $\F$ of $\Az$-modules,
a Weil structure on the perverse sheaf
$\PPhi_{\Az}\F$ is the same as an isomorphism
$\F\con {\bf q}^*(\F)$.
In particular this shows
that  any $\fGm$-equivariant structure on $\F$ defines
a Weil structure on $\PPhi_{\Az}\F$.
Notice that the resulting functor
from $D^b(Coh^{G\times \Gm}(\Nt))$ to Weil complexes on $\Fl$
sends the twist by the tautological $\Gm$ character $M\mapsto M(1)$
to the square root of Weil twist $\F\mapsto F(\frac{1}{2})$
acting on Weil sheaves.

$M(1)$ here stands for the graded module
$M(1)^i=M^{i+1}$. This functor is compatible
 with the functor $\F\mapsto \F(\frac{1}{2})$ on
Weil perverse sheaves under the equivalence
  \eqref{ABeq}.

It is shown in \cite{CohTiltMod} that when $\PPhi_{\Az}\F$ is an irreducible
perverse sheaf, any $\fGm$-equivariant structure on $\F$ induces a pure
Weil structure on $\PPhi_{\Az}\F$ and there is a unique $\fGm$-equivariant
structure on $\F$ such that the corresponding Weil structure on
$\PPhi_{\Az}\F$ is pure of weight zero.\fttt{
\cite{CohTiltMod}
  provides also a more direct way to describe the resulting $G\times \fGm$
equivariant sheaves.
}
It is also proven in {\em loc. cit.} that for $\F,\G\in D^b(Coh^{G\times
\fGm}(\Nt))$ the isomorphism
$$
\Hom_{D^b(Coh^G(\Nt))}(\F,\G)
\ \cong\
\Hom _{D^b(Perv_\Fl)}(\PPhi(\F),\PPhi(\G))
$$
takes the grading induced by the $\fGm$-equivariant structure into
the  grading by Frobenius weights.

\sss{
Property $\osttt$ and  Purity Theorem
}
Purity Theorem of \cite{BBD}
implies that
$\Ext^1$ between two
pure weight zero Weil sheaves in $Perv_\Fl$
has weights $< 0$.
Thus  for the above graded lifts $\tii L$
of
irreducible
equivariant perverse coherent
sheaves
of $\Az$-modules,
$\Ext^1$  between two such objects
has weights $< 0$.
It is not hard to check that this property
is inherited by a quotient category, thus
property
$\osttt$
follows.

\sss{
Property
$\tsttt$ and definition of cells
}
\lab{Property tsttt and definition of cells}
Property
$\tsttt$
says that
for any
$L$, $L'$ in $Irr^{G_e}(\Az_e)$,
there exists a sequence of
irreducible objects
$L_0=L, L_1,\dots, L_n=L'$ such that
the component of degree $-1$ in
$\Ext^1_{mod^{G_e}(\Az_e)}(L_{i-1},L_i)$ is nontrivial.

Recall that by \cite[Theorem 4(a)]{AB}
the support filtration on $D^b(Coh^G(\Nt))$
is identified with the (left) cell filtration on $D^b(Perv_\Fl)$.
In particular the irreducible
objects in the subquotient piece of the
filtration corresponding to a given nilpotent
orbit $O_e$ are in bijection with elements in a
{\em canonical left cell}
in the two-sided cell in $\Waff$ attached to $O_e$.
Furthermore, the definition
of a left cell implies the following.
For any two
irreducible objects
$L$, $L'$
in the same left cell,
there exists a sequence of
irreducible objects
$L_0=L, L_1,\dots, L_n=L'$ such that for any step
$M=L_{i-1}$ and $N=L_i$ in the chain, there is a simple
affine root $\al$ such that
$N$ is a direct summand  in the
perverse sheaf
$\pi^*\pi_*M[1]$,
where $\pi$ stands for the projection
$\Fl\to \Fl_\al$
to the partial flag variety of the corresponding type.

This implies that  $\pi_* M$ is a
semisimple perverse sheaf on $\Fl_{\al}$
and that  the relation of $\al$ to $M$
is such that we have
a canonical extension  of Weil perverse sheaves
 \begin{equation}
 \label{alpha small}
0\to \pi^*\pi_*
M[1](\frac{1}{2} )
\ \to\
\F\ \to\
M\to 0
.\end{equation}
Here $\F=\ J_{s_\al}^*\star M$
where\
$\star$\ denotes the  convolution
of constructible sheaves on $\Fl$ and
$J_{s_\al}^*$
is the $*$ extension
of the (pure weight zero perverse)
constant sheaf on the Schubert cell corresponding to $s_\al$.


So, it suffices to  see that in each of the above steps
the
component of degree $-1$ in $\Ext^1_{Perv_\Fl^e}(M,N)$ is nontrivial, where
$Perv_\Fl^e
\dff
Perv_\Fl/Perv_\Fl^{<e}
$
for the Serre
subcategory
$Perv_\Fl^{<e}$
generated by irreducible objects belonging to smaller cells.
Exact sequence \eqref{alpha small}
gives
$$
\Hom_{Perv_\Fl^e}
[\F,N]
@>>>
\Hom_{Perv_\Fl^e}
[\pi^*\pi_*M[1](\ha),N]
@>>>
\Ext^1_{Perv_\Fl^e}
[M,N]
.$$
The middle term is nonzero since
$N$ is a summand of $\pi^*\pi_*M[1]$. It
has weight $-1$ because $M,N$ are pure of weight zero,
hence
$\pi^*\pi_*M[1]
$
and
$\Hom_{Perv_\Fl^e}
[\pi^*\pi_*M[1],N]
$
are
also pure of weight zero.
So it suffices
to see that
\begin{equation}
 \label{homnul}
  \Hom_{Perv_\Fl^e} (\F, N)=0
.\end{equation}

To check \eqref{homnul} notice that  a nonzero element of the Hom space
  corresponds to a quotient $\F'$ of $\F$ in $Perv_\Fl$, such that
the only irreducible
  constituent of $\F'$ which does not belong to
$Perv_{<e}(\Fl)$ is $N$. Since $M$ is not in
$Perv_{<e}(\Fl)$,
the
  exact sequence
\eqref{alpha small}
shows that such quotient $\FF'$ is necessarily
of the form $\pi^*\F''[1]$
for some semi-simple perverse sheaf $\F''$ on the partial affine flag variety
$\Fl_\al$.
We have:
$$
\Hom (\F, \pi^*\F''[1])
=\
\Hom(\F, \pi^!\F''[-1])
=\
\Hom(\pi_*\F, \F''[-1])
$$
$$
=\
  \Hom(\pi_*M[1], \F''[-1])
,$$
  where we used the identity
$
\pi_*(J_{s_\al}^*\star \GG)
=\pi_*\GG[1]
$\ for $\GG=M$.
Finally,
since $\pi_*M$ and  $\F'$
  are perverse sheaves
$\Hom(\pi_*M[1], \F''[-1])
=\Ext^{-2}
(\pi_*M, \F'')=0$.

\medskip

\sss{
Property $\zsttt$ and Gabber's theorem
}

Property $\zsttt$
claims that the class of semisimple objects of
$mod^{\fG_e,fg}(\Az_e)$
whose irreducible constituents are the particular
lifts
$\tii L$
(chosen in
\ref{Translation between perverse coherent and constructible})
of irreducibles
$L$ in $mod^{G_e}(\Az_e)$,
is invariant under the action of
$Rep^{ss}(G_e)$.

It is explained in \cite{AB}
that under the equivalence $
\PPhi:\ D^b(Coh^G(\Nt))\iso
D^b(Perv_\Fl)
$,
the action $\GG\mm V\ten\GG$ of $V\in Rep(G)$ on
the source,
corresponds on the target  to the action of
a {\em central functor}  $\ZZ_V$ described in  \cite{KGB}.
This is then also true for the
equivalence $\PPhi_{\Az}: D^b(mod^{G,fg}(\Az))\iso
D^b(Perv_\Fl)
$.
Since the central functors are defined by means of a nearby cycles
functor, thus they carry the canonical {\em monodromy} automorphism $\fM$.

To
any $M\in mod^{G_e,fg}(\Az_e)$
one associates a
$G$-equivariant vector bundle
$\MM$ on the nilpotent orbit
$O_e$
and its intersection cohomology extension
$\IC(\MM)$ which lies in the
heart of the perverse
t-structure
of middle perversity on $D^b(mod^{G,fg}(\Az))$
and has support
$\barr{O_e}$ (see \cite{DimaRoma}). We will denote $\IC(\MM)$ just by $\IC(M)$,
then  $M\mm \IC(M)$ is a bijection of irreducibles in
$mod^{G_e,fg}(\Az_e)$
and those   irreducibles in the heart of the perverse
t-structure
that have  support
$\barr{O_e}$ (ibid).

For $V\in Rep(G)$ we have
$V\otimes \IC(M)=\IC(V|_{G_e}\ten M)$.
Moreover,
for any semisimple subquotient $\rho$ of $V|_{G_e}$, the tensor product
$\rho\ten M$ is semisimple, so
$\IC({\rho\otimes M})$ is semi-simple.
It is also a subquotient of $V\otimes \IC(M)$ (ibid).

The $G_e$-module $V|_{G_e}$ carries a nilpotent endomorphism given by the
action of $e$ and we denote by
$F^i(V)$ the corresponding Jacobson-Morozov-Deligne
filtration, and
$gr_i(V)=F^i(V)/F^{i+1}(V)$.
By definition of this filtration the graded pieces
$gr_i(V)$ are semisimple $G_e$-modules,
thus $\IC({gr_i(V)\otimes M})$ is a semisimple subquotient of
$V\otimes \IC(M)$.

Now we pass to $\fG_m$-equivariant objects. By the same formalism,
if we start with
$\tii M \in mod^{\fG_e,fg}(\Az_e)$
with the underlying object $M$ in $mod^{G_e,fg}(\Az_e)$,
we get  a graded
lift
$\IC({\tii M})$
of $\IC(M)$ that lies in the perverse heart of
$D^b(mod^{G\times \fGm,fg}(\Az))$.
As was explained in
\ref{Translation between perverse coherent and constructible},
the $\fGm$-equivariant
structure $\tii M$
induces a
Weil structure on the  perverse sheaf
$\PPhi_{\Az}
\big(
\IC(M)
\big)$;\
we will denote the corresponding Weil sheaf by
$\PPhi_{\Az}\big(\IC(\tii M)\big)$.
We will combine this with the action
of semisimple representations
$\rho$ of $G_e$ on $mod^{\fG_e,fg}(\Az_e)$ in order to produce
Weil sheaves
$\PPhi_{\Az}\big(\IC(\rho\otimes \tii M)\big)$.
Now property
$\zsttt$ is the part b) of the following Lemma.

\slem {\em
Let
$\tii M \in mod^{\fG_e,fg}(\Az_e)$ be such that the Weil structure on
$\PPhi_{\Az}\big(\IC(\tii M)\big)$ is pure of weight zero.

a)
For any $V\in Rep(G)$,
the Weil structure on
$\PPhi_{\Az}\big(\IC({ gr_i(V)\otimes \tii M })\big)$
is pure of weight $i$.

b)
For any semisimple representation $\rho$ of $G_e$
the Weil structure on
$\PPhi_{\Az}\big(\IC({ \rho\otimes \tii M })\big)$
is pure of weight zero.

}            \pf
a)
We consider the
nilpotent endomorphism $\ud{e}$ of the $G_e$-module
$V|_{G_e}\otimes M$  given by the action of $e$ on
$V|_{G_e}$.
It induces a  nilpotent endomorphism
of $V\ten \IC(M)=
\IC(V\ten M)$
which can be used to define a
Deligne-Jacobson-Morozov filtration on
$V\ten \IC(M)$.
The induced filtration on the fiber
$V\otimes \IC(M)|_e\cong V|_{G_e}\otimes M$
is just the Deligne-Jacobson-Morozov filtration
for $\ud{e}$
because formation
of Deligne-Jacobson-Morozov filtration commutes with exact functors
and
restriction to $e\in \barr{O_e}$ is exact on perverse sheaves supported
in $\barr{O_e}$,
Thus the
semi-simple subquotient
$\IC({ gr_i(V)\otimes \tii M })$
of $V\ten \IC(\tii M)$
is actually a subquotient of $gr_i(V\ten \IC(\tii M))$.

According to \cite{AB}, $\ud e$
induces on
$$
\PPhi_{\Az}\big(\IC({V|_{G_e}\otimes \M})\big)
=\ \PPhi_{\Az}\big(V\ten \IC(M)\big)
=\
\ZZ_V\big(\PPhi_{\Az}\ \IC(M)\big)
$$
the endomorphism given by the action
of the logarithm of monodromy $\log \fM$ on the functor
${\ZZ_V}$.
Now the Lemma follows from
Gabber's Theorem asserting that the monodromy filtration (i.e., the
Deligne-Jacobson-Morozov filtration for the logarithm of monodromy)
coincides with the
weight filtration
on the nearby cycles of a pure weight zero sheaf, cf. \cite{BeiBer}.

b)
Any irreducible representation $\rho$ of $G_e$ is a subquotient of
$V|_{G_e}$ for some
$V\in Rep(G)$, hence a subquotient of some
$gr_i(V|_{\fG_e})$.
The definition  of Deligne-Jacobson-Morozov filtration implies
that the natural
 $\Gm$ action on $gr_i(V|_{\fG_e})$ is by the character $t\mapsto t^i$,
so part a)     implies that the Weil sheaf
$\IC({\rho\otimes \tii M})(i)$
is a subquotient in $gr_i(V\ten \IC(\tii M))$.
Thus $\PPhi_{\Az}\big( \IC({\rho\otimes \tii M})(i)\big)$ has weight $i$ and
then $\PPhi_{\Az}\big(\IC({\rho\otimes \tii M})\big)$ has weight zero.
\qed

\appendix
\se{
Involutions on homology of Springer fibers
}
\lab{Involutions on homology of Springer fibers}

Our goal here is to prove equality
\eqref{Ups2} from
\ref{Three Involutions}.
The result can be viewed as a generalization of the fact
that a Chevalley involution (i.e. an involution
which sends every element of some Cartan subgroup to its inverse)
sends every irreducible representation
of an algebraic group to its dual.

%


\sus{
Cohomology of a Springer fiber as a module
for the extended centralizer
}

All cohomology spaces in this subsection are taken with coefficients
in $\Ce$ in the classical topology
 or coefficients in $\Qlb$ in the $l$-adic setting.

Let $\iota$ be an involution of $G$ which induces conjugation
with $w_0$ on the abstract Weyl group (e.g. a Chevalley involution).
Let $\Gbu$ denote the semi-direct product
$\{1,\io\} \ltimes G$.
It is well known that $\iota$ as above is unique up to composition with an
inner automorphism, thus the group $\Gbu$ is defined uniquely
 (up to an isomorphism).

Fix a nilpotent $e\in \g$
and set $d_e=\dim (\BB_e)$. Let $G_e$ be the centralizer of $e$ in $G$
and $\Ztil$ be the stabilizer of $e$ in $\Gbu$. Set $\Gamma=\pi_0(G_e)$
and
$\Gatil=\pi_0(\Ztil)$.
It is easy to see that $\Ztil$ intersects the non-identity component of
$\Gbu$, thus $\Gatil/\Ga\cong \Zet_2$.
 Let $\eps$ be the nontrivial character of
 $\Gatil/\Ga$.

The group $\Ztil$ acts on the Springer fiber $\BBke$ thus
$\Gatil$ acts on its cohomology.
We denote this action by $\eta$.
 We consider also another
action of $\Gatil$ on $H^*(\BBke)$: the two actions coincide on the subgroup
$\Gamma\sub \Gatil$, while on elements of
$\Gatil\setminus \Ga$
they differ by the action of $w_0\in W$ (where $W$ acts via the
Springer representation).
We denote
 this new action of $\Gatil$ on $H^{\bu}(\BBke)$
by $\psi$.

Notice that unlike the original action, $\psi$
commutes with the action of $W$ in all cases.


\spro {\em
\label{signe}
Let $\rho$ be an irreducible constituent of
the $\Gatil$-module $(H^{2i}(\BBke), \psi)$.
Then $\rho \otimes \eps^{d_e-i}$ is a constituent of
$ (H^{2d_e}(\BBke),\psi)$.

}


\srems
(1) Validity of the proposition for the
groups such that $w_0$ is central in $W$
is equivalent to the result of \cite{Sp}.
The method of \cite{Sp} is
based on Shoji's orthogonality formula for Green functions and is quite
different from the present one.

(2) After the paper has been submitted we have learned of a recent
preprint \cite{Kato} where it is shown that homology of a Springer
fiber is generated by its top degree component as a module
over cohomology of the flag variety. This result yields an alternative
proof of Proposition \ref{signe}.

            \pf
It is well-known that any irreducible representation of $\Ga$
which occurs in $H^i(\BBke)$
for some $i$ occurs also in $H^{d_e}(\BBke)$.
Thus the Proposition follows from the following

\slem {\em
The extension $\Gatil$ acts on $(H^{2i}\otimes H^{2j})^\Ga$
by
the character $\eps^{\ten\ i+j}$.

}            \pf
We will deduce the lemma  from some known properties
of equivariant Borel-Moore homology of the Steinberg variety of triples
$St\dff\ \Nt\times_\g \Nt$.
Let $H_{BM}^\bu$ denote Borel-Moore homology, i.e. derived global sections
of the Verdier dualizing sheaf
(for convenience
we use cohomological grading despite the term ``homology'').\fttt{
Since
the Verdier dualizing sheaf admits a canonical lifting to the equivariant
derived category, equivariant Borel-Moore homology is also defined
(cf. \cite[1.1]{Lgr} for a slightly more elementary definition).
}

It is well known (see e.g. \cite{Lgr} Corollary 6.4 for a much stronger
result) that
$$H_{BM}^{2i,G}(St)\cong \Ce[W]\otimes Sym^{i+d}(\fh^*),$$
where $d=2\dim \BB$ and odd degree homology vanishes.

On the right hand side of the last isomorphism we have a natural action
of $W$ (by conjugation on the first
factor and by the reflection representation on $\fh^*$)
 and of the group of outer automorphisms of $G$.
Standard considerations show that the automorphism
$\iota\circ w_0$ acts trivially on $\Ce[W]$ and by $(-1)^i$
on $Sym^i(\fh^*)$.

Let $\varpi:St\to \N$ be the projection.
Let $O$ be the $G$-orbit of $e$. We reduce the equivariance
of Borel-Moore homology from $G$ to $G_e$
$$
H_{BM}^{G,i}(\varpi^{-1}O) = H_{BM}^{G_e,i+2(d-2d_e)}(\BBke^2)
,$$
and then to the
maximal torus $C$ in  the
 identity component
$G_e^0$ of $G_e$,
$$
H_{BM}^{G_e,j}({\BBke}^2)=
H_{BM}^{G_e^0,j}(\BBke^2)^\Ga
\aand
H_{BM}^{G_e^0,j}(\BBke^2)=H_{BM}^{C,j}(\BBke^2)^{W(G_e^0)}
.$$
Here $W(G_e^0)$ is the Weyl group of $G_e^0$.
 Also, $H^{2k+1}_{BM}(\BBke^2)=0$ and
$$
H^*(\BBke^2)
\ =\
H^*_C(\BBke^2)\otimes _{H^\bu_C(pt)} H^0(pt)
.$$
The last two isomorphisms follow from the existence of
a $C$-invariant stratification of $\BBke$ where each stratum $X_i$ is a
$C$ equivariant vector
bundle over a space $Z_i$ such that $C$ acts trivially on $Z_i$ and
$H^{2k+1}(Z_i)=0$ \cite{DLP}.

In particular, odd degree components in
$H_{BM}^G(\varpi^{-1}O)$ vanish.
This argument applies to other orbits,
thus we see that
the Cousin spectral sequence for $H_{BM}^G(St)$ corresponding to the
stratification  by the preimages of $G$-orbits
under $\varpi$ degenerates, thus we get a canonical filtration on
 $H_{BM}^G(St)$ whose associated graded pieces are equivariant
Borel-Moore homology spaces of the preimages of $G$-orbits under $\varpi$.

In particular, one of the pieces is $H_{BM}^G(\varpi^{-1}O)$.
The above isomorphisms show that
$H_{BM}(\BBke^2)^\Ga =
[H_{BM}(\BBke)^{\otimes 2}]^\Ga$ is  naturally a
quotient of  $H_{BM}^G(\varpi^{-1}O)$. Thus
$[H_{BM}(\BBke)^{\otimes 2}]^\Ga$ is a subquotient of $H_{BM}^G(St)$.

For $s\in \Ztil\setminus G_e$ the action of $\psi(s)$
on  $H_{BM}(\BBke^2)^\Ga$ is clearly compatible with the action
of $w_0\circ \iota$ on $H_{BM}^G(St)$.
Thus the restriction of this action to
$[H_{BM}^{2i}(\BBke)\otimes H_{BM}^{2j}(\BBke)]^\Ga$
equals $(-1)^{i+j+d-2d_e+d}=(-1)^{i+j}$.

However, since $\BBke$ is compact, Borel-Moore homology coincides
with
homology
$$
H_{BM}^{-k}(\BBke)=H_k(\BBke)=H^k(\BBke)^*
,$$
which yields
\begin{equation}
\label{VV}
\psi(\gamma)|_{[H^{2i}(\BBke)\otimes H^{2j}(\BBke)]^\Ga}
\ =\
\eps^{\ten i+j}(\ga)\cdot  Id \ \ \ \ {\rm for} \ \ \ga
\in \Gatil.
\end{equation}
and thereby finishes the proof. \qed

\sus{
The proof of \eqref{Ups2} for distinguished nilpotents
}

In this subsection we assume that $e$ is distinguished. 
In this case the torus $C$ is trivial; thus we are
dealing with the group
$K^0(Coh(\BBke))$.
The result of \cite{DLP} implies that it is a free abelian group and
 the Chern character map
induces an isomorphism $ch: K^0(Coh(\BBke))\otimes \Ce\iso H_{BM}^\bu(\BBke)$.

\lem {\em 
\lab{Chern chracter}
a) The Chern character map intertwines
Grothendieck-Serre duality $\D$
on
$K^0(Coh(\BBke))$
and the involution $\si$ on $H_{BM}^\bu(\BBke)$ such that
$\si=(-1)^i$ on ${H^{2i}_{BM}}$,
$$
ch\ \ci\ \D
\ =\ \si\ \ci\ ch
.$$

b) The action of
$\B\subset \Baff$ on
$K^0(Coh(\BBke))=K^0(Coh_{\BBke}(\Nt))$ induced by the action
of
$\Baff$ on the category
$D^b[Coh_{\BBke}(\Nt)]$
factors through $W$ and corresponds under $ch$ to the Springer action.

}            \pf
(a) follows from  triviality of the canonical class.
(b) is clear from Theorem \ref{Action of braid group on  base changes}(b) above. \qed

\sss{}
Recall that $\overline{\Upsilon} = \sum_{s=1}^l \
a_{s} g_s^*$  with
$g_s\in A(\Ctil ,\k \cd e)$,\
${\rm ord}(g_s)<\infty$
and $
a_s\in \Qu$\
\eqref{Ups1}.
It is immediate from the definition of
$\Upsilon$ in \cite{Kth2} that the automorphisms $g_s$ of $G$
lie in the outer class of the Chevalley involution.
Thus $g_s$ can be considered as an element in $\Ztil\setminus
G_e$.

Let $\overline{g_s}\in \Gatil\setminus \Ga$
 denote the image of $g_s$ in $\Gatil$,
and set $\upsilon = \sum a_s \overline{g_s}\in \Qu [\Gatil]$.
It is clear from the definitions that the Chern character map
$ch$ intertwines  $\overline{\Upsilon}$ with
$\eta(\upsilon)$, the natural action of $\upsilon$ on $H^*(\BBke)$.

The definition of the modified action $\psi$ and the fact
that $\overline{g_s}\in \Gatil \setminus \Ga$ show that
$\psi(\upsilon)=w_0 \cdot \eta(\upsilon)$, where $w_0$ acts via
the Springer action. By Lemma \ref{Chern chracter}(b)
the endomorphism $\psi(\upsilon)$ is compatible with $T_{w_0}^{-1}
\cdot \overline{\Upsilon}$ under the Chern character map.
Thus, in view of Lemma \ref{Chern chracter}(a),
 we will be done if we check that

\begin{equation}
\label{sumsum}
 \psi(\upsilon)
\ =\ (-1)^i
\ \ \ \ \ \ \ \ \
\tx{on}\
H^{2i}(\BBke)
.
\end{equation}

 Notice that
Proposition \ref{signe} shows that
\eqref{sumsum} holds for all $i$ provided that it holds for $i=d_e$.

This latter fact has almost been checked by Lusztig.
More precisely, \cite{Kth2} implies that
$$
\psi(\upsilon) |_{H^{2d_e}(\BBke)}= \pm 1
.$$
When $e$ is not of type $E_8(b_6)$, then this
is clear from Proposition 5.2 and definition
of $\Upsilon$ in 5.7.
If $e$ is of type $E_8(b_6)$, then this
follows from part IV of the proof
of Proposition 5.2 and definition in 5.7 (all references are to \cite{Kth2}).
Thus it remains to show that
the sign in the last displayed equality equals $(-1)^{d_e}$.

To see this observe that the homomorphism
 $(\BBke\aa{i}\to\BB)^*:\ H^{2d_e}(\BB)\to H^{2d_e}(\BBke)$
is nonzero because the cohomology class of an algebraic cycle is nonzero.
The map $i^*$ is obviously
equivariant with respect to the action of  automorphisms preserving $e$,
and it is well known that this map is $W$-equivariant.

Thus it intertwines $\psi(\upsilon)$ with
 $w_0\cdot \sum a_s g_s$, where $w_0$ acts via
the canonical (Springer) action of $W$ on $H^*(\BB)$ and the action of $g_s$
comes from its action of $\BB$.  So
we will be done if we show that this endomorphism
coincides with $(-1)^i$ on $H^{2i}(\BB)$.

Since $G$ is connected, each $g_s$ acts in fact by the identity map.
Also it is well known that $w_0$ acts by $(-1)^i $ on $H^{2i}$.
So, we are done because we find in
\eqref{Ups1} that $\sum\limits_s a_s=1$.

\subsection{The general case}
Let now $e\in \g$ be an arbitrary nilpotent.
We fix an $sl(2)$ triple $(e,h,f)$ containing $e$ and let
 $\varphi:SL(2)\to G$
be a homomorphism such that the image of $d\varphi$ is spanned by $(e,h,f)$.
We can and will assume that $Im(\varphi)$ commutes with $C$.
There exists an element $\sigma$ in the image of $\varphi$ such that
$Ad(\sigma): e\mapsto -e$.

Recall that $K^C(\BBke)$ is a free module over $K^0(Rep(C))=\Zet[X^*(C)]$
and $K(\BBke)\cong K^C(\BBke)\otimes_{K^0(Rep(C))} K^0(Vect)$.
So, an involution of a free $\Zet[X^*(C)]$-module $M$
which induces identity on the
quotient $M\otimes_{\Zet[X^*(C)]} \Zet$ is itself equal to identity.
Thus it is enough to check that an analogue of \eqref{Ups2} holds in the non-equivariant
K-group.

Furthermore, it suffices to check that this identity holds when the base field $\k$ is of positive
characteristic $p>h$.
In this case the equivalence of \cite{BMR1} provides
 an isomorphism
$K(\BBke) \cong K^0(mod_e^{0,fg}(U)).$

We will identify the two groups by means of this isomorphism.
By the result of \cite[$\S 3$]{BMR2},
the involution  $T_{w_0}\circ \D$
on the left hand side corresponds to the map
$
[M]\mapsto \sigma^*[M^*]
$
on the
right hand side, where for $M\in mod_e^{0,fg}(U)$ we let
$M^*$ denote the dual
$\g$ module (which happens to lie in $mod_{-e}^{0,fg}(U)$).
Thus we are reduced to showing the equality in $K^0(mod_e^{0,fg}(U))$:
\begin{equation}
\label{uu}
[\sigma^* (M^*)]=\Upsilon [M],
\end{equation}
where we set
$\Upsilon [M]=\sum a_s [g_s^*(M)]$,
with $a_s$, $g_s$ being  as in \eqref{Ups1}.

We will actually show an equality stronger than \eqref{uu}. Namely,
consider the
category
$mod^{C,0,fg}_e(U)$ of modules equipped with a compatible
grading by the weights of $C$. We will show that for $M$ in this
category equality \eqref{uu}
holds in $K^0( mod^{C,0,fg}_e(U))$.

We have the Levi subalgebra
$\fl=\fz(C)\subset \g$
such that $e\in \fl$ is distinguished.
By the previous subsection we can assume that the equality is known for
$(e,\fl)$.
We claim that the restriction functor from
$ mod^{C,0,fg}_e(U)$ to
$mod_e^{C.0,fg}(U(\fl))$
induces an injective
map on K-groups. This follows from the well-known fact that an irreducible
module in
$mod_e^{C,0,fg}(U(\fg))$
is uniquely determined by its highest weight component
which is an irreducible object
in $mod_e^{C,0}(U(\fl))$.
[We use an ordering on weights
corresponding to a choice of a parabolic with Levi $L$].

It is clear that this restriction functor
is compatible with the duality functor.
It is also immediate from the definition in
\cite[5.7]{Kth2} that it is compatible
with the involution $\Upsilon$.
Thus \eqref{Ups2} for $e\in \g$ follows from
\eqref{Ups2} for $e\in \fl$.
\qed

\newcommand{\cF}{\mathcal{F}}
\newcommand{\nilrad}{\mathfrak{n}}
\newcommand{\dualnil}{\Langdual\nilcone}
\newcommand{\local}{\Weights^{loc}_{+}}
\newcommand{\complex}{\mathbf C}
\newcommand{\zz}{\mathbf Z}
\newcommand{\closure}{\overline{\orbit}}
\renewcommand{\rank}{\mbox{rank}}

\newcommand{\ph}{\Tilde{\Pi}}
\newcommand{\posroots}{\ro^+}
\newcommand{\ideal}{\mathcal I}
\newcommand{\coal}{{\alpha^{\scriptscriptstyle\vee}}}
\newcommand{\coroot}{{\scriptscriptstyle\vee}}
\newcommand{\ro}{{\Phi}}
\newcommand{\co}{{\Phi^{\scriptscriptstyle\vee}}}
\newcommand{\aw}{{{W_a}}}
\newcommand{\colat}{{Q^{\scriptscriptstyle\vee}}}
\renewcommand{\supp}{\text{ supp}}
\newcommand{\orbit}{\mathcal O}
\newcommand{\levi}{\mathfrak l}

\renewcommand{\wr}{\hat{W}}
\newcommand{\ar}{{\hat{A}}(e)}
\newcommand{\aF}{{\bar{\mathbb{F}}}}
\newcommand{\flag}{{\mathcal B}}

\newcommand{\dualp}{\lambda'}

\newcommand{\rep}{\chi^{\alpha, \beta}}

\se{A result on component groups,\ by Eric Sommers}
\label{Eric}

Here, $G$ is a reductive algebraic group over the algebraically closed field $\Bbbk$ and $\mathfrak{g}$ its Lie algebra.
As in Section 5.2.2, we are given a homomorphism
$$\varphi: SL_2(\Bbbk) \to G$$  
and the characteristic of $\Bbbk$ is at least $3h-3$.

Let $$s = \varphi \left( \begin{smallmatrix}  -1 & 0  \\ 0 & -1 \end{smallmatrix} \right)$$ 
and $$e =  d\varphi \left( \begin{smallmatrix}  0 & 1  \\ 0 & 0 \end{smallmatrix} \right).$$  It is clear that $s \in G_e$  and $e \in \mathfrak{g}_s$.

Recall that $\phi: {\mathbb G}_m \to G$ is defined as $\phi(t) = \varphi \left( \begin{smallmatrix}  t & 0  \\ 0 & t^{-1} \end{smallmatrix} \right)$.

\pro {\em
If $G$ is semisimple and adjoint, then $s$ belongs to the identity component of $G_e$.
}

\rem After this appendix was written and made available in a preprint form, we learned that the result was also proved
by A.~Premet \cite[Lemma 2.1]{Premet1}

\pf
Let $x \in G$ be an arbitrary semisimple element commuting with $e$.   The conjugacy class of the image of $x$ in $G_e/(G_e)^0$ is determined by the
$G$-orbit of the pair $(e, \mathfrak{l''})$ where ${\mathfrak{l''}}$ is any Levi subalgebra of ${\mathfrak{g}}_x$ such that $e \in  \mathfrak{l''}$ is distinguished.
More precisely, two semisimple elements commuting with $e$ have conjugate image in
$G_e/(G_e)^0$ if and only if the corresponding pairs as above are $G$-conjugate.
This result is true in any good characteristic by \cite{ms}, \cite{premet}.
In the case where $x=1$, the $G$-orbit of such pairs includes $(e, \mathfrak{l})$ where $\mathfrak{l}$ is a Levi subalgebra of $\mathfrak{g}$ such that
 $e \in  \mathfrak{l}$ is distinguished.
Hence, an arbitrary $x$ as above lies in the identity component of $G_e$  if and only if $e$ is distinguished in  ${\mathfrak{l''}} \subset  {\mathfrak{g}}_x$ where
$\mathfrak{l''}$ is a Levi subalgebra of $\mathfrak{g}$ (and not only  of ${\mathfrak{g}}_x$).

Now as in Section 5.2.2, let $C$ be a maximal torus in the centralizer of the image of $\varphi$ in $G$.  Then
with the assumption on the characteristic of $\Bbbk$,
$C$ is a maximal torus of $G_e$ and thus $e$ is distinguished in
the Levi subalgebra ${\mathfrak{l}} = Z_{\mathfrak{g}}(C)$ of ${\mathfrak{g}}$ (see \cite{carter}).   We then also have that
the orbit of $e$ in $\fl$ is an even nilpotent orbit.
In other words, if we pick a maximal torus of $L =Z_{G}(C)$ containing the image of $\phi$, then each root of $L$ paired with the co-character $\phi$ is an even integer. Thus $s=\phi(-1)$ acts trivially on $\mathfrak{l}$,
and hence ${\mathfrak{l}} \subset {\mathfrak{l'}}:={\mathfrak{g}}_s$.

On the other hand, since $C \subset G_s$
we have that $Z_{\mathfrak{l'}}(C)$ is a Levi subalgebra of $\mathfrak{l'}$.
But by the previous paragraph,
${\mathfrak{l}} = Z_{\mathfrak{g}}(C) \subset \mathfrak{l}'$, so ${\mathfrak{l}} = Z_{\mathfrak{l'}}(C)$. Therefore
${\mathfrak{l}}$ is a Levi subalgebra of both ${\mathfrak{g}}_s$ and ${\mathfrak{g}}$, and we can conclude by the first paragraph
that $s$ lies in the identity component of $G_e$.
\qed

\rem
A similar result holds in all good characteristics for $s = \phi(-1)$, where $\phi$ is an associated co-character of a nilpotent element $e$.   In this case, $C$ is defined to be the maximal torus in the simultaneous centralizer in $G$ of $e$ and the image of $\phi$.   Then $e$ is distinguished in $Z_{\mathfrak{g}}(C)$ as before
and by \cite{premet} or \cite{jantzen}
$\phi$ corresponds to a weighted Dynkin diagram arising in characteristic zero for a distinguished element for the corresponding Levi subalgebra.   Therefore,
it remains true that $s$ acts trivially on $Z_{\mathfrak{g}}(C)$ and the proof goes through. 

\cor {\em 
For reductive $G$, $sz \in C$ for some $z \in Z(G)$, where $C$ is as above.
}         

   \pf
As $G/Z(G)$ is semisimple and adjoint, it amounts to showing that 
$s \in C$ when $G$ is semisimple and adjoint.  Assume the latter.
We know that $s$ centralizes $C$ by definition.
Then since $s$ is in the identity component of $G_e$ by the proposition, we know that $s$ belongs to the centralizer of $C$ in the identity component of $G_e$.
That centralizer is equal to $C$ itself, being the centralizer of a maximal torus in a connected group.  Hence $s \in C$.

I thank David Vogan, Jr.\ for helpful conversations.

\bibliographystyle{pnaplain}

\end{document}